\documentclass[12pt,reqno,basicdelimiters]{NumPDEsArticle}

\usepackage[utf8]{inputenc}
\usepackage[T1]{fontenc}
\usepackage[english]{babel}

\usepackage{NumPDEsMacros}

\usepackage{csquotes}
\usepackage{enumitem}
\usepackage{amsmath,amssymb}
\usepackage{bbm}
\usepackage{microtype}
\usepackage{appendix}
\usepackage{xcolor}
\usepackage{bm}
\usepackage{stackengine}
\usepackage{subfig}
\usepackage{float}
\usepackage[backend=biber]{biblatex}

\makeatletter
\renewcommand{\enddoc@text}{}
\makeatother

\newcommand{\Linf}{L^{\infty}(\Omega)}

\newcommand{\gen}[1]{\mathrm{gen}(#1)}

\newcommand{\supp}[1]{\mathrm{supp}(#1)}

\newcommand{\JJJ}{\mathfrak{J}}

\newcommand{\PPP}{\mathfrak{P}}

\newcommand{\MMM}{\mathfrak{M}}

\newcommand{\SSS}{\mathfrak{S}}

\newcommand{\RRR}{\mathfrak{R}}

\newcommand{\QQQ}{\mathfrak{Q}}

\allowdisplaybreaks
\title{Convergence of an adaptive multilevel stochastic Galerkin FEM driven by a residual-based error estimator}

\author{Alex Bespalov}
\address{School of Mathematics, University of Birmingham, Edgbaston, Birmingham B15 2TT, UK}
\email{a.bespalov@bham.ac.uk}

\author{Maximilian Brunner}
\address{TU Wien, Institute for Analysis and Scientific Computing, Wiedner Hauptstraße 8-10, 1040 Wien}
\email{maximilian.brunner@asc.tuwien.ac.at}

\author{Alexander Freiszlinger}
\email{alexander.freiszlinger@asc.tuwien.ac.at}

\author{Dirk Praetorius}
\email{dirk.praetorius@asc.tuwien.ac.at}

\thanks{This research was funded by the Austrian Science Fund (FWF) projects
\href{https://www.fwf.ac.at/en/research-radar/10.55776/F65}{10.55776/F65} (SFB
F65 ``Taming complexity in PDE systems''),
\href{https://www.fwf.ac.at/en/research-radar/10.55776/I6802}{10.55776/I6802}
(international project I6802 ``Functional error estimates for PDEs on unbounded
domains''), and
\href{https://www.fwf.ac.at/en/research-radar/10.55776/P33216}{10.55776/P33216}
(standalone project P33216 ``Computational nonlinear PDEs''). \\
The generative AI tool Claude was used to check spelling and grammar in the manuscript, but was not involved in the development of the mathematical content or the software used in the numerical simulations. The authors take full responsibility for the content of the publication.}

\begin{document}

\maketitle

\begin{abstract}
	In this paper we analyze a residual-based \textsl{a~posteriori} error estimation strategy and adaptive multilevel algorithms for a class of parametric elliptic PDEs, where the diffusion coefficient depends on countably many uncertainty parameters. 
  Under a saturation assumption on the parametric approximation space, we propose a computable residual-based error estimator.
   We also introduce the notion of $K$-compatibility of finite element meshes, which ensures that only the triangulations with comparable mesh sizes are used, and incorporate this concept into the adaptive algorithm. 
   We prove that the estimator is reliable and satisfies $K$-dependent stability and reduction properties.
   These results underpin our convergence analysis of the adaptive multilevel algorithm: we show that the algorithm drives both the \textsl{a~posteriori} error estimator and the true approximation error to zero. 
   Numerical experiments illustrate and confirm the theoretical findings.
\end{abstract}


\section{Introduction}

\subsection{Model problem and state of the art}
Let $\Omega \subset \R^d$, $d \geq 2$, be a bounded domain with Lipschitz boundary $\partial \Omega$ and $\Lambda \coloneqq \prod_{m \in \N} [-1,1]$ be the infinite-dimensional hypercube. 
We consider the following parametric model problem: 
Given a volume force $f\colon \Omega \to \R$, find $u^\star \colon \Omega \times \Lambda \to \R$ such that 
\begin{equation} \label{eq:modelproblem}
 \begin{split}
 -\div(a \nabla u^\star) 
 &= f \quad \text{in } \Omega \times \Lambda \qquad 
 \text{subject to} \qquad 
 u^\star = 0 \quad \text{on } \partial \Omega \times \Lambda,
 \end{split}
\end{equation}
where the diffusion coefficient $a \colon \Omega \times \Lambda \to \R$ is affine with respect to $\vec{y} = (y_m)_{m \in \N} \in \Lambda$, i.e.,
\begin{equation} \label{eq:affine}
 \begin{split}
  a(x,\vec{y}) 
  &= a_0(x) + \sum_{m \in \N} a_m(x) y_m.
 \end{split}
\end{equation}

The boundary value problem~\eqref{eq:modelproblem} and other PDE-based models with stochastic inputs arise in various fields of science and engineering. 
Efficient numerical algorithms for solving these PDEs are paramount in the context of \emph{uncertainty quantification}. 
The construction of such algorithms has attracted considerable attention over the last two decades in connection with, e.g., reduced basis schemes (e.g., \cite{Haasdonk2013}), collocation methods (e.g., \cite{Babuska2007}) and (quasi-)Monte Carlo sampling (e.g., \cite{Charrier2013}).

Stochastic Galerkin finite element methods (SGFEMs) are a useful alternative to the above mentioned sampling methods; we refer to \cite{Xiu2009,Baeck2011,Gunzburger2014} for details of the differences between SGFEM and other methods. 
One of the main advantages of \textit{intrusive} Galerkin methods over \textit{non-intrusive} sampling or collocation methods is the significantly lower number of function evaluations and solver steps; see, e.g., \cite{Herzog2008} and~\cite[Section~6]{Xiu2009}.

SGFEM approximations combine \emph{generalized polynomial chaos} (gPC) expansions in the stochastic variables/parameters with finite element discretizations of the expansion coefficients. The approximations (and the algorithms that generate them) are termed \textit{single-level} if the same finite element mesh is used for all expansion coefficients, and \textit{multilevel} otherwise. 
Single-level algorithms have emerged in \cite{Eigel2014a}, while multilevel approximations have been considered in \cite{Cohen2010,Gittelson2013}, among other works. 
Flexibility is the main advantage of multilevel methods, especially in the context of adaptivity, since each spatial coefficient in the gPC expansion is allowed to reside in a finite element space associated with its own appropriately chosen mesh. 
Thus, multilevel algorithms can successfully tackle problems where some gPC expansion coefficients require a highly refined mesh to approximate the solution accurately, while a coarser mesh suffices for other gPC coefficients.

\textsl{A~posteriori} error estimation strategies that steer adaptive multilevel refinement and the practical realization of the multilevel SGFEM have been considered, e.g., in~\cite{Eigel2014,Crowder2019,Bespalov2021a,Bespalov2021}. 
In the latter work, under a saturation assumption on the \emph{full} approximation space, convergence and rate-optimality of an adaptive multilevel algorithm are proven; in that work, the adaptive algorithm employs hierarchical \emph{parametric} error indicators and two-level \emph{spatial} error indicators together with \textit{combined} D\"orfler marking. 
Very recently, in \cite{Bachmayr2024}, convergence was proven for an adaptive algorithm that utilized multilevel expansions of the diffusion coefficient and residual error indicators based on finite element frames. 
However, convergence of an adaptive multilevel algorithm employing residual-based error indicators as proposed in \cite{Eigel2014} has so far remained an open problem.

\subsection{Contributions}

In this paper, we follow \cite{Eigel2014} and employ residual-based error indicators to steer mesh refinement and parametric enrichment in an adaptive multilevel algorithm. 
At present, this type of \textsl{a~posteriori} estimation for the \emph{parametric} error component remains largely theoretical and is not directly accessible in computations.
Under a saturation assumption on the \emph{parametric} approximation space, we overcome this limitation.

The main difficulty in the analysis is that the spatial error indicator associated with a coefficient in the gPC expansion also depends on other gPC coefficients, which, in general, do not reside in the same finite element space. 
While this is not an issue for single-level approximations~\cite{Eigel2014a,Bespalov2019a}, the analysis of multilevel methods requires greater care. 
Indeed, to ensure stability of the spatial estimator, one needs to control the ratio of the mesh sizes of triangulations associated with neighboring gPC coefficients. 
In \cite{Eigel2014}, this issue is addressed by introducing projection operators between different finite element spaces. 
However, this leads to  additional projection errors exhibiting unstable behavior. 
We propose an alternative way to guarantee stability of the spatial estimator by introducing the notion of $K$-compatibility for a set of triangulations, which ensures that their mesh sizes are compatible. 
We then follow~\cite{Eigel2014a} and prove that the residual-based estimator is \emph{reliable} and satisfies $K$-dependent versions of the stability and reduction properties (see, e.g.,~\cite{Cascon2008,Carstensen2014}). 
Furthermore, for a suitable sequence $(K_\ell)_{\ell \in \N_0}$ of compatibility degrees, we develop a novel adaptive algorithm with built-in $K_\ell$-compatibility of all meshes employed at the $\ell$-th iteration of the algorithm.
In particular, we prove that the underlying multilevel refinement strategy steered by \emph{combined} D\"orfler marking maintains $K_\ell$-compatibility of meshes throughout the adaptive loop.
Finally, we show that the developed algorithm drives both the \textsl{a~posteriori} error estimator and the true approximation error to zero.

\subsection{Outline}

The paper is organized as follows. 
Section~\ref{sec:prelim} collects all necessary preliminaries and notation on function spaces, triangulations, multilevel structures, and SGFEM.
In Section~\ref{sec:apost}, we propose a computable and reliable residual-based \textsl{a~posteriori} error estimator.
We also propose an adaptive algorithm and, subsequently, formulate Theorem~\ref{thm:main}, the main result of this work. It states that the overall error and, under stronger assumptions on $(K_\ell)_{\ell \in \N_0}$, also the error estimator converge to zero.
Sections~\ref{sec:proof}--\ref{sec:stabred} are devoted to the proof of Theorem~\ref{thm:main}. 
In Section~\ref{sec:algo}, we develop multilevel refinement strategies which maintain $K_\ell$-compatibility throughout the adaptive loop. 
Section~\ref{sec:numerics} concludes the work with numerical experiments.

\section{Preliminaries and Notation} \label{sec:prelim}

\subsection{General notation}

Throughout this paper, let $\Omega$ be a bounded domain in $\R^d$, $d \geq 2$, with Lipschitz boundary $\partial \Omega$. By $\abs{\, \cdot \,}$ we denote, without any ambiguity, the absolute value of a scalar, the $d$-dimensional Lebesgue volume of a measurable set, the $(d-1)$-dimensional Hausdorff measure of a surface (piece), or the cardinality of a finite set. 
We write $\int_{U} \cdot \d x$ for integration over a Lebesgue measurable set $U \subseteq \Omega$. 
For a Lipschitz domain $U \subseteq \Omega$, we denote by $\bm{n}_{\partial U}$ the outward pointing unit normal vector on $\partial U$.

Throughout, we denote by $\Lambda \coloneqq \prod_{m \in \N} [-1,1]$ the infinite-dimensional hypercube and by $\BB(\Lambda)$ the Borel $\sigma$-algebra on $\Lambda$. 
For $m \in \N$, let $\pi_m$ be a symmetric Borel probability measure on $[-1,1]$ and set $\pi \coloneqq \prod_{m \in \N} \pi_m$.
Integration over $[-1,1]$ with respect to $\pi_m$ is denoted by $\int_{[-1,1]} \cdot \d \pi_m(y_m)$ and integration over $\Lambda$ with respect to $\pi$ is denoted by $\int_{\Lambda} \cdot \d \pi(\vec{y})$.

For a Banach space $\XX$, we denote by $\norm{\, \cdot\, }_\XX$ the norm on $\XX$ and by $\XX^*$ the topological dual space of $\XX$. The duality pairing is denoted by $\dual{\cdot}{\cdot}_{\XX^*}$.

For a countable set $I$, we introduce the standard space $\ell^2(I)$ of square-summable sequences with the associated norm
\begin{equation*}
 \norm{(x_i)_{i \in I}}_{\ell^2(I)} 
 \coloneqq \Bigl(\sum_{i \in I} x_i^2\Bigr)^{1/2} 
 \quad \text{for all } (x_i)_{i \in I} \in \ell^2(I).
\end{equation*}

Discrete quantities are indicated by a subscript, where related objects are indexed by the same subscript, e.g., $\TT_\ell$ denotes a triangulation of $\Omega$ with corresponding finite element functions $\varphi_\ell$ etc. 
Similarly, we use the same subscript for quantities related to a multi-index, e.g., $w_{\ell,\nu}$ is a discrete coefficient corresponding to a mesh $\TT_{\ell,\nu}$ and so forth. 
To abbreviate notation in proofs, we write $A \lesssim B$ if and only if there holds $A \leq cB$ with some constant $c > 0$ that does not depend on the iteration counter $\ell \in \N_0$. 
We write $A \simeq B$ if $A \lesssim B$ and $B \lesssim A$. 
Finally, we introduce the Kronecker delta
\begin{equation*}
 \delta_{AB} \coloneqq 
 \begin{cases}
  1 & \text{if } A = B, \\
  0 & \text{if } A \neq B,
 \end{cases} 
 \quad \text{for any quantities } A, B.
\end{equation*}

\subsection{Function spaces}

For a Lipschitz domain $U \subseteq \Omega$, we denote by $L^2(U)$ the usual space of square-integrable functions on $U$ with the corresponding norm $\norm{\, \cdot \,}_{L^2(U)}$ and by $L^\infty(U)$ the space of all measurable functions $v$ on $U$ such that $\norm{v}_{L^\infty(U)} \coloneqq \mathrm{ess} \sup\limits_{x \in U} \, \abs{v(x)} < +\infty$.
If it is clear from the context, we may write $\norm{\cdot}_{U}$ instead of $\norm{\cdot}_{L^2(U)}$. 
For vector-valued spaces $[L^2(U)]^n$ or $[L^\infty(U)]^n$ with $n \in \N$, we will omit the superscript in unambiguous cases. 
We introduce the standard first-order Sobolev space
\begin{equation*}
 H^1(\Omega) 
 \coloneqq \set{\varphi \in L^2(\Omega) \given \nabla \varphi \in L^2(\Omega)}
 \quad \text{with the norm} \quad \norm{\varphi}_{H^1(\Omega)} \coloneqq (\norm{\varphi}_{\Omega}^2 + \norm{\nabla \varphi}_{\Omega}^2)^{1/2}.
\end{equation*}
Furthermore, we define the space of $H^1(\Omega)$-functions with vanishing trace on $\partial \Omega$ as
\begin{equation*}
 H_0^1(\Omega) 
 \coloneqq \set{\varphi \in H^1(\Omega) \given \varphi|_{\partial \Omega} = 0}.
\end{equation*}
Additionally, we introduce the space
\begin{equation*}
 W^{1,\infty}(\Omega) 
 \coloneqq \set{\varphi \in L^\infty(\Omega) \given \nabla \varphi \in L^\infty(\Omega)}
\end{equation*}
with the associated norm $\norm{\varphi}_{W^{1,\infty}(\Omega)} \coloneqq \max\set{\norm{\varphi}_{L^\infty(\Omega)}, \norm{\nabla \varphi}_{L^\infty(\Omega)}}$.

\subsection{Bochner spaces and polynomial chaos} \label{subsec:gpc}

Given a Hilbert space $\HH$ with norm $\norm{\, \cdot \,}_\HH$, we denote by $L^2_{\pi}(\Lambda;\HH)$ the Bochner space of $\HH$-valued functions on $\Lambda$ with respect to the measure $\pi$, i.e., 
\begin{equation*}
 L^2_{\pi}(\Lambda;\HH) 
 \coloneqq \Big\{ v \colon \Lambda \to \HH \ | \ \norm{v}_{L^2_\pi(\Lambda;\HH)} 
 \coloneqq \Bigl(\int_{\Lambda} \norm{v(\vec{y})}_{\HH}^2 \d \pi(\vec{y}) \Bigr)^{1/2} < +\infty\Big\}.
\end{equation*}%
If $\HH = \R$, we write $L^2_{\pi}(\Lambda)$ instead of $L^2_{\pi}(\Lambda;\R)$. 
An orthonormal basis of $L^2_\pi(\Lambda)$ is given by the following construction$\colon$ 
For $m \in \N$, let $(P_m^n)_{n \in \N_0}$ be the sequence of normalized orthogonal polynomials on $[-1,1]$ with respect to $\pi_m$, i.e., the degree of $P_m^n$ is $n$ and 
\begin{equation} \label{eq:ortho}
 \int_{[-1,1]} P_m^n(y_m) P_m^k(y_m) \d \pi_m(y_m) 
 = \delta_{nk}.
\end{equation}
Setting $P_m^{-1} \coloneqq 0$, there exist unique positive constants $(\beta_m^n)_{n \in \N_0}$ such that there holds the three-term recurrence
\begin{equation} \label{eq:recurrence}
 \beta_m^nP_m^{n+1}(y_m) 
 = y_m P_m^{n}(y_m) - \beta_m^{n-1} P_m^{n-1}(y_m) 
 \quad \text{for all } y_m \in [-1,1] \text{ and } n \in \N_0.
\end{equation}
Let $\N_0^{\N} \coloneqq \set{\nu = (\nu_m)_{m \in \N} \given \nu_m \in \N_0 \text{ for all } m \in \N}$ be the set of all multi-indices and define $\mathrm{supp}(\nu) \coloneqq \set{m \in \N \given \nu_m \neq 0}$ as well as 
\begin{equation*}
 |\nu| 
 \coloneqq 
 \sum_{m \in \supp{\nu}} \nu_m \in \N_0 \cup \set{+\infty}
 \quad \text{for all } \nu \in \N_0^\N.
\end{equation*}
We define the set of \emph{multi-indices} with finite support as
\begin{equation*}
 \JJJ 
 \coloneqq 
 \set{\nu \in \N_0^{\N} \given \abs{\mathrm{supp}(\nu)} < +\infty} 
 = 
 \set{\nu \in \N_0^{\N} \given \abs{\nu} < +\infty}.
\end{equation*}
By defining $\eps_m \coloneqq (\delta_{mk})_{k \in \N} \in \JJJ$ for all $m \in \N$, every multi-index $\nu \in \mathfrak{J}$ can be written as
\begin{equation*}
 \nu 
 = \sum_{m \in \N} \nu_m \eps_m 
 = \sum_{m \in \supp{\nu}} \nu_m \eps_m.
\end{equation*}
Defining the multivariate polynomial
\begin{equation} \label{eq:multivariate}
 P_{\nu}(\vec{y})
 \coloneqq
 \prod_{m \in \N} P_m^{\nu_m}(y_m) 
 =
 \prod_{m \in \mathrm{supp}(\nu)} P_m^{\nu_m}(y_m)
 \quad \text{for all } \nu \in \mathfrak{J} \text{ and all } \vec{y} = (y_m)_{m \in \N} \in \Lambda,
\end{equation} 
the set $\set{P_{\nu} \given \nu \in \mathfrak{J}}$ is an orthonormal basis of $L^2_{\pi}(\Lambda)$, see, e.g., \cite{Ernst2012}, \cite{Eigel2014}. 
Finally, we note that it is well known that, for separable Hilbert spaces $\HH$, the Bochner space $L^2_{\pi}(\Lambda;\HH)$ is isometrically isomorphic to the tensor product space $\HH \otimes L^2_{\pi}(\Lambda)$, i.e., for every $v \in L^2_{\pi}(\Lambda;\HH)$, there exist unique coefficients $(v_\nu)_{\nu \in \JJJ} \subset \HH$ such that
\begin{equation*}
 v(\vec{y}) 
 = \sum_{\nu \in \JJJ} v_\nu P_\nu(\vec{y}) 
 \quad \text{for all } \vec{y} \in \Lambda.
\end{equation*}

\subsection{Triangulations} \label{subsec:meshes}

For $k \leq d$, a $k$-dimensional simplex $T \subset \R^d$ is a compact set with positive $k$-dimensional measure $\abs{T} > 0$ and $T = \mathrm{conv}(\NN_T)$, where $\NN_T = \set{z_0, \ldots, z_k}$ are the vertices of $T$. 
A conforming mesh $\TT$ on $\Omega$ is a partition of $\overline{\Omega}$ into $d$-dimensional simplices (called elements) such that the intersection $T \cap T'$ of two different elements $T,T' \in \TT$ is either empty, a common vertex, or the joint $k$-dimensional simplex $T \cap T' = \mathrm{conv}(\NN_T \cap \NN_{T'})$ with $1 \leq k \leq d -1$.

Throughout this paper, we assume that all considered meshes are obtained by a finite number of newest-vertex bisection (NVB) refinements of some initial conforming mesh $\TT_0$. 
In the special case $d=2$, NVB is based on choosing a reference edge $E_T \subset \partial T$ for every $T \in \TT$, which, if marked for refinement, is bisected by connecting the midpoint of $E_T$ with the opposite vertex in $T$.  
For the newly created triangles $T_1,T_2 \subseteq T$, the reference edges $E_{T_1},E_{T_2}$ are chosen as the edges opposite to the new vertex.
We assume that $\TT_0$ is \emph{admissible}, i.e., for any two triangles $T,T' \in \TT_0$ that share an edge $E = T \cap T'$, there holds either $E_T = T \cap T' = E_{T'}$ or $E_T \neq T \cap T' \neq E_{T'}$. 
We refer to~\cite{Stevenson2007} for NVB for $d \geq 2$ with admissible $\TT_0$, as well as to~\cite{Karkulik2012} ($d = 2$) and~\cite{Diening2025} ($d \geq 2$) for NVB with non-admissible $\TT_0$. 

It is well known that NVB preserves uniform $\kappa$-shape regularity of the mesh, i.e., there exists $\kappa > 0$ depending only on $\TT_0$ such that each mesh $\TT$ obtained by a finite number of NVB refinements of $\TT_0$ satisfies 
\begin{equation*}
 \sup_{T \in \TT} \frac{\mathrm{diam}(T)}{\abs{T}^{1/d}} 
 \leq \kappa.
\end{equation*}
We denote by $\NN$, $\EE$, and $\FF$ the sets of ($0$-dimensional) vertices, ($1$-dimensional) edges, and ($(d-1)$-dimensional) facets of the mesh $\TT$, respectively. 
Moreover, given $T \in \TT$, we define $\EE_T$ and $\FF_T$ as the edges and facets contained in $T$, respectively.
For given meshes $\TT_1$ and $\TT_2$, we write $\TT_1 \succeq \TT_2$ or $\TT_2 \preceq \TT_1$, if for all triangles $T \in \TT_1$ and $T' \in \TT_2$ with $\abs{T \cap T'} > 0$, there holds $T' \subseteq T$. We then say that $\TT_2$ is a refinement of $\TT_1$. 
We define the piecewise constant mesh-size function $h \in L^\infty(\Omega)$ by $h|_T \coloneqq h_T \coloneqq \abs{T}^{1/d}$. The quantities $h_E \coloneqq \abs{E}$ and $h_F \coloneqq \abs{F}^{1/(d-1)}$ for edges $E \in \EE$ and facets $F \in \FF$, respectively, are defined analogously. Mesh-size functions associated with different meshes will be distinguished by the corresponding subscript, e.g., $h_{\ell,\nu}$ is the mesh-size function associated with the mesh $\TT_{\ell,\nu}$.

Given a mesh $\TT$ and a triangle $T \in \TT$, there exists a unique triangle $T_0 \in \TT_0$ such that $ T \subseteq T_0$. We define the \emph{level} or \emph{generation} of $T$ as
\begin{equation*}
 \gen{T} 
 \coloneqq \log_2 \abs{T_0} - \log_2 \abs{T}
 \in \N_0
\end{equation*}%
and note that there exists a constant $\const{C}{gen} > 0$ depending only on $\TT_0$ and $d$ such that 
\begin{equation} \label{eq:levelequiv}
 \const{C}{gen}^{-1} \, 2^{-\gen{T}/d} \leq h_T 
 \leq \const{C}{gen} \, 2^{-\gen{T}/d} 
 \quad \text{for all triangulations } \TT \text{ and } T \in \TT;
\end{equation}
see, e.g.,~\cite{Karkulik2012}. Although the estimates~\eqref{eq:levelequiv} are proven in~\cite{Karkulik2012} only for $d=2$, we remark that the proof transfers verbatim to the case $d \geq 2$.

We introduce the notion of \emph{$K$-compatibility} for triangulations $\TT_1$ and $\TT_2$ as follows: We say that $\TT_1$ and $\TT_2$ are \emph{$K$-compatible} for some $K \in \N$, if
\begin{equation} \label{eq:defLcomp}
 \abs{\gen{T_1} - \gen{T_2}} 
 \leq K 
 \quad \text{for all } T_1 \in \TT_1 \text{ and } T_2 \in \TT_2 \text{ with } \abs{T_1 \cap T_2} > 0.
\end{equation}
In this case, we will call $K$ the \emph{compatibility degree} of the meshes $\TT_1$ and $\TT_2$.
Because of~\eqref{eq:levelequiv}, $K$-compatibility is equivalent to
\begin{equation} \label{eq:defLcomp2}
 \const{C}{gen}^{-2} 2^{-K/d} h_{T_1} 
 \leq h_{T_2} 
 \leq \const{C}{gen}^2 2^{K/d} h_{T_1}
 \quad \text{for all } T_1 \in \TT_1 \text{ and } T_2 \in \TT_2 \text{ with } \abs{T_1 \cap T_2} > 0.
\end{equation}

Finally, given two triangulations $\TT_1$ and $\TT_2$, let 
\begin{equation*}
 \TT_1 \oplus \TT_2 
 \coloneqq \set{T_1 \cap T_2 \given T_1 \in \TT_1 \text{ and } T_2 \in \TT_2 \text{ with } \abs{T_1 \cap T_2} > 0}
\end{equation*}%
denote the \emph{overlay} of $\TT_1$ and $\TT_2$ (i.e., the coarsest common refinement). 
Due to the set operation $\cap$ being associative, $\oplus$ is also associative. 
Thus, we define the overlay of a finite number of triangulations $\TT_1, \ldots, \TT_n$ as
\begin{equation*}
 \bigoplus_{k = 1}^n \TT_k 
 \coloneqq \Big\{\bigcap_{k = 1}^n T_k \ | \ T_k \in \TT_k \text{ for all } k = 1, \ldots, n \text{ with } \Big| \bigcap_{k = 1}^n T_k \Big| > 0 \Big\}.
\end{equation*}

\subsection{Discrete spaces and multilevel structures} \label{subsec:discspaces}

For a given set $U \subseteq \R^d$ and $p \in \N$, we denote by $\P^p(U)$ the set of all polynomials of degree at most $p$ on $U$. 
For a given mesh $\TT_\ell$ on $\Omega$, we define the conforming Lagrange finite element space as
\begin{equation*}
 \SS^p_0(\TT_\ell) 
 \coloneqq \set{v \in H^1_0(\Omega) \given v|_T \in \P^p(T) \text{ for all } T \in \TT_\ell}.
\end{equation*}
For polynomial approximations in $\vec{y} \in \Lambda$, we introduce an \emph{index set} $\PPP_\ell$ as a finite subset of $\JJJ$.
The components of any multi-index $\nu \in \PPP_\ell$ represent polynomial degrees in the multivariate polynomial $P_\nu(\vec{y})$ (see~\eqref{eq:multivariate}).
We define
\begin{equation*}
 \supp{\PPP_\ell} 
 \coloneqq \bigcup_{\nu \in \PPP_\ell} \mathrm{supp}(\nu) 
 = 
 \set{m \in \N \given \exists \, \nu \in \PPP_\ell: \nu_m > 0}
\end{equation*} 
as the set of \emph{active components} of all multi-indices in $\PPP_\ell$ and $M_\ell \coloneqq \abs{\mathrm{supp}(\PPP_\ell)} \in \N_0$ as the number of these active components.
Recall that $\eps_m = (\delta_{km})_{k \in \N}$ for $m \in \N$. We denote by $\vec{0} \coloneqq (0, \, 0 , \, \ldots)$ the zero index and consider an initial index set $\PPP_0$ that always contains $\vec{0}$ and may also contain multi-indices $\varepsilon_m$, $m = 1,2,\ldots, M_0$, for some $M_0 \in \N$ (e.g., $\PPP_0 \coloneqq \set{\vec{0}}$ or $\PPP_0 \coloneqq \set{\vec{0}, \eps_1, \eps_2, \ldots , \eps_{10}}$).
Next, we define the \emph{detail index set} as
\begin{equation} \label{eq:detaildef}
 \begin{split}
  \mathfrak{Q}_\ell 
  \coloneqq \big\{ \nu \in \JJJ \setminus \PPP_\ell \; \vert \; &\exists \, m \in \set{1,2,\ldots,M_\ell} \text{ such that }  \\
  &\nu + \eps_m \in \PPP_\ell \text{ or } \nu - \eps_m \in \PPP_\ell \big\} \cup \set{\eps_{M_\ell + 1}}.
 \end{split}
\end{equation}
Note that all the multi-indices in $\QQQ_\ell$ have an increased or decreased component compared to the multi-indices in $\PPP_\ell$, and the new component is activated by adding only the corresponding unit multi-index.
Note that all the multi-indices in $\QQQ_\ell$ have one more active component than the multi-indices in $\PPP_\ell$, and the new component is activated by adding only the corresponding unit multi-index. 
Throughout the paper, in order to refine polynomial approximations on $\Lambda$ we enrich the index set $\PPP_\ell$ by adding to it a set $\MMM_\ell \subseteq \QQQ_\ell$ of marked multi-indices.
For $\nu \in \JJJ$, we denote by $\PPP_\ell[\nu]$ the set of all \emph{neighbors} of $\nu$ in $\PPP_\ell$, i.e.,
\begin{equation*}
 \PPP_\ell[\nu] 
 \coloneqq \set{\mu \in \PPP_\ell \given \exists \, m \in \N \colon \mu = \nu + \eps_m \text{ or } \mu = \nu - \eps_m}.
\end{equation*}

Finally, given a family of meshes $(\TT_{\ell,\nu})_{\nu \in \PPP_\ell}$, we call the pair $[\PPP_\ell,(\TT_{\ell,\nu})_{\nu \in \PPP_\ell}]$ a \emph{multilevel structure} and define the initial multilevel structure as $[\set{\vec{0}},(\TT_0)]$. 
We say that the multilevel structure $[\PPP_\ell,(\TT_{\ell,\nu})_{\nu \in \PPP_\ell}]$ is $K$-compatible, if $\TT_{\ell,\nu}$ and $\TT_{\ell,\mu}$ are $K$-compatible for all $\nu \in \PPP_\ell$ and $\mu \in \PPP_\ell[\nu]$.

\subsection{Multilevel stochastic Galerkin FEM} \label{subsec:sgfem}

In line with~\cite{Schwab2011}, we make the following assumptions on the expansion coefficients in the representation~\eqref{eq:affine} of the diffusion coefficient $a\colon$ $a_m \in W^{1,\infty}(\Omega)$ for all $m \in \N_0$ and there exist positive constants $0 < \const{a}{min} \leq \const{a}{max} < +\infty$ such that
\begin{equation} \label{eq:a0bounds}
 \const{a}{min} 
 \leq 
 a_0(x)
 \leq 
 \const{a}{max} 
 \qquad \text{for almost all } x \in \Omega.
\end{equation} 
Moreover, we require that
\begin{equation} \label{eq:asumbound}
 \gamma 
 \coloneqq 
 a_{\mathrm{min}}^{-1} \, \Big \lVert \sum_{m \in \N} \abs{a_m}\Big \rVert_{L^\infty(\Omega)} 
 < 1, \qquad 
 \sum_{m \in \N} \norm{a_m}_{L^\infty(\Omega)} 
 < +\infty, 
\end{equation}%
and that
\begin{equation} \label{eq:ampliesti}
 \sum_{m \in \N_0} \norm{\nabla a_m}_{L^{\infty}(\Omega)} 
 < 
 + \infty.
\end{equation}

\begin{remark}
 Assumption~\eqref{eq:ampliesti} will be used for proving stability of the spatial error estimator in Section~\ref{sec:stabred}. An alternative assumption involving mesh parameters was used for the same purpose in~\cite{Eigel2014a}; see Lemma~5.5 therein.
\end{remark}

We set $\X \coloneqq H^1_0(\Omega)$ with the corresponding norm
\begin{equation*}
 \norm{\varphi}_\X 
 \coloneqq \Big(\int_\Omega a_0(x) \abs{\nabla \varphi(x)}^2 \d x \Big)^{1/2}
 \quad \text{for all } \varphi \in \X.
\end{equation*}
Note that~\eqref{eq:a0bounds} implies
\begin{equation} \label{eq:normequivX}
 a_{\mathrm{min}} \, \norm{\nabla \varphi}_\Omega^2 
 \leq \norm{\varphi}_\X^2 
 \leq a_{\mathrm{max}} \, \norm{\nabla \varphi}_\Omega^2 
 \quad \text{for all } \varphi \in \X.
\end{equation}
We define the underlying Bochner space $\V \coloneqq L^2_\pi(\Lambda; \X)$ as well as the bilinear forms
\begin{equation*}
 B_0(v,w) 
 \coloneqq \int_{\Lambda} \int_{\Omega} a_0(x) \nabla v(x,\vec{y}) \cdot \nabla w(x,\vec{y}) \d x \d \pi(\vec{y}) 
 \quad \text{for all } v,w \in \V,
\end{equation*} 
and
\begin{equation*}
 B(v,w) 
 \coloneqq B_0(v,w) + \sum_{m \in \N} \int_\Lambda \int_{\Omega} y_m a_m(x) \nabla v(x,\vec{y}) \cdot \nabla w(x,\vec{y}) \d x \d \pi(\vec{y})
 \quad \text{for all } v,w \in \V.
\end{equation*} 
We note that by~\eqref{eq:a0bounds} and~\eqref{eq:asumbound} the induced (energy) norms $\enorm{\, \cdot \,} \coloneqq B(\cdot,\cdot)^{1/2}$, $\enorm{\, \cdot \,}_0 \coloneqq B_0(\cdot,\cdot)^{1/2}$ satisfy the equivalence
\begin{equation} \label{eq:normequiv}
 (1 - \gamma) \, \enorm{v}_0^2 \leq \enorm{v}^2 \leq (1 + \gamma) \, \enorm{v}_0^2 \qquad \text{for all } v \in \V.
\end{equation}%
Moreover, orthonormality of the basis $(P_\nu)_{\nu \in \JJJ}$ of $L^2_\pi(\Lambda)$ (cf. Section~\ref{subsec:gpc}) implies that
\begin{equation} \label{eq:parseval}
 \enorm{v}_0^2 
 = \sum_{\nu \in \JJJ} \norm{v_\nu}_\X^2
 \qquad \text{for all } v = \sum_{\nu \in \JJJ} v_\nu P_\nu \in \V.
\end{equation}%
The assumptions \eqref{eq:a0bounds} and \eqref{eq:asumbound} ensure that $B_0$ and $B$ are continuous and elliptic on $\V$ (see, e.g.,~\cite{Eigel2014,Bespalov2019}). 
Hence, by requiring $f \in L^2(\Omega)$ and setting
\begin{equation*}
 F(v) 
 \coloneqq
 \int_{\Lambda} \int_\Omega f(x) v(x,\vec{y}) \d x \d \pi(\vec{y})
 \quad \text{for all } v \in \V,
\end{equation*} 
the variational formulation 
\begin{equation*} 
 B(u^{\star},v) = F(v) \qquad \text{for all } v \in \V
\end{equation*} 
of~\eqref{eq:modelproblem} has a unique solution $u^\star \in \V$ by the Lax--Milgram lemma.

Let $\PPP_\ell \subseteq \JJJ$ be a finite index set and let $(\TT_{\ell,\nu})_{\nu \in \PPP_\ell}$ be a family of meshes. 
Recall the orthonormal basis $(P_\nu)_{\nu \in \JJJ}$ of $L^2_\pi(\Lambda)$ and define the semidiscrete space
\begin{equation*}
 \V_{\PPP_\ell} \coloneqq \X \otimes \mathrm{span}\set{P_\nu \given \nu \in \PPP_\ell} = \Big\{\sum_{\nu \in \PPP_\ell} v_{\nu} P_\nu \ | \ v_\nu \in \X \text{ for } \nu \in \PPP_\ell \Big\} \subset \V
\end{equation*}%
as well as the fully discrete space
\begin{equation*}
 \V_\ell \coloneqq \bigoplus_{\nu \in \PPP_\ell} \bigl(\SS_0^p(\TT_{\ell,\nu}) \otimes \mathrm{span}\set{P_\nu}\bigr) = \Big\{\sum_{\nu \in \PPP_\ell} v_{\ell,\nu} P_\nu \ | \ v_{\ell,\nu} \in \SS_0^p(\TT_{\ell,\nu}) \text{ for } \nu \in \PPP_\ell\Big\} \subset \V_{\PPP_\ell} \subset \V.
\end{equation*}%
Then, the Lax--Milgram lemma implies that for every multilevel structure $[\PPP_\ell,(\TT_{\ell,\nu})_{\nu \in \PPP_\ell}]$ there are unique solutions $u^\star_{\PPP_\ell} \in \V_{\PPP_\ell}$ and $u^\star_\ell \in \V_\ell$ to the semidiscrete and fully discrete problems
\begin{subequations}
\begin{align}
 B(u_{\PPP_\ell}^\star,v_{\PPP_\ell}) &= F(v_{\PPP_\ell}) \qquad \text{for all } v_{\PPP_\ell} \in \V_{\PPP_\ell},  \label{eq:weakformsemi} \\
 B(u_\ell^\star,v_\ell) &= F(v_\ell) \qquad \ \, \text{for all } v_\ell \in \V_\ell, \label{eq:weakformdisc}
\end{align}
\end{subequations}
respectively.

\section{A posteriori error estimation and main result} \label{sec:apost}

\noindent The error in SGFEM approximations consists of two components$\colon$ the parametric error and the spatial error. 
We follow the approach of \cite{Eigel2014,Eigel2014a} and consider residual-based \textsl{a~posteriori} indicators for the parametric and spatial error components. 
However, the parametric part of the error estimator proposed in \cite{Eigel2014,Eigel2014a} is not computable in general. 
We address this by introducing a saturation assumption in the parameter space, yielding computable error estimators for the energy error.

\subsection{Saturation assumption in the parameter space}

Given a finite index set $\PPP_\ell \subset \mathfrak{J}$, the corresponding semidiscrete space $\V_{\PPP_\ell}$ (cf. Section~\ref{subsec:sgfem}), and the associated detail index set $\mathfrak{Q}_{\ell}$, we consider the \textit{uniformly enriched index set} $\widehat{\mathfrak{P}}_{\ell} \coloneqq \PPP_\ell \cup \QQQ_{\ell}$ and the associated semidiscrete space $\widehat{\V}_{\PPP_\ell} \coloneqq \V_{\widehat{\mathfrak{P}}_\ell}$. 
Let $u_{\PPP_\ell}^{\star} \in \V_{\PPP_\ell}$ be the solution to the semidiscrete problem~\eqref{eq:weakformsemi} and let $\widehat{u}_{\PPP_\ell}^{\star} \coloneqq u_{\widehat{\mathfrak{P}}_\ell}^{\star} \in \widehat{\V}_{\PPP_\ell}$ solve~\eqref{eq:weakformsemi}, where $\V_{\PPP_\ell}$ is replaced by $\widehat{\V}_{\PPP_\ell}$.
We suppose that there is a uniform constant $0 < \const{q}{sat} < 1$ such that there holds the \textit{saturation assumption}
\begin{equation} \label{eq:sat}
 \enorm{u^{\star} - \widehat{u}^{\star}_{\PPP_\ell}}
 \leq \const{q}{sat} \, \enorm{u^{\star} - u_{\PPP_\ell}^{\star}}.
\end{equation}%
Galerkin orthogonality and the resulting Pythagorean identity imply
\begin{equation*}
 \enorm{u^{\star} - u_{\PPP_\ell}^{\star}}^2 
 = \enorm{u^{\star} - \widehat{u}_{\PPP_\ell}^{\star}}^2 + \enorm{\widehat{u}_{\PPP_\ell}^{\star} - u_{\PPP_\ell}^{\star}}^2
 \eqreff*{eq:sat}{\leq} \const{q}{sat}^2 \, \enorm{u^{\star} - u_{\PPP_\ell}^{\star}}^2 + \enorm{\widehat{u}_{\PPP_\ell}^{\star} - u_{\PPP_\ell}^{\star}}^2,
\end{equation*}%
which yields the norm equivalence
\begin{equation} \label{eq:satequiv}
 \enorm{\widehat{u}_{\PPP_\ell}^{\star} - u_{\PPP_\ell}^{\star}}^2
 \leq \enorm{u^{\star} - u_{\PPP_\ell}^{\star}}^2 
 \leq \frac{1}{1-\const{q}{sat}^2} \, \enorm{\widehat{u}_{\PPP_\ell}^{\star} - u_{\PPP_\ell}^{\star}}^2.
\end{equation}
For the discrete solution $u_\ell^{\star} \in \V_\ell$ satisfying~\eqref{eq:weakformdisc}, this yields
\begin{align}
  \notag \enorm{\widehat{u}_{\PPP_\ell}^\star - u_\ell^\star}^{2} &\leq
  \enorm{u^{\star} - u_\ell^{\star}}^2 
  = \enorm{u^{\star} - u_{\PPP_\ell}^{\star}}^2 + \enorm{u_{\PPP_\ell}^{\star} - u_\ell^{\star}}^2 \\
  \notag &\eqreff*{eq:satequiv}{\leq} (1 - \const{q}{sat}^2)^{-1} \, \enorm{\widehat{u}_{\PPP_\ell}^{\star} - u_{\PPP_\ell}^{\star}}^2 + \enorm{u_{\PPP_\ell}^{\star} - u_\ell^{\star}}^2 \\ 
  &\leq (1 - \const{q}{sat}^2)^{-1} \, \enorm{\widehat{u}_{\PPP_\ell}^ {\star} - u_\ell^{\star}}^2. \label{eq:fullerreq}
\end{align}
This equivalence allows us to consider $\enorm{\widehat{u}_{\PPP_\ell}^{\star} - u_{\ell}^{\star}}$ when deriving reliable error estimators for $\enorm{u^\star - u_\ell^\star}$. 

\begin{remark}
	The saturation assumption is a common requirement in the context of a~posteriori error estimation for elliptic PDEs, which may fail in general; see \cite[Proposition~2.2]{Bornemann1996}. Supposing~\eqref{eq:sat} for all semidiscrete subspaces $\V_{\PPP_\ell}$ is a strong premise, however, in the present work it is only required for a nested sequence of finite sets $\PPP_\ell$ generated by an adaptive loop.
\end{remark}

\subsection{Expansion of the semidiscrete residual}

Let $[\PPP_\ell,(\TT_{\ell,\nu})_{\nu \in \PPP_\ell}]$ be a multilevel structure and let $u_\ell^\star \in \V_\ell$ be the solution to~\eqref{eq:weakformdisc}. 
We follow the approach of \cite{Eigel2014,Eigel2014a}, but instead of the continuous residual, we consider the semidiscrete residual $\RR_\ell \in \V^*$ defined by 
\begin{equation*}
 \dual{\RR_\ell}{v}_{\V^*} 
 \coloneqq B(\widehat{u}_{\PPP_\ell}^{\star}-u_\ell^\star,v)
 \quad \text{for all } v \in \V.
\end{equation*} 
Similarly to \cite[Section~5]{Eigel2014}, elementary calculations show that the residual can be decomposed into contributions associated with individual multi-indices, i.e., for every $\nu \in \JJJ$, there exists a coefficient $r_{\nu} \in \X^*$ such that
\begin{equation*}
 \dual{\RR_\ell}{v}_{\V^*}
 = \sum_{\nu \in \JJJ} \dual{r_{\nu}}{v_{\nu}}_{\X^*}
 \quad \text{for all } v = \sum_{\nu \in \JJJ} v_{\nu} P_{\nu} \in \V.
\end{equation*}
Let $v \in \V$ and denote by $\widehat{G}_\ell \colon \V \to \widehat{\V}_{\PPP_\ell}$ the orthogonal Galerkin projection with respect to the inner product $B(\cdot,\cdot)$. We obtain $\enorm{\widehat{G}_\ell v} \leq \enorm{v}$ and $\widehat{G}_\ell(\widehat{u}_{\PPP_\ell}^{\star} - u_\ell^\star) = \widehat{u}_{\PPP_\ell}^{\star} - u_\ell^{\star}$. Hence, we obtain
\begin{equation*}
 \begin{split}
  \frac{\abs{\dual{\RR_\ell}{v}_{\V^*}}}{\enorm{v}} 
  &= \frac{\abs{B(\widehat{u}_{\PPP_\ell}^{\star}-u_\ell^\star,v)}}{\enorm{v}}
  = \frac{\abs{B(\widehat{u}_{\PPP_\ell}^{\star}-u_\ell^\star,\widehat{G}_\ell v)}}{\enorm{v}}
  \leq \frac{\abs{B(\widehat{u}_{\PPP_\ell}^{\star}-u_\ell^\star,\widehat{G}_\ell v)}}{\enorm{\widehat{G}_\ell v}}.
 \end{split}
\end{equation*}
This, together with the estimates~\eqref{eq:parseval} and~\eqref{eq:normequiv}, yields
\begin{align} 
 \notag \enorm{\widehat{u}_{\PPP_\ell}^{\star}-u_\ell^\star} 
 &= \sup_{v \in \V \setminus \set{0}} \frac{\abs{\dual{\RR_\ell}{v}_{\V^*}}}{\enorm{v}}
 \leq \sup_{v \in \V \setminus \set{0}} \frac{\abs{\dual{\RR_\ell}{\widehat{G}_\ell v}_{\V^*}}}{\enorm{\widehat{G}_\ell v}}
 \leq \sup_{v \in \widehat{\V}_{\PPP_\ell} \setminus \set{0}} \frac{\abs{\dual{\RR_\ell}{v}_{\V^*}}}{\enorm{v}} \\
 \notag &\leq \sup_{v \in \widehat{\V}_{\PPP_\ell} \setminus \set{0}} \enorm{v}^{-1} \sum_{\nu \in \widehat{\PPP}_\ell} \abs{\dual{r_{\nu}}{v_{\nu}}_{\X^*}}
 \leq \sup_{v \in \widehat{\V}_{\PPP_\ell} \setminus \set{0}} \enorm{v}^{-1} \sum_{\nu \in \widehat{\PPP}_\ell} [\norm{r_{\nu}}_{\X^*} \norm{v_{\nu}}_{\X}] \\
 \notag &\leq \sup_{v \in \widehat{\V}_{\PPP_\ell} \setminus \set{0}} \enorm{v}^{-1} \Bigl(\sum_{\nu \in \widehat{\PPP}_\ell} \norm{r_{\nu}}_{\X^*}^2\Bigr)^{1/2} \Bigl(\sum_{\nu \in \widehat{\PPP}_\ell} \norm{v_{\nu}}_{\X}^2\Bigr)^{1/2} \\
 &\eqreff*{eq:parseval}{\leq} \sup_{v \in \widehat{\V}_{\PPP_\ell} \setminus \set{0}} \enorm{v}^{-1} \enorm{v}_0 \Bigl(\sum_{\nu \in \widehat{\PPP}_\ell} \norm{r_{\nu}}_{\X^*}^2\Bigr)^{1/2} 
 \eqreff{eq:normequiv}{\leq} (1 - \gamma)^{-1/2} \Bigl(\sum_{\nu \in \widehat{\PPP}_\ell} \norm{r_{\nu}}_{\X^*}^2\Bigr)^{1/2}, \label{eq:resnorm}
\end{align}
which allows us to consider the (finite) sum $\sum_{\nu \in \widehat{\PPP}_\ell} \norm{r_\nu}_{\X^*}^2$ when deriving reliable error estimators for $\enorm{\widehat{u}_{\PPP_\ell}^{\star}-u_\ell^\star}^2$.

\subsection{Computable residual-based error estimator} \label{sec:esti}
Let $\nu \in \PPP_\ell$. Set 
\begin{equation} \label{eq:Toverlaydef}
 \TT_\ell[\nu] 
 \coloneqq \TT_{\ell,\nu} \oplus \bigoplus_{\mu \in \PPP_\ell[\nu]} \TT_{\ell,\mu}.
\end{equation}%
Note that for all $w_\ell = \sum_{\nu \in \PPP_\ell} w_{\ell,\nu} P_\nu \in \V_\ell$ and $\mu \in \PPP_\ell[\nu] \cup \set{\nu}$, one has $w_{\ell,\mu} \in \SS_0^p(\TT_\ell[\nu])$.
Given $w_\ell = \sum_{\nu \in \PPP_\ell} w_{\ell,\nu} P_\nu \in \V_\ell$, we define for $\nu \in \PPP_\ell$ and $T \in \TT_{\ell,\nu}$ the \emph{spatial error indicator} as 
\begin{equation} \label{eq:spatialdef}
 \tau_{\ell}(\nu,T;w_\ell) 
 \coloneqq \Bigl(\sum_{\substack{T' \in \TT_\ell[\nu] \\ T' \subseteq T}} \Bigl[ h_T^2 \norm{\delta_{\vec{0} \nu} f + \div \sigma_\nu(w_\ell)}_{T'}^2 + h_T \norm{\jump{\sigma_\nu(w_\ell)} \cdot \vec{n}_{\partial T'}}_{\partial T'\cap \Omega}^2 \Bigr]\Bigr)^{1/2},
\end{equation}
where, with the constants $\beta_m^n$ from~\eqref{eq:recurrence},
\begin{equation} \label{eq:sigmadef}
 \sigma_{\nu}(w_{\ell}) 
 \coloneqq a_0\nabla w_{\ell,\nu} + \sum_{m \in \N} a_m \bigl( \beta^{\nu_m}_m \nabla w_{\ell,\nu + \eps_m} + \beta^{\nu_m-1}_m \nabla w_{\ell,\nu - \eps_m} \bigr).
\end{equation} 
Here, we set $w_{\nu} \coloneqq 0$ for $\nu \notin \JJJ$. Furthermore, given $\nu \in \QQQ_\ell$, we define the computable \emph{parametric error indicator} as
\begin{equation} \label{eq:paramdef}
 \zeta_\ell(\nu;w_\ell)
 \coloneqq \sum_{m \in \N} \norm{a_m}_{W^{1,\infty}(\Omega)} \bigl(\beta_m^{\nu_m} \norm{\nabla w_{\ell,\nu + \eps_m}}_\Omega + \beta_m^{\nu_m-1} \norm{\nabla w_{\ell,\nu - \eps_m}}_{\Omega}\bigr).
\end{equation}
Moreover, given $\SS_{\ell,\nu} \subseteq \TT_{\ell,\nu}$ for $\nu \in \PPP_\ell$, $\mathfrak{U}_\ell \subseteq \PPP_\ell$, and $\mathfrak{W}_\ell \subseteq \QQQ_\ell$, we set 
\begin{equation} \label{eq:estnot}
 \begin{split}
  \tau_{\ell}(\nu,\SS_{\ell,\nu};w_\ell)^2 
  &\coloneqq \sum_{T \in \SS_{\ell,\nu}} \tau_\ell(\nu,T;w_\ell)^2, \quad \tau_\ell((\SS_{\ell,\nu})_{\nu \in \PPP_\ell};w_\ell)^2 
  \coloneqq \sum_{\nu \in \PPP_\ell} \tau_\ell(\nu,\SS_{\ell,\nu};w_\ell)^2, \\
  \tau_\ell(\nu;w_\ell)^2 
  &\coloneqq \tau_\ell(\nu,\TT_{\ell,\nu};w_\ell)^2, 
  \qquad \qquad \quad \ \ \,  \tau_\ell(\mathfrak{U}_\ell;w_\ell)^2 
  \coloneqq \sum_{\nu \in \mathfrak{U}_\ell} \tau_\ell(\nu;w_\ell)^2,\\ 
  \tau_\ell(w_\ell)^2 
  &\coloneqq \tau_{\ell}(\PPP_\ell;w_\ell)^2, 
  \qquad \qquad \qquad \ \  \zeta_\ell(\mathfrak{W}_\ell;w_\ell)^2 
  \coloneqq \sum_{\nu \in \mathfrak{W}_\ell} \zeta_\ell(\nu;w_\ell)^2,  \\
  \zeta_\ell(w_\ell)^2 &\coloneqq \zeta_\ell(\QQQ_\ell;w_\ell)^2, 
  \qquad \qquad \qquad \qquad \ \, \eta_\ell(w_\ell)^2 
  \coloneqq \tau_\ell(w_\ell)^2 + \zeta_\ell(w_\ell)^2.
 \end{split}
\end{equation}
For $w_\ell = u_\ell^\star$ we omit the last argument and simply write $\tau_\ell(\nu,T) \coloneqq \tau_\ell(\nu,T;u_\ell^\star)$, $\zeta_\ell(\mathfrak{W}_\ell) \coloneqq \zeta_\ell(\mathfrak{W}_\ell;u_\ell^\star)$, $\eta_\ell \coloneqq \eta_\ell(u_\ell^\star)$  and so forth. Analogous arguments as in \cite[Section~6]{Eigel2014} lead to 
\begin{equation} \label{eq:locrel}
 \begin{split}
  \norm{r_{\nu}}_{\X^*} 
  &\leq a_{\mathrm{min}}^{-1/2} \const{\widetilde{C}}{rel} \, \tau_\ell(\nu)
  \quad \text{for all } \nu \in \PPP_\ell, \\
  \norm{r_{\nu}}_{\X^*} 
  &\leq a_{\mathrm{min}}^{-1/2} \, \zeta_\ell(\nu)   \quad \text{for all } \nu \in \QQQ_\ell,
 \end{split}
\end{equation}
where $\const{\widetilde{C}}{rel} > 0$ depends only on $\kappa$-shape regularity of $\TT_0$. The estimate~\eqref{eq:fullerreq} together with~\eqref{eq:resnorm} and~\eqref{eq:locrel} yield \emph{reliability} of the error estimator $\eta_\ell$, i.e.,
\begin{equation} \label{eq:reliability}
 \enorm{u^{\star} - u_\ell^{\star}} 
 \leq \const{C}{rel} \eta_\ell
 \quad \text{with} \quad 
 \const{C}{rel} 
 \coloneqq \Bigl(\frac{\max \set{1,\const{\widetilde{C}}{rel}^2}}{a_{\mathrm{min}}(1 - \const{q}{sat}^2)(1 - \gamma)}\Bigr)^{1/2}.
\end{equation}\subsection{Adaptive multilevel algorithm and main result}

Throughout the remainder of this paper, we denote by $(K_\ell)_{\ell \in \N_0}$ a sequence of compatibility degrees such that $0 <  K_{\ell - 1} \leq K_\ell$ for all $\ell \in \N$. We require $(K_\ell)_{\ell \in \N_0}$ to not grow too fast in the sense that 
\begin{equation} \label{eq:sumdiv}
 \sum_{\ell \in \N_0} 2^{-2K_\ell} = +\infty.
\end{equation} 
\begin{remark}
 We purposefully exclude the case $K_\ell = 0$ from the analysis (but not from the numerics), since treating it within the present multilevel framework would necessitate additional case distinctions, considerably lengthening the proofs. As plain convergence---and even R-linear convergence---of the adaptive algorithm for $K_\ell = 0$ (single-level SGFEM) has already been established in~\cite{Eigel2014a}, this restriction entails no essential loss of generality.
\end{remark}
Additionally, we suppose that there exists a multilevel refinement routine $\mathtt{refine}$ which satisfies the following properties$\colon$ Given a $K_\ell$-compatible multilevel structure $[\PPP_\ell,(\TT_{\ell,\nu})_{\nu \in \PPP_\ell}]$, a family of marked elements $(\MM_{\ell,\nu})_{\nu \in \PPP_\ell}$ with $\MM_{\ell,\nu} \subseteq \TT_{\ell,\nu}$ for all $\nu \in \PPP_\ell$, and a set of marked multi-indices $\MMM_\ell \subseteq \QQQ_\ell$, the routine $\mathtt{refine}$ returns a \emph{refined} multilevel structure $[\PPP_{\ell+1},(\TT_{\ell+1,\nu})_{\nu \in \PPP_{\ell+1}}]$ such that 
\begin{itemize}
 \item $\PPP_\ell \cup \MMM_\ell \subseteq \PPP_{\ell + 1} \subseteq \PPP_\ell \cup \QQQ_\ell$,
 \medskip
 \item for every $\nu \in \PPP_\ell$, the mesh $\TT_{\ell + 1,\nu}$ is a refinement of $\TT_{\ell,\nu}$ such that at least all elements in $\MM_{\ell,\nu}$ have been refined, and
 \item for every $\nu \in \PPP_{\ell+1}$, the mesh $\TT_{\ell+1,\nu}$ is the coarsest triangulation such that the meshes $(\TT_{\ell+1,\nu})_{\nu \in \PPP_{\ell + 1}}$ are $K_{\ell + 1}$-compatible. 
\end{itemize}
A construction of such a refinement routine is given in Section~\ref{sec:algo}.

Based on the aforementioned refinement routine and the error indicators introduced in Section~\ref{sec:esti}, we propose the following adaptive algorithm.
\begin{algorithm}[Adaptive algorithm] \label{alg:adapt}
\textbf{Input:} Initial mesh $\TT_0$, initial multilevel structure $[\set{\vec{0}},(\TT_0)]$, polynomial degree $p \in \N$ for the FEM, sequence of compatibility degrees $(K_\ell)_{\ell \in \N_0}$ satisfying~\eqref{eq:sumdiv}, marking parameter $0 < \theta \leq 1$, diffusion coefficients $a_m \in W^{1,\infty}(\Omega)$, $m \in \N_0$, and the right-hand side $f \in L^2(\Omega)$. \\
\textbf{Loop:} For all $\ell = 0,1,2,\ldots,$ repeat the following steps:
\begin{enumerate}[label = \rm(\roman*)]
 \item Compute the discrete solution $u_\ell^\star \in \V_\ell$ to~\eqref{eq:weakformdisc}.

 \item Compute the spatial indicators $\tau_\ell(\nu,T)$ given by~\eqref{eq:spatialdef} for all $\nu \in \PPP_\ell$ and $T \in \TT_{\ell,\nu}$.

 \item Compute the parametric indicators $\zeta_\ell(\nu)$ given by~\eqref{eq:paramdef} for all $\nu \in \QQQ_\ell$.

 \item Determine sets of marked elements $\MM_{\ell,\nu} \subseteq \TT_{\ell,\nu}$ for all $\nu \in \PPP_\ell$ and a set of marked multi-indices $\MMM_\ell \subseteq \QQQ_\ell$ (usually of minimal cardinality $\bigl(\sum_{\nu \in \PPP_\ell} \abs{\MM_{\ell,\nu}}\bigr) + \abs{\MMM_\ell}$) such that
 \begin{equation} \label{eq:doerfler}
  \tau_\ell\bigl((\MM_{\ell,\nu})_{\nu \in \PPP_\ell};u_\ell^\star\bigr)^2 + \zeta_\ell(\MMM_\ell;u_\ell^\star)^2 
  \geq \theta \, \eta_\ell^2.
 \end{equation}%

 \item Generate a $K_{\ell + 1}$-compatible multilevel structure 
 \begin{equation} \label{eq:refine}
  [\PPP_{\ell + 1},(\TT_{\ell + 1,\nu})_{\nu \in \PPP_{\ell + 1}}] 
  = \mathtt{refine}([\PPP_\ell,(\TT_{\ell,\nu})_{\nu \in \PPP_\ell}],(\MM_{\ell,\nu})_{\nu \in \PPP_\ell},\MMM_\ell,K_{\ell + 1}).
 \end{equation}

\item Increase the counter $\ell$ by one and go to step (i).
\end{enumerate}
\end{algorithm}

The following convergence theorem is the main result of this work. The subsequent section is devoted to its proof.
\begin{theorem}[Plain convergence] \label{thm:main}
Under the foregoing assumptions, Algorithm~\ref{alg:adapt} guarantees plain convergence of the error 
\begin{equation}
 \enorm{u^{\star} - u_\ell^{\star}} 
 \xrightarrow{\ell \to \infty} 0.
\end{equation}
Moreover, if $(K_\ell)_{\ell \in \N_0}$ is bounded, i.e., there exists $K > 0$ such that $K_\ell \leq K$ for all $\ell \in \N_0$, then there holds estimator convergence
\begin{equation}
 \eta_\ell(u_\ell^\star) 
 \xrightarrow{\ell \to \infty} 0.
\end{equation}
\end{theorem}
\begin{remark}[Extensions of the analysis]
 Theorem~\ref{thm:main} still remains valid in more general settings$\colon$
 \begin{itemize}
  \item Uncertain right-hand side$\colon$ 
  \begin{equation*}
   f(x,\vec{y}) 
   = f_0(x) + \sum_{m \in \N} y_m f_m(x) 
   \in L^2(\Lambda;L^2(\Omega)).
  \end{equation*}
  In this case, the error indicators~\eqref{eq:spatialdef} and~\eqref{eq:paramdef} need to be adjusted accordingly.
  \item Piecewise-$W^{1,\infty}(\Omega)$ diffusion coefficients$\colon$
  For any $m \in \N_0$, the expansion coefficient $a_m \in L^\infty(\Omega)$ in~\eqref{eq:affine} may be piecewise $W^{1,\infty}(\Omega)$, i.e., there exists a partition $(\Omega^m_j)_{j = 1}^{J_m}$ of $\Omega$ into Lipschitz domains such that $a_m|_{\Omega^m_j} \in W^{1,\infty}(\Omega^m_j)$ for all $j = 1,\ldots,J_m$. 
  We emphasize that the partitions $(\Omega^m_j)_{j = 1}^{J_m}$ may depend on $m \in \N_0$ and that assumption~\eqref{eq:ampliesti} needs to be adjusted accordingly.
  \item Extended detail index set$\colon$ 
  Given a fixed number $N \in \N$, the detail index set $\QQQ_\ell$ in~\eqref{eq:detaildef} may be replaced by
  \begin{equation*}
   \begin{split}
    \mathfrak{Q}_\ell 
    \coloneqq \big\{ \nu \in \JJJ \setminus \PPP_\ell \; \vert \; &\exists \, m \in \set{1,2,\ldots,M_\ell} \text{ such that } \nonumber \\
    &\nu + \eps_m \in \PPP_\ell \text{ or } \nu - \eps_m \in \PPP_\ell \big\} \cup \set{\eps_{M_\ell + 1}, \ldots, \eps_{M_\ell + N}}.
   \end{split}
  \end{equation*}
  \item Separate marking$\colon$ 
  Instead of~\eqref{eq:doerfler}, one may consider separate marking for the spatial and parametric error indicators, i.e., given $0 < \theta_\tau,\theta_\zeta \leq 1$, choose $\MM_{\ell,\nu} \subseteq \TT_{\ell,\nu}$ for all $\nu \in \PPP_\ell$ and $\MMM_\ell \subseteq \QQQ_\ell$ (usually of minimal cardinality) such that
  \begin{equation*}
   \tau_\ell((\MM_{\ell,\nu})_{\nu \in \PPP_\ell};u_\ell^\star)^2 
   \geq \theta_\tau \tau_\ell(u_\ell^\star)^2
   \quad \text{and} \quad
   \zeta_\ell(\MMM_\ell;u_\ell^\star)^2 
   \geq \theta_\zeta \zeta_\ell(u_\ell^\star)^2.
  \end{equation*}
 \end{itemize}
 Although these extensions require only minor modifications in the proofs, we refrain from presenting the details for the sake of clarity and readability.
\end{remark}
\section{Proof of Theorem~\ref{thm:main} (plain convergence)} \label{sec:proof}

\noindent The proof of Theorem~\ref{thm:main} crucially hinges on the \emph{stability} and \emph{reduction} properties of the error estimator $\eta_\ell$ (see, e.g.,~\cite{Carstensen2014}). 
Subsequently, we consider the quasi-error
\begin{equation} \label{eq:auxquasi}
 \mathsf{H}_\ell^2 
 \coloneqq \enorm{u^\star - u_\ell^\star}^2 +  2^{ - 4 K_\ell/d} \, \tau_\ell^2 + \mu \, 2^{ - 2 K_\ell/d} \, \zeta_\ell^2
\end{equation}
with a free parameter $\mu > 0$ and we show that there exists a suitable choice for $\mu > 0$ such that $\mathsf{H}_\ell$ is indeed quasi-contractive up to a perturbation.
This will eventually lead to the proof of Theorem~\ref{thm:main}.

Throughout the present section, we assume that all involved quantities (discrete solutions, error estimators, multilevel structures etc.) are generated by Algorithm~\ref{alg:adapt}.

The proof of Theorem~\ref{thm:main} requires the following proposition, which is proved in Section~\ref{sec:stabred} but stated here for convenience.
\begin{proposition}[Stability and reduction] \label{prop:stabred}
{\bf Stability.} Under the assumptions of Section~\ref{sec:apost}, there exist constants $C_\zeta > 0$ depending only on the diffusion coefficient $a$ and $C_\tau > 0$ depending only on $\Omega$, the dimension $d$, the initial mesh $\TT_0$, the polynomial degree $p$, and the diffusion coefficient $a$ such that for all $\ell \in \N_0$ there holds
\begin{subequations} \label{eq:stab}
 \begin{align}
  \tau_\ell(w_\ell) 
  &\leq \tau_\ell(v_\ell) + C_\tau 2^{K_\ell/d} \, \enorm{w_\ell - v_\ell} \quad \text{for all } w_\ell, v_\ell \in \V_\ell, \label{eq:taustab} \\
  \zeta_\ell(w_\ell) 
  &\leq \zeta_\ell(v_\ell) + C_\zeta \, \enorm{w_\ell - v_\ell} 
  \qquad \quad \, \text{for all } w_\ell, v_\ell \in \V_\ell. \label{eq:zetastab}
 \end{align}\end{subequations}
{\bf Reduction.} Moreover, there exists a constant $\const{C}{red} > 0$ depending only on $\Omega$, the dimension $d$, the initial mesh $\TT_0$, the polynomial degree $p$, and the diffusion coefficient $a$ such that, for $\const{q}{red} = 1 - 2^{-1/d}$, there holds, for all $\ell \in \N_0$ and $w_\ell \in \V_\ell$, that
\begin{subequations} \label{eq:red}
\begin{equation} \label{eq:taured}
 \tau_{\ell+1}(w_\ell)^2 
 \leq \tau_\ell(w_\ell)^2 + \const{C}{red} 2^{2K_{\ell+1}/d} \, \zeta_\ell(\QQQ_\ell \cap \PPP_{\ell+1};w_\ell)^2 - \const{q}{red} \tau_\ell((\TT_{\ell,\nu} \setminus \TT_{\ell+1,\nu})_{\nu \in \PPP_\ell};w_\ell)^2
\end{equation}%
as well as
\begin{equation} \label{eq:zetared}
 \zeta_{\ell + 1}(w_\ell)^2 
 = \zeta_\ell(w_\ell)^2 + \alpha_{\ell}(w_\ell)^2 - \zeta_\ell(\QQQ_\ell \cap \PPP_{\ell+1};w_\ell)^2,
\end{equation}
\end{subequations}
where, for all $\ell \in \N_0$, $w_\ell \in \V_\ell$, and $M_\ell = \abs{\supp{\PPP_\ell}}$,
\begin{equation} \label{eq:alphadef}
 \alpha_\ell(w_\ell) 
 \coloneqq 
 \begin{cases} 
  \beta_{M_{\ell + 1} + 1}^0 \norm{a_{M_{\ell + 1} + 1}}_{W^{1,\infty}(\Omega)} \norm{\nabla w_{\ell,\vec{0}}}_\Omega 
  & \text{if } M_{\ell + 1} > M_{\ell}, \\
  0 
  & \text{if } M_{\ell + 1} = M_{\ell}.
 \end{cases}
\end{equation}

\end{proposition}

The proof of plain convergence relies on the following criterion, which results from summability properties of perturbed contractions.
\begin{lemma}[Summability criterion] \label{lem:sumcrit}
Let $(a_\ell)_{\ell \in \N_0}, (b_\ell)_{\ell \in \N_0}, (c_\ell)_{\ell \in \N_0}$ be non-negative sequences. 
Assume that
\begin{equation} \label{eq:contract}
 \begin{split}
  a_{\ell+1} 
  &\leq (1 - c_\ell) \, a_\ell + b_\ell 
  \quad  \text{for all } \ell \in \N_0,
 \end{split}
\end{equation}
and
\begin{equation} 
 \sum_{\ell = 0}^\infty b_\ell 
 < +\infty 
 \quad \text{and} \quad 
 \sum_{\ell = 0}^\infty c_\ell 
 = +\infty.
\end{equation}
Then, $(a_\ell)_{\ell \in \N_0}$ converges to $0$, i.e.,
\begin{equation}
 a_{\ell} 
 \xrightarrow{\ell \to \infty} 0.
\end{equation}
\end{lemma}

\begin{proof}
We first note that, if $c_\ell \geq 1$ for finitely many $\ell \in \N_0$, the convergence behavior of the involved sequences remains unchanged. If $c_\ell \geq 1$ for infinitely many $\ell \in \N_0$, we have $0 \leq a_{\ell + 1} \leq b_\ell$.
Since $(b_\ell)_{\ell \in \N_0}$ is summable, $b_\ell \xrightarrow{\ell \to \infty} 0$ and thus also $a_{k(\ell)} \xrightarrow{\ell \to \infty} 0$, where $\set{k(\ell) \given \ell \in \N_0} = \set{\ell \in \N \given c_{\ell-1} \geq 1}$. It remains to show that the complementary subsequence $(a_{m(\ell)})_{\ell \in \N_0}$, where $\set{m(\ell) \given \ell \in \N_0} = \set{\ell \in \N \given c_{\ell - 1} < 1}$, converges to zero.
Thus, without loss of generality, suppose that $c_\ell < 1$ for all $\ell \in \N_0$.

Let $n \in \N$. There holds
\begin{equation*}
 \begin{split}
  \sum_{\ell = 1}^n a_\ell 
  &= \sum_{\ell = 0}^{n-1} a_{\ell+1}
  \eqreff{eq:contract}{\leq} \sum_{\ell = 0}^{n-1} \bigl((1 - c_\ell) a_\ell + b_\ell\bigr)
  \leq \sum_{\ell = 0}^n (1 - c_{\ell}) a_{\ell} + \sum_{\ell = 0}^\infty b_\ell,
 \end{split}
\end{equation*}%
which implies 
\begin{equation} \label{eq:acsum}
 \sum_{\ell = 1}^n c_{\ell} a_\ell 
 \leq (1 - c_0) a_0 + \sum_{\ell = 0}^\infty b_\ell 
 < +\infty 
 \quad \text{for all } n \in \N.
\end{equation}

We argue by contradiction and assume that $(a_\ell)_{\ell \in \N_0}$ does not converge to $0$. Then, there exists $\eps > 0$ such that, for all $n \in \N_0$, there exists $\ell_n \geq n$ with $a_{\ell_n} \geq \eps$. Since $(b_\ell)_{\ell \in \N_0}$ is summable, there exists $n_0 \in \N_0$ such that 
\begin{equation*}
 \sum_{\ell = n_0}^\infty b_\ell 
 < \eps/2.
\end{equation*}
For any $n \geq n_0$, the recursion~\eqref{eq:contract} and the fact that $c_\ell < 1$ for all $\ell \in \N_0$ imply 
\begin{equation*}
 \begin{split}
  \eps 
  \leq a_{\ell_n} 
  &\leq \Bigl( \prod_{k = n}^{\ell_n - 1} (1 - c_k) \Bigr) a_n + \sum_{k = n}^{\ell_n - 1} \Bigl( \prod_{j = n+1}^{k} (1 - c_{\ell_n - j + n}) \Bigr) b_{\ell_n - 1 - k + n} \\
  &\leq a_n + \sum_{k = n}^{\ell_n - 1} b_{\ell_n - 1 - k + n} 
  \leq a_n + \sum_{\ell = n}^\infty b_\ell 
  < a_n + \eps/2,
 \end{split}
\end{equation*}
which concludes that $a_n > \eps/2$ for all $n \geq n_0$. With~\eqref{eq:acsum}, we obtain that
\begin{equation*}
 \sum_{\ell = n_0}^\infty c_\ell 
 < 2 \, \eps^{-1} \sum_{\ell = n_0}^\infty c_\ell \, a_\ell 
 < +\infty,
\end{equation*} 
which contradicts 
\begin{equation*}
 \sum_{\ell = 0}^\infty c_\ell 
 = +\infty.
\end{equation*}	
This concludes the proof of $a_\ell \xrightarrow{\ell \to \infty} 0$.
\end{proof}
By means of the foregoing results, we are now able to prove that the quasi-error $\mathsf{H}_\ell^2$ from~\eqref{eq:auxquasi} is quasi-contractive in the sense of~\eqref{eq:contract}.
\begin{lemma}[Quasi-error is quasi-contractive] \label{lem:quasi}
 There are constants $\delta, \const{C}{per} > 0$ and a free parameter $\mu > 0$ depending only on the initial mesh $\TT_0$, the polynomial degree $p$, the diffusion coefficient $a$, the saturation constant $\const{q}{sat}$, the marking parameter $\theta$, the dimension $d$, and $K_0$ such that the quasi-error from~\eqref{eq:auxquasi} with $\alpha_\ell(u_\ell^\star)$ as in~\eqref{eq:alphadef} satisfies
 \begin{equation} \label{eq:auxconv}
  \mathsf{H}_{\ell + 1}^2 
  \leq (1 - \delta \, 2^{-4K_{\ell + 1}/d}) \mathsf{H}_{\ell}^2 + \const{C}{per} \big[\alpha_{\ell}(u_{\ell}^\star)^2 + \enorm{u_{\ell+1}^\star - u_{\ell}^\star}^2\big]			
  \quad \text{for all } \ell \in \N_0.
 \end{equation} 
\end{lemma}

\begin{proof}
The proof is split into four steps. 

\textbf{Step~1 (Definition of \boldsymbol{$\delta$}).}
Let $\mu > 0$ satisfy
\begin{equation} \label{eq:muchoice}
 \mu 
 > \const{C}{red} + \const{q}{red}.
\end{equation}
Since $\theta \, \const{q}{red} < 1$ for any $0 < \theta \leq 1$, we can choose $\eps > 0$ to satisfy
\begin{equation*}
 \eps 
 < \theta \, \const{q}{red} \min \set{(1 - \theta \, \const{q}{red})^{-1},\mu^{-1}}.
\end{equation*}
Hence, $(1 + \eps)(1 - \theta \, \const{q}{red}) < 1$ and $\theta \, \const{q}{red} > \eps \mu$. Finally, let $\rho > 0$ be such that
\begin{equation*}
 \rho 
 < \const{C}{rel}^{-2} \min \set{1 - (1 + \eps)(1 - \theta \, \const{q}{red}), \theta \, \const{q}{red} - \eps \mu}.
\end{equation*}
We define
\begin{equation} \label{eq:deltadef}
 \delta 
 \coloneqq \min \set{\rho, \mu^{-1}(\theta \const{q}{red} - \eps \mu - \rho\const{C}{rel}^2),2^{4K_0/d} \, [1 - (1 + \eps)(1 - \theta \const{q}{red}) - \rho\const{C}{rel}^2]}.
\end{equation}%
Due to our choices of $\mu,\eps,\rho$, there holds $\delta > 0$.

\textbf{Step~2 (Estimating the weighted estimator).}
We abbreviate
\begin{equation} \label{eq:labbrev}
 \mu_\ell 
 \coloneqq \mu \, 2^{-2K_\ell/d},
 \quad \eps_\ell 
 \coloneqq \eps \, 2^{-2K_\ell/d}
 \quad \text{and} \quad 
 \rho_\ell 
 \coloneqq \rho \, 2^{-4K_\ell/d}
 \quad \text{for all } \ell \in \N_0,
\end{equation}
and
\begin{equation} \label{eq:etatildedef}
 \widetilde{\eta}_\ell(w_\ell)^2 
 \coloneqq 2^{-4K_\ell/d} \, \tau_\ell(w_\ell)^2 + \mu_\ell \, \zeta_\ell(w_\ell)^2 
 \quad \text{for all } \ell \in \N_0 \text{ and } w_\ell \in \V_\ell,
\end{equation}%
where we again use the convention that for $w_\ell = u_\ell^\star$, we simply write $\widetilde{\eta}_\ell \coloneqq \widetilde{\eta}_\ell(u_\ell^\star)$.
Stability~\eqref{eq:stab}, Young's inequality, and reduction~\eqref{eq:red} yield
\begin{align} \label{eq:taustabred}
 \begin{split}
  \tau_{\ell+1}^2 
  &\eqreff*{eq:taustab}{\leq} (\tau_{\ell+1}(u_{\ell}^\star) + C_\tau \, 2^{K_{\ell+1}/d} \, \enorm{u_{\ell+1}^\star - u_{\ell}^\star})^2 \\
  &\leq (1 + \eps_{\ell+1}) \, \tau_{\ell+1}(u_{\ell}^\star)^2 + (1 + \eps_{\ell+1}^{-1}) \, C_\tau^2 \, 2^{2K_{\ell+1}/d} \, \enorm{u_{\ell+1}^\star - u_{\ell}^\star}^2 \\
  &\eqreff*{eq:taured}{\leq} (1 + \eps_{\ell+1}) \, \tau_{\ell}^2 - (1 + \eps_{\ell+1})\, \const{q}{red} \, \tau_\ell((\TT_{\ell,\nu} \setminus \TT_{\ell+1,\nu})_{\nu \in \PPP_{\ell}})^2 \\
  &\qquad  + (1 + \eps_{\ell+1}) \, \const{C}{red} \, 2^{2K_{\ell+1}/d} \, \zeta_\ell(\QQQ_{\ell} \cap \PPP_{\ell+1})^2 \\
  &\qquad + (1 + \eps_{\ell+1}^{-1}) \, C_\tau^2 \, 2^{2K_{\ell+1}/d} \,\enorm{u_{\ell+1}^\star - u_{\ell}^\star}^2
 \end{split}
\end{align}
and
\begin{equation} \label{eq:zetastabred}
 \begin{split}
  \zeta_{\ell+1}^2 
  &\eqreff*{eq:zetastab}{\leq} (\zeta_{\ell+1}(u_{\ell}^\star) + C_\zeta \, \enorm{u_{\ell+1}^\star - u_{\ell}^\star})^2 \\
  &\leq (1 + \eps_{\ell+1}) \, \zeta_{\ell+1}(u_{\ell}^\star)^2 + (1 + \eps_{\ell+1}^{-1}) \, C_\zeta^2 \, \enorm{u_{\ell+1}^\star - u_{\ell}^\star}^2 \\
  &\eqreff*{eq:zetared}{=} (1 + \eps_{\ell+1}) \, \zeta_{\ell}^2 - (1 + \eps_{\ell+1}) \, \zeta_\ell(\QQQ_{\ell} \cap \PPP_{\ell+1})^2 \\
  &\qquad + (1 + \eps_{\ell+1}) \, \alpha_{\ell}(u_{\ell}^\star)^2 + (1 + \eps_{\ell+1}^{-1}) \, C_\zeta^2 \, \enorm{u_{\ell+1}^\star - u_{\ell}^\star}^2.
 \end{split}
\end{equation}
We note that~\eqref{eq:muchoice} implies $\mu_\ell = \mu \, 2^{-2K_\ell/d} > \const{C}{red} \, 2^{-2K_\ell/d}$ and $\mu_\ell - \const{C}{red} \, 2^{-2K_\ell/d} > \const{q}{red} \, 2^{-4K_\ell/d}$ for all $\ell \in \N_0$ since $K_\ell \in \N_0$ and hence $2^{-2K_\ell/d} \leq 1$. 
This, together with $\MM_{\ell,\nu} \subseteq \TT_{\ell,\nu} \setminus \TT_{\ell+1,\nu}$ for all $\nu \in \PPP_{\ell}$, $\MMM_{\ell} \subseteq \QQQ_{\ell} \cap \PPP_{\ell+1}$, and the  marking strategy~\eqref{eq:doerfler} shows that
\begin{equation} \label{eq:markapply}
 \begin{split}
  &\const{q}{red} \, 2^{-4K_{\ell + 1}/d} \, \tau_\ell((\TT_{\ell,\nu} \setminus \TT_{\ell + 1,\nu})_{\nu \in \PPP_{\ell}})^2 + (\mu_{\ell + 1} - \const{C}{red} \, 2^{-2K_{\ell + 1}/d}) \, \zeta_\ell(\QQQ_\ell \cap \PPP_{\ell + 1})^2 \\
  &\quad \geq  \const{q}{red} \, 2^{-4K_{\ell + 1}/d} \, \tau_\ell((\MM_{\ell,\nu})_{\nu \in \PPP_\ell})^2 + (\mu_{\ell + 1} - \const{C}{red} \, 2^{-2K_{\ell + 1}/d}) \, \zeta_\ell(\MMM_\ell)^2 \\
  &\quad \geq \const{q}{red} \, 2^{-4K_{\ell + 1}/d} \big[ \tau_\ell((\MM_{\ell,\nu})_{\nu \in \PPP_\ell})^2 + \zeta_\ell(\MMM_\ell)^2 \big] 
  \eqreff{eq:doerfler}{\geq} \theta \, \const{q}{red} \, 2^{-4K_{\ell + 1}/d} \, \eta_\ell^2.
 \end{split}
\end{equation} 
We set $\const{C}{stab} \coloneqq (1 + \eps^{-1}) \, (C_\tau^2 + \mu C_\zeta^2)$. We note that $2^{2K_\ell/d} > 1$ yields $1 + \eps_\ell^{-1} = 1 + \eps^{-1} \, 2^{2K_\ell/d} < 2^{2K_\ell/d} \, (1 + \eps^{-1})$. This, together with the definition of $\widetilde{\eta}_\ell$ from~\eqref{eq:etatildedef}, the estimates~\eqref{eq:taustabred}--\eqref{eq:zetastabred}, $2^{2K_\ell/d} > 1$ for all $\ell \in \N_0$, and estimate~\eqref{eq:markapply} verifies
 \begin{equation} \label{eq:etatildeesti}
 \begin{split}
  \widetilde{\eta}_{\ell+1}^2 
  &\leq (1 + \eps_{\ell+1}) \, \big[ 2^{-4K_{\ell + 1}/d} \, \tau_{\ell}^2 + \mu_{\ell+1} \, \zeta_{\ell}^2 + \mu_{\ell+1} \, \alpha_{\ell}(u_\ell^\star)^2   \\
  &\qquad \qquad \quad \ \  - 2^{-4K_{\ell + 1}/d} \, \const{q}{red} \, \tau_\ell((\TT_{\ell,\nu} \setminus \TT_{\ell + 1,\nu})_{\nu \in \PPP_{\ell}})^2 \\
  &\qquad \qquad \quad \ \ - (\mu_{\ell+1} - \const{C}{red} \, 2^{-2K_{\ell + 1}/d}) \, \zeta_\ell(\QQQ_\ell \cap \PPP_{\ell + 1})^2 \big] \\
  &\qquad \qquad \quad \ \ + (1 + \eps_{\ell + 1}^{-1}) \, (C_{\tau}^{2} \, 2^{-2K_{\ell + 1}/d }  + \mu \, C_\zeta^{2} \, 2^{-2K_{\ell+1}/d}) \, \enorm{u_{\ell+1}^\star - u_{\ell}^\star}^2 \\
  &\eqreff*{eq:markapply}{\leq} (1 + \eps_{\ell + 1}) \, \big[ 2^{-4K_{\ell + 1}/d}\tau_{\ell}^2 + \mu_{\ell+1} \, \zeta_{\ell}^2 +  \mu_{\ell + 1} \alpha_{\ell}(u_\ell^\star)^2 \\
  &\qquad \qquad \quad \ \ - \theta \, \const{q}{red} \, 2^{-4K_{\ell + 1}/d}\eta_\ell^2 \big] + \const{C}{stab} \, \enorm{u_{\ell + 1}^\star - u_\ell^\star}^2 \\
  &\leq (1 + \eps_{\ell + 1}) \, \big[ (1 - \theta \, \const{q}{red}) \, 2^{-4K_{\ell + 1}/d} \, \tau_\ell^2 \\
  &\qquad \qquad \quad \ \ + (1 -  \mu_{\ell + 1}^{-1} \, \theta \, \const{q}{red} \, 2^{-4K_{\ell + 1}/d}) \, \mu_{\ell+1} \, \zeta_\ell^2 \\
  &\qquad \qquad \quad \ \ + \mu_{\ell + 1} \, \alpha_{\ell}(u_\ell^\star)^2 \big] + \const{C}{stab} \, \enorm{u_{\ell + 1}^\star - u_\ell^\star}^2. \\
 \end{split}
 \end{equation}
 \textbf{Step~3 (Estimating \boldsymbol{$\mathsf{H}_{\ell + 1}^2$}).}
We start by estimating the factor in front of $2^{-4K_\ell/d} \,\tau_\ell(u_\ell^\star)^2$.
Due to $K_\ell \leq K_{\ell + 1}$, we obtain $2^{-4K_{\ell + 1}/d} \leq 2^{-4K_\ell/d}$ and $\mu_{\ell + 1} \leq \mu_\ell$ for all $\ell \in \N_0$. 
This, together with Galerkin orthogonality~\eqref{eq:weakformdisc}, estimate~\eqref{eq:etatildeesti}, and reliability~\eqref{eq:reliability}, yield
\begin{equation} \label{eq:Hesti}
 \begin{split}
  \mathsf{H}_{\ell + 1}^2 
  &= \enorm{u^\star - u_{\ell + 1}^\star}^2 + \widetilde{\eta}_{\ell + 1}^2
  \eqreff{eq:weakformdisc}{=} \enorm{u^\star - u_\ell^\star}^2 - \enorm{u_{\ell + 1}^\star - u_\ell^\star}^2 + \widetilde{\eta}_{\ell + 1}(u_{\ell + 1}^\star)^2 \\
  &\eqreff*{eq:reliability}{\leq} (1 - \rho_{\ell + 1}) \, \enorm{u^\star - u_\ell^\star}^2 + \rho_{\ell + 1} \, \const{C}{rel}^2 \, \eta_\ell^2 - \enorm{u_{\ell + 1}^\star - u_\ell^\star}^2 + \widetilde{\eta}_{\ell+1}^2 \\
  &\eqreff*{eq:etatildeesti}{\leq} (1 - \rho_{\ell+1}) \, \enorm{u^\star - u_\ell^\star}^2 \\
  &\qquad + \big[\rho_{\ell + 1}  \, \const{C}{rel}^2 \, 2^{4K_{\ell + 1}/d} + (1 + \eps_{\ell + 1})\, (1 - \theta \, \const{q}{red})\big] \, 2^{-4K_{\ell + 1}/d} \, \tau_{\ell}^2 \\
  &\qquad + \big[\rho_{\ell + 1} \, \mu_{\ell + 1}^{-1} \, \const{C}{rel}^2 + (1 + \eps_{\ell + 1})\, (1 - \mu_{\ell + 1}^{-1} \, \theta \, \const{q}{red} \, 2^{-4K_{\ell + 1}/d})\big] \, \mu_{\ell + 1} \, \zeta_\ell^2 \\
  &\qquad + (1 + \eps_{\ell + 1}) \, \mu_{\ell + 1} \, \alpha_{\ell}(u_\ell^\star)^2 + (\const{C}{stab} - 1) \, \enorm{u_{\ell + 1}^\star - u_\ell^\star}^2 \\
  &\leq (1 - \rho_{\ell+1}) \, \enorm{u^\star - u_\ell^\star}^2 \\
  &\qquad + \big[\rho_{\ell + 1}  \const{C}{rel}^2 \, 2^{4K_{\ell + 1}/d} + (1 + \eps_{\ell + 1})\, (1 - \theta \,  \const{q}{red})\big] \, 2^{-4K_\ell/d} \, \tau_{\ell}^2 \\
  &\qquad + \big[\rho_{\ell + 1} \, \mu_{\ell  + 1}^{-1} \, \const{C}{rel}^2 + (1 + \eps_{\ell + 1}) \, (1 - \mu_{\ell  + 1}^{-1} \, \theta \, \const{q}{red}2^{-4K_{\ell + 1}/d})\big] \, \mu_{\ell} \, \zeta_\ell^2 \\
  &\qquad + (1 + \eps_{\ell + 1}) \, \mu_{\ell  + 1} \, \alpha_{\ell}(u_\ell^\star)^2 + (\const{C}{stab} - 1) \, \enorm{u_{\ell + 1}^\star - u_\ell^\star}^2.
 \end{split}
\end{equation}\textbf{Step~4 (Proof of the quasi-contraction (\ref{eq:auxconv})).}
Recall $\delta$ from~\eqref{eq:deltadef}. We want to prove that the sum of the first three terms on the right-hand side of~\eqref{eq:Hesti} is bounded by $(1 - \delta \, 2^{-4K_{\ell + 1}/d}) \, \mathsf{H}_\ell^2$. 
We start by dealing with the factor in front of $\enorm{u^\star - u_\ell^\star}^2$. 
From~\eqref{eq:deltadef}, we get that $\delta \leq \rho$ and with~\eqref{eq:labbrev} we obtain $1 - \rho_{\ell + 1} = 1 - \rho \, 2^{-4K_{\ell + 1}/d} \leq 1 - \delta \, 2^{-4K_{\ell + 1}/d}$. 
To estimate the factor in front of $2^{-4K_\ell/d} \, \tau_\ell^2$, we use~\eqref{eq:labbrev},~\eqref{eq:deltadef}, and $0 < K_\ell \leq K_{\ell + 1}$ for all $\ell \in \N_0$ to obtain
\begin{equation*}
 \begin{split}
  &\rho_{\ell + 1} \, 2^{4K_{\ell + 1}/d} \, \const{C}{rel}^2 + (1 + \eps_{\ell + 1}) \, (1 - \theta \, \const{q}{red}) 
  \eqreff{eq:labbrev}{=} \rho \, \const{C}{rel}^2 + (1 + \eps \, 2^{-2K_{\ell + 1}/d}) \, (1 - \theta \, \const{q}{red}) \\
  &\qquad \leq 1 + \eps + \rho \, \const{C}{rel}^2 - (1 + \eps) \, \theta \, \const{q}{red} 
  \eqreff{eq:deltadef}{\leq} 1 - \delta \, 2^{-4K_0/d} 
  \leq 1 - \delta \, 2^{-4K_{\ell + 1}/d}.
 \end{split}
\end{equation*} 
Finally, we estimate the factor in front of $\mu_\ell \,\zeta_\ell^2$. We use~\eqref{eq:labbrev} and $K_\ell \geq 0$ to obtain
\begin{equation*}
 \begin{split}
   &\rho_{\ell + 1} \, \mu_{\ell + 1}^{-1} \, \const{C}{rel}^2 + (1 + \eps_{\ell + 1}) \, (1 - \mu_{\ell + 1}^{-1} \, \theta \, \const{q}{red} \, 2^{-4K_{\ell + 1}/d}) \\
   &\qquad \eqreff*{eq:labbrev}{=} \rho \, \mu^{-1} \, \const{C}{rel}^2 \, 2^{-2K_{\ell + 1}/d} + (1 + \eps \, 2^{-2K_{\ell + 1}/d}) \, (1 - \mu^{-1} \, \theta \, \const{q}{red} \, 2^{-2K_{\ell + 1}/d}) \\
   &\qquad = 1 - \mu^{-1} \, (\theta \, \const{q}{red}  - \mu \, \eps - \rho \, \const{C}{rel}^2) \, 2^{-2K_{\ell + 1}/d} - \eps \, \mu^{-1} \, \theta \, \const{q}{red} \, 2^{-4K_{\ell + 1}/d} \\
   &\qquad \leq 1 - \mu^{-1} \, (\theta \, \const{q}{red}  - \mu \, \eps - \rho \, \const{C}{rel}^2) \, 2^{-4K_{\ell + 1}/d} 
   \leq 1 - \delta \, 2^{-4K_{\ell + 1}/d}.
 \end{split}
\end{equation*}
We set $\const{C}{per} \coloneqq \max \set{(1 + \eps)\,\mu,\const{C}{stab} - 1}$ and obtain
\begin{equation*}
 \begin{split}
 &(1 + \eps_{\ell + 1}) \, \mu_{\ell + 1} \, \alpha_{\ell}(u_\ell^\star)^2 + (\const{C}{stab} - 1) \, \enorm{u_{\ell + 1}^\star - u_\ell^\star}^2 \\
 &\qquad = (1 + \eps 2^{-2K_{\ell + 1}/d}) \, \mu \, 2^{-2K_{\ell + 1}/d} \, \alpha_{\ell}(u_\ell^\star)^2 + (\const{C}{stab} - 1) \, \enorm{u_{\ell + 1}^\star - u_\ell^\star}^2 \\
 &\qquad \leq (1 + \eps) \, \mu \, \alpha_{\ell}(u_\ell^\star)^2 + (\const{C}{stab} - 1) \, \enorm{u_{\ell + 1}^\star - u_\ell^\star}^2 
 \leq \const{C}{per} \, (\alpha_{\ell}(u_\ell^\star)^2 + \enorm{u_{\ell + 1}^\star - u_\ell^\star}^2).
 \end{split}
\end{equation*} 
Hence, we showed that 
\begin{equation*}
 \mathsf{H}_{\ell + 1}^2 
 \leq (1 - \delta \, 2^{-4K_{\ell + 1}/d}) \, \mathsf{H}_\ell^2 + \const{C}{per} \, (\alpha_{\ell}(u_\ell^\star)^2 + \enorm{u_{\ell + 1}^\star - u_\ell^\star}^2),
\end{equation*}
which concludes the proof.
\end{proof}

Finally, the foregoing results allow us to prove Theorem~\ref{thm:main}.

\begin{proof}[\textbf{Proof of Theorem~\ref{thm:main}}]
The proof is split into two steps.

\textbf{Step~1 (Perturbation is summable).}
We want to prove that 
\begin{equation} \label{eq:perturb}
 \sum_{\ell = 0}^\infty \enorm{u_{\ell + 1}^\star - u_\ell^\star}^2 < +\infty \quad \text{and} \quad \sum_{\ell = 0}^\infty  \alpha_{\ell}(u_\ell^\star)^2 < +\infty .
\end{equation}
To show that the first sum in~\eqref{eq:perturb} converges, we use a Pythagorean identity to obtain, for all $N \in \N_0$,
\begin{equation*}
 \sum_{\ell = 0}^N \, \enorm{u_{\ell + 1}^\star - u_\ell^\star}^2 
 = \enorm{u_{N + 1}^\star - u_0^\star}^2 
 \leq\enorm{u^\star - u_{N + 1}^\star}^2 + \enorm{u_{N+ 1}^\star - u_0^\star}^2
 = \enorm{u^\star - u_0^\star}^2 
 < +\infty.
\end{equation*}It remains to estimate the second sum in~\eqref{eq:perturb}. We note that if there exists a constant $C > 0$ such that $M_\ell \leq C$ (cf. Section~\ref{subsec:discspaces}) for all $\ell \in \N_0$, then $\alpha_{\ell}(u_\ell^\star)^2$ vanishes for almost every $\ell \in \N_0$, which concludes that the sum of all $\alpha_{\ell}(u_\ell^\star)^2$ is finite.
We now consider the case that $M_\ell \to +\infty$ as $\ell \to \infty$. 
We first show that $\beta_m^0 \leq 1$ for all $m \in \N$.
Due to $P_m^0 = 1$ and $P_m^{-1} = 0$, the orthonormality condition~\eqref{eq:ortho} and the recurrence relation~\eqref{eq:recurrence} imply
\begin{equation} \label{eq:betabound}
 \begin{split}
  1 
  &\eqreff{eq:ortho}{=} \int_{[-1,1]} P_m^1(y_m)^2 \d \pi_m(y_m)
  \eqreff{eq:recurrence}{=} \int_{[-1,1]} (\beta_m^0)^{-2} y_m^2 P_m^0(y_m)^2 \d \pi_m(y_m) \\
  &= (\beta_m^0)^{-2} \int_{[-1,1]} y_m^2 \d \pi_m(y_m) \leq (\beta_m^0)^{-2} \pi_m([-1,1]) = (\beta_m^0)^{-2}.
 \end{split}
\end{equation}
This yields that $\beta_m^0 \leq 1$ for all $m \in \N$. Using the estimates~\eqref{eq:normequivX}--\eqref{eq:parseval}, and the continuous dependence on the data (cf.~\eqref{eq:weakformdisc}) by the Lax--Milgram lemma, we obtain
\begin{equation} \label{eq:uellbound}
 \norm{\nabla u_{\ell,\vec{0}}^{\star}}_\Omega^2 
 \eqreff{eq:normequivX}{\lesssim} \norm{u_{\ell,\vec{0}}^\star}_\X^2
 \leq \sum_{\nu \in \JJJ} \norm{u_{\ell,\nu}^\star}_\X^2
 \eqreff{eq:parseval}{=} \enorm{u_{\ell}^\star}_0^2
 \eqreff{eq:normequiv}{\lesssim} \enorm{u_{\ell}^\star}^2
 \eqreff{eq:weakformdisc}{\lesssim} \norm{f}_\Omega^2
 \quad \text{for all } \ell \in \N_0,
\end{equation}
where the hidden constants depend only on the diffusion coefficient $a$.
For $m \in \N$, let $\ell(m) \in \N_0$ be such that $M_{\ell(m)+1} + 1 = m$. The estimates~\eqref{eq:betabound} and~\eqref{eq:uellbound} as well as the assumptions on $a$ from~\eqref{eq:asumbound} and~\eqref{eq:ampliesti} imply 
\begin{equation*}
 \begin{split}
 \sum_{\ell \in \N_0} \alpha_{\ell}(u_\ell^\star)^2 
 &\leq \sum_{m \in \N} (\beta_m^0)^2 \norm{a_m}^2_{W^{1,\infty}(\Omega)} \norm{\nabla u_{\ell(m),\vec{0}}^\star}_\Omega^2 
 \lesssim \sum_{m \in \N} \norm{a_m}_{W^{1,\infty}(\Omega)}^2 \\
 &\leq \Bigl(\sum_{m \in \N} \norm{a_m}_{W^{1,\infty}(\Omega)} \Bigr)^2 \quad \, 
 \overset{\makebox[0pt]{\mbox{\normalfont\tiny\sffamily \eqref{eq:asumbound}+\eqref{eq:ampliesti}}}}{<} \quad \, +\infty,
 \end{split}
\end{equation*}
which concludes the proof of~\eqref{eq:perturb}.

\textbf{Step~2 (Convergence of error and estimator).}
Combining Lemma~\ref{lem:quasi} and Lemma~\ref{lem:sumcrit} with the summability~\eqref{eq:perturb} of the perturbations, we obtain that 
\begin{equation} \label{eq:Hlconv}
0 \leq \enorm{u^\star - u_\ell^\star}^2 \leq \mathsf{H}_\ell^2 \xrightarrow{\ell \to \infty} 0.
\end{equation}
If $K_\ell \leq K$ for all $\ell \in \N_0$, then $2^{2K_\ell/d} \leq 2^{4K_\ell/d} \leq 2^{4K/d}$ for all $\ell \in \N_0$. As in~\eqref{eq:Hlconv} and with the definition of $\widetilde{\eta}_\ell$ from~\eqref{eq:etatildedef}, we obtain
\begin{equation*}
 0 
 \leq \eta_\ell^2 
 \leq 2^{4K/d} \max \set{ 1 , \mu^{-1}} \, \widetilde{\eta}_\ell^2
 \lesssim \mathsf{H}_\ell^2
 \xrightarrow{\ell \to \infty} 0.
\end{equation*}
This concludes the proof.
\end{proof}

\section{Proof of Proposition~\ref{prop:stabred} (stability and reduction)} \label{sec:stabred}

\noindent A key difficulty in the analysis of Algorithm~\ref{alg:adapt} is the behavior of the spatial error indicator $\tau_\ell$ within the adaptive loop, since it involves interactions between piecewise polynomial functions with respect to different triangulations. 
The goal of this section is to prove Proposition~\ref{prop:stabred}, which is a main ingredient for the proof of Theorem~\ref{thm:main}.

The proof of Proposition~\ref{prop:stabred} requires the following result from~\cite[Lemma~3.1]{Eigel2014a}.
\begin{lemma} \label{lem:multiplic}
The operator $L_m \colon \ell^2(\JJJ) \rightarrow \ell^2(\JJJ)$ defined by
\begin{equation} \label{eq:Lmdef}
 L_m\bigl( (x_{\nu})_{\nu \in \JJJ} \bigr) 
 \coloneqq 
 (\beta^{\nu_m}_m x_{\nu + \eps_m} + \beta^{\nu_m -1}_m x_{\nu - \eps_m})_{\nu \in \JJJ},
\end{equation}
satisfies
\begin{equation} \label{eq:Kmnorm}
 \norm{L_m\bigl((x_\nu)_{\nu \in \JJJ}\bigr)}_{\ell^2(\JJJ)} 
 \leq \norm{(x_\nu)_{\nu \in \JJJ}}_{\ell^2(\JJJ)} 
 \quad \text{for all } m \in \N \text{ and all } (x_\nu)_{\nu \in \JJJ} \in \ell^2(\JJJ). \qed
\end{equation}\end{lemma}

Throughout, given $m \in \N$, $\nu \in \JJJ$, and $w = \displaystyle{\sum_{\mu \in \JJJ}} w_\mu P_\mu \in \V$, we will abbreviate 
\begin{subequations} \label{eq:Kmnot}
\begin{equation} 
 L_{m,\nu}(w) 
 \coloneqq \beta^{\nu_m}_m w_{\nu + \eps_m} + \beta^{\nu_m - 1}_m w_{\nu - \eps_m}, 
 \quad L_{m,\nu}(\norm{w}_\X) 
 \coloneqq \bigl(L_m((\norm{w_\nu}_\X)_{\nu \in \JJJ})\bigr)_\nu,
\end{equation}
and
\begin{equation} \label{eq:L0def}
 L_{0,\nu}(w) 
 \coloneqq w_\nu, 
 \quad L_{0,\nu}(\norm{w}_\X) 
 \coloneqq \norm{w_\nu}_\X.
\end{equation}
\end{subequations}
The quantity $L_{m,\nu}(\norm{\nabla w}_\Omega)$ is defined analogously.
Note that with the definition~\eqref{eq:Lmdef}, we can rewrite $\sigma_\nu$ from~\eqref{eq:sigmadef} as
\begin{equation*}
 \sigma_\nu(w) 
 = \sum_{m \in \N_0} a_m L_{m,\nu}(w).
\end{equation*}

The following lemma is key to proving~\eqref{eq:taustab} and~\eqref{eq:taured}. 
Its proof crucially relies on $K_\ell$-compatibility of the meshes $(\TT_{\ell,\nu})_{\nu \in \PPP_\ell}$, which allows us to compare local mesh sizes of neighboring meshes (see~\eqref{eq:defLcomp2}). 
Given $\nu \in \PPP_\ell$, $T \in \TT_{\ell,\nu}$, and $w_\ell \in \V_\ell$, we define an auxiliary spatial error indicator
\begin{equation} \label{eq:tautildedef}
 \widetilde{\tau}_\ell(\nu,T;w_\ell) 
 \coloneqq \Bigl( \sum_{\substack{T' \in \TT_\ell[\nu] \\ T' \subseteq T}} \big[h_T^2 \norm{\div \sigma_\nu(w_\ell)}_{T'}^2 + h_T\norm{\jump{\sigma_\nu(w_\ell)} \cdot \vec{n}_{\partial T'}}_{\partial T'\cap \Omega}^2\big] \Bigr)^{1/2}.
\end{equation}
The quantities $\widetilde{\tau}_\ell(\nu;w_\ell)$, $\widetilde{\tau}_\ell(w_\ell)$, and $\widetilde{\tau}_\ell$ are defined analogously to~\eqref{eq:estnot}. 
\begin{lemma} \label{lem:tech}
 Under the assumptions of Section~\ref{sec:apost}, there exists a constant $\widetilde{C}_\tau > 0$ depending only on $\Omega$, the polynomial degree $p$, the dimension $d$, and the initial mesh $\TT_0$ such that 
 \begin{equation} \label{eq:tech}
  \widetilde{\tau}_\ell(\nu;w_\ell)
  \leq \widetilde{C}_\tau \, 2^{K_\ell/d} \sum_{m \in \N_0} \norm{a_m}_{W^{1,\infty}(\Omega)} \, L_{m,
  \nu}(\norm{\nabla w_\ell}_\Omega) 
  \quad \text{for all } w_\ell \in \V_\ell \text{ and } \nu \in \PPP_\ell.
 \end{equation} 
 Moreover, there exists a constant $\const{C}{nrg} > 0$ depending only on the diffusion coefficient $a$ such that
 \begin{equation} \label{eq:stabnrg}
  \Bigl(\sum_{\nu \in \JJJ} \Bigl[\sum_{m \in \N_0} \norm{a_m}_{W^{1,\infty}(\Omega)} L_{m,\nu}(\norm{\nabla w_\ell}_\Omega) \Bigr]^2\Bigr)^{1/2} 
  \leq \const{C}{nrg} \, \enorm{w_\ell}
  \quad \text{for all } w_\ell \in \V_\ell.
 \end{equation}
\end{lemma}
\begin{proof}
Let $w_\ell \in \V_\ell$ and $\nu \in \PPP_\ell$. The proof is split into three steps.

\textbf{Step~1 (Estimate for the volume term).}
First, we show that there exists a constant $C_{\tau,1} > 0$ depending only on the polynomial degree $p$, and the dimension $d$ such that 
\begin{equation} \label{eq:tauesti1}
 \Bigl(\sum_{T \in \TT_{\ell,\nu}} \sum_{\substack{T' \in \TT_\ell[\nu] \\ T' \subseteq T}} h_T^2  \norm{ \div \sigma_\nu(w_\ell)}_{T'}^2\Bigr)^{1/2} 
 \leq C_{\tau,1}\, 2^{K_\ell/d}\sum_{m \in \N_0} \norm{a_m}_{W^{1,\infty}(\Omega)} L_{m,\nu}(\norm{\nabla w_\ell}_\Omega).
\end{equation}
Since $L_{m,\nu}(w_\ell)|_{T'} \in \P^p(T')$ for all $T' \in \TT_\ell[\nu]$, the product rule yields
\begin{equation} \label{eq:prodrule}
 \div \sigma_\nu(w_\ell) 
 = \sum_{m \in \N_0} \div[a_m \nabla L_{m,\nu}(w_\ell)] 
 = \sum_{m \in \N_0} [a_m \Delta L_{m,\nu}(w_\ell) + \nabla a_m \cdot \nabla L_{m,\nu}(w_\ell)]
\end{equation}
on $T'$ for all $T' \in \TT_\ell[\nu]$.
We note that, since the family $(\TT_{\ell,\nu})_{\nu \in \PPP_\ell}$ is $K_\ell$-compatible, the triangulations $\TT_{\ell,\nu}$ and $\TT_\ell[\nu]$ from~\eqref{eq:Toverlaydef} are also $K_\ell$-compatible. 
This establishes
\begin{equation} \label{eq:overlaycomp}
 h_T 
 \leq \const{C}{gen}^2 2^{K_\ell/d} h_{T'} 
 \quad \text{for all } T \in \TT_{\ell,\nu} \text{ and } T' \in \TT_\ell[\nu] \text{ with } T' \subseteq T.
\end{equation} 
Recall the inverse inequality
\begin{equation} \label{eq:invineqvol}
 \norm{\Delta \varphi_\ell}_{T'} 
 \leq \const{C}{inv}h_{T'}^{-1} \norm{\nabla \varphi_\ell}_{T'} 
 \quad \text{for all } T' \in \TT_\ell[\nu] \text{ and } \varphi_\ell \in \P^p(T'),
\end{equation}
where the constant $\const{C}{inv}$ depends only on the polynomial degree $p$ and the dimension $d$.
Altogether, we obtain
\begin{equation} \label{eq:invvolappl}
 \begin{split}
  &\sum_{T \in \TT_{\ell,\nu}} \sum_{\substack{T' \in \TT_\ell[\nu] \\ T' \subseteq T}} h_T^2 \norm{ a_m \Delta L_{m,\nu}(w_\ell)}_{T'}^2 
  \leq \norm{a_m}_{\Linf}^2 \sum_{T \in \TT_{\ell,\nu}} h_T^2 \sum_{\substack{T' \in \TT_\ell[\nu] \\ T' \subseteq T}} \norm{\Delta L_{m,\nu}(w_\ell)}_{T'}^2 \\
  &\qquad \eqreff*{eq:invineqvol}{\lesssim} \norm{a_m}_{\Linf}^2 \sum_{T \in \TT_{\ell,\nu}} h_T^2 \sum_{\substack{T' \in \TT_\ell[\nu] \\ T' \subseteq T}} h_{T'}^{-2} \norm{\nabla L_{m,\nu}(w_\ell)}_{T'}^2 \\
  &\qquad \eqreff*{eq:overlaycomp}{\lesssim} 2^{2K_\ell/d}\norm{a_m}_{\Linf}^2 \sum_{T \in \TT_{\ell,\nu}} \sum_{\substack{T' \in \TT_\ell[\nu] \\ T' \subseteq T}} \norm{\nabla L_{m,\nu}(w_\ell)}_{T'}^2 \\
  &\qquad = 2^{2K_\ell/d} \norm{a_m}_{\Linf}^2 \norm{\nabla L_{m,\nu}(w_\ell)}_\Omega^2.
 \end{split}
\end{equation}
From~\eqref{eq:prodrule}, we deduce with estimate~\eqref{eq:invvolappl}, $h_T = \abs{T}^{1/d} \leq \abs{\Omega}^{1/d}$, and $2^{K_\ell/d} \geq 1$ for all $\ell \in \N_0$ that
\begin{equation*}
 \begin{split}
  &\Bigl(\sum_{T \in \TT_{\ell,\nu}} \sum_{\substack{T' \in \TT_\ell[\nu] \\ T' \subseteq T}} h_T^2 \norm{\div \sigma_\nu(w_\ell)}_{T'}^2 \Bigr)^{1/2} 
  \leq \sum_{m \in \N_0} \Bigl(\sum_{T \in \TT_{\ell,\nu}} \sum_{\substack{T' \in \TT_\ell[\nu] \\ T' \subseteq T}} h_T^2 \norm{\div [a_m \nabla L_{m,\nu}(w_\ell)]}_{T'}^2 \Bigr)^{1/2} \\
  &\qquad \eqreff*{eq:prodrule}{\leq} \sum_{m \in \N_0} \Bigl(\sum_{T \in \TT_{\ell,\nu}} \sum_{\substack{T' \in \TT_\ell[\nu] \\ T' \subseteq T}}  h_T^2\big[ \norm{a_m \Delta L_{m,\nu}(w_\ell)}_{T'} + \norm{\nabla a_m \cdot \nabla L_{m,\nu}(w_\ell)}_{T'}\big]^2\Bigr)^{1/2} \\
  &\qquad \eqreff*{eq:invvolappl}{\lesssim} 2^{K_\ell/d}\sum_{m \in \N_0} \norm{a_m}_{\Linf} \norm{\nabla L_{m,\nu}(w_\ell)}_\Omega + \sum_{m \in \N_0} \norm{\nabla a_m}_{\Linf} \norm{\nabla L_{m,\nu}(w_\ell)}_{\Omega} \\
  &\qquad \leq  2^{K_\ell/d}\sum_{m \in \N_0} \norm{a_m}_{W^{1,\infty}(\Omega)} \norm{\nabla L_{m,\nu}(w_\ell)}_{\Omega}.
 \end{split}
\end{equation*}
Noting that $\norm{\nabla L_{m,\nu}(w_\ell)}_\Omega \leq L_{m,\nu}(\norm{\nabla w_{\ell}}_\Omega)$ concludes the proof of~\eqref{eq:tauesti1}.

\textbf{Step~2 (Estimate for the jump term).}
We want to show that there exists a constant $C_{\tau,2} > 0$ depending only on the polynomial degree $p$, and the dimension $d$ such that
\begin{equation} \label{eq:jumpesti}
 \begin{split}
 \Bigl(\sum_{T \in \TT_{\ell,\nu}} \sum_{\substack{T' \in \TT_\ell[\nu] \\ T' \subseteq T}} h_T &\norm{\jump{\sigma_\nu(w_\ell)} \cdot \vec{n}_{\partial T'}}_{\partial T' \cap \Omega}^2 \Bigr)^{1/2} \\
 &\qquad \leq C_{\tau,2} \, 2^{K_\ell/(2d)} \sum_{m \in \N_0} \norm{a_m}_{W^{1,\infty}(\Omega)} L_{m,\nu}(\norm{\nabla w_\ell}_\Omega).
 \end{split}
\end{equation}
To this end, we introduce the patch 
\begin{equation*}
 \omega_\ell[\nu,T'] 
 \coloneqq \bigcup \set{T'' \in \TT_\ell[\nu] \given T'' \cap T' \neq \emptyset},
\end{equation*}
i.e., the domain generated by $T'$ and all elements $T'' \in \TT_\ell[\nu]$ that touch $T'$.
It is well known that there is a constant $\const{C}{patch} > 0$ depending only on the initial mesh $\TT_0$ such that 
\begin{equation} \label{eq:patch}
 \sum_{T' \in \TT_\ell[\nu]} \norm{\varphi}_{\omega_\ell[\nu,T']}^2 \leq \const{C}{patch}^2 \norm{\varphi}_\Omega^2 
 \quad \text{for all } \varphi \in L^2(\Omega).
\end{equation}

Moreover, recall the discrete trace inequality
\begin{equation} \label{eq:invineqtrace}
 \norm{\varphi_\ell}_{\partial T'} 
 \leq \const{C}{trace} \, h_{T'}^{-1/2} \norm{\varphi_\ell}_{\omega_\ell[\nu,T']} 
 \quad \text{for all } T' \in \TT_\ell[\nu] \text{ and } \varphi_\ell \in \SS^{p-1}_0(\TT_\ell[\nu]),
\end{equation}
where the constant $\const{C}{trace}$ depends only on the polynomial degree $p$, the dimension $d$, and $\kappa$-shape regularity. 
Given $T' \in \TT_\ell[\nu]$, Morrey's inequality applied to $a_m \in W^{1,\infty}(\Omega)$ (see, e.g.,~\cite[Theorem~5.4]{Adams2003}) and the discrete trace inequality yield
\begin{equation} \label{eq:invineqtraceappl}
 \begin{split}
  \norm{\sigma_\nu(w_\ell)}_{\partial T' \cap \Omega} 
  &= \Big \lVert \sum_{m \in \N_0} a_m \cdot \nabla L_{m,\nu}(w_\ell) \Big \rVert_{\partial T' \cap \Omega} \\
  &\leq \sum_{m \in \N_0} \norm{a_m}_{L^\infty(\partial T')} \norm{\nabla L_{m,\nu}(w_\ell)}_{\partial T' \cap \Omega} \\
  &\leq \sum_{m \in \N_0} \norm{a_m}_{L^\infty(\Omega)} \norm{\nabla L_{m,\nu}(w_\ell)}_{\partial T'} \\
  &\eqreff*{eq:invineqtrace}{\lesssim} h_{T'}^{-1/2}\sum_{m \in \N_0} \norm{a_m}_{W^{1,\infty}(\Omega)} \norm{\nabla L_{m,\nu}(w_\ell)}_{\omega_\ell[\nu,T']}.
 \end{split}
\end{equation}
Given a facet $F$ of $\TT_\ell[\nu]$ and elements $T_1,T_2 \in\TT_\ell[\nu]$ with $T_1 \cap T_2 = F$, the jump $\jump{\sigma_\nu(w_\ell)}|_F$ is just the difference of the traces of $\sigma_\nu(w_\ell)|_{T_1}$ and $\sigma_\nu(w_\ell)|_{T_2}$ on $F$. Hence, $\abs{\vec{n}_{\partial T'}} = 1$, together with the estimate~\eqref{eq:invineqtraceappl}, $K_\ell$-compatibility~\eqref{eq:overlaycomp}, estimate~\eqref{eq:patch}, and the triangle inequality yield
\begin{equation*}
 \begin{split}
  &\Bigl(\sum_{T \in \TT_{\ell,\nu}} \sum_{\substack{T' \in \TT_\ell[\nu] \\ T' \subseteq T}} h_T \norm{\jump{\sigma_\nu(w_\ell)} \cdot \vec{n}_{\partial T'}}_{\partial T'\cap \Omega}^2\Bigr)^{1/2} 
  \leq \Bigl(\sum_{T \in \TT_{\ell,\nu}} \sum_{\substack{T' \in \TT_\ell[\nu] \\ T' \subseteq T}} h_T \norm{\jump{\sigma_\nu(w_\ell)}}_{\partial T'\cap \Omega}^2\Bigr)^{1/2} \\ 
  &\qquad \lesssim \Bigl(\sum_{T \in \TT_{\ell,\nu}} \sum_{\substack{T' \in \TT_\ell[\nu] \\ T' \subseteq T}} h_T \norm{\sigma_\nu(w_\ell)}_{\partial T' \cap \Omega}^2 \Bigr)^{1/2}
  = \Bigl(\sum_{T \in \TT_{\ell,\nu}} h_T \sum_{\substack{T' \in \TT_\ell[\nu] \\ T' \subseteq T}} \norm{\sigma_\nu(w_\ell)}_{\partial T' \cap \Omega}^2\Bigr)^{1/2} \\
  &\qquad \eqreff*{eq:invineqtraceappl}{\lesssim} \Bigl(\sum_{T \in \TT_{\ell,\nu}} h_T \sum_{\substack{T' \in \TT_\ell[\nu] \\ T' \subseteq T}} h_{T'}^{-1} \Bigl[\sum_{m \in \N_0} \norm{a_m}_{W^{1,\infty}(\Omega)} \norm{\nabla L_{m,\nu}(w_\ell)}_{\omega_\ell[\nu,T']}\Bigr]^2\Bigr)^{1/2}  \\
  &\qquad \eqreff*{eq:overlaycomp}{\lesssim} 2^{K_\ell/(2d)} \Bigl(\sum_{T \in \TT_\ell[\nu]} \Bigl[\sum_{m \in \N_0} \norm{a_m}_{W^{1,\infty}(\Omega)} \norm{\nabla L_{m,\nu}(w_\ell)}_{\omega_\ell[\nu,T]}\Bigr]^2\Bigr)^{1/2} \\
  &\qquad \leq 2^{K_\ell/(2d)} \sum_{m \in \N_0} \norm{a_m}_{W^{1,\infty}(\Omega)} \Bigl(\sum_{T \in \TT_\ell[\nu]} \norm{\nabla L_{m,\nu}(w_\ell)}_{\omega_\ell[\nu,T]}^2\Bigr)^{1/2} \\
  &\qquad \eqreff*{eq:patch}{\lesssim} 2^{K_\ell/(2d)} \sum_{m \in \N_0} \norm{a_m}_{W^{1,\infty}(\Omega)} \norm{\nabla L_{m,\nu}(w_\ell)}_{\Omega} \\
  &\qquad \leq 2^{K_\ell/(2d)} \sum_{m \in \N_0} \norm{a_m}_{W^{1,\infty}(\Omega)} \, L_{m,\nu}(\norm{\nabla w_\ell}_\Omega).
 \end{split}
\end{equation*}
This concludes the proof of~\eqref{eq:jumpesti}.
From the estimates~\eqref{eq:tauesti1} and~\eqref{eq:jumpesti} as well as $2^{K_\ell/(2d)} \leq 2^{K_\ell/d}$ for all $\ell \in \N_0$, we conclude~\eqref{eq:tech}.

\textbf{Step~3 (Proof of~(\ref{eq:stabnrg})).}
The triangle inequality with respect to the norm $\norm{\, \cdot \,}_{\ell^2(\JJJ)}$, the estimate~\eqref{eq:Kmnorm} from Lemma~\ref{lem:multiplic} as well as the properties~\eqref{eq:asumbound} and~\eqref{eq:ampliesti} of $a$ yield
\begin{equation} \label{eq:fintautildestab}
 \begin{split}
  &\Bigl(\sum_{\nu \in \JJJ} \Bigl[ \sum_{m \in \N_0} \norm{a_m}_{W^{1,\infty}(\Omega)} L_{m,\nu}(\norm{\nabla w_\ell}_\Omega) \Bigr]^2 \Bigr)^{1/2} \\
  &\qquad \leq \sum_{m \in \N_0} \Bigl(\sum_{\nu \in \JJJ} \norm{a_m}_{W^{1,\infty}(\Omega)}^2 L_{m,\nu}(\norm{\nabla w_\ell}_\Omega)^2 \Bigr)^{1/2} \\
  &\qquad = \sum_{m \in \N_0} \norm{a_m}_{W^{1,\infty}(\Omega)} \big \lVert L_m\bigl((\norm{\nabla w_{\ell,\nu}}_\Omega)_{\nu \in \JJJ}\bigr)\big \rVert_{\ell^2(\JJJ)} \\
  &\qquad \eqreff*{eq:Kmnorm}{\leq} \sum_{m \in \N_0} \norm{a_m}_{W^{1,\infty}(\Omega)} \big \lVert (\norm{\nabla w_{\ell,\nu}}_\Omega)_{\nu \in \JJJ}\big \rVert_{\ell^2(\JJJ)}
  \qquad 
 \overset{\makebox[0pt]{\mbox{\normalfont\tiny\sffamily \eqref{eq:asumbound} + \eqref{eq:ampliesti}}}}{\leq} \qquad \big \lVert (\norm{\nabla w_{\ell,\nu}}_\Omega)_{\nu \in \JJJ}\big \rVert_{\ell^2(\JJJ)}.
 \end{split}
\end{equation}
By the norm equivalences~\eqref{eq:normequivX}--\eqref{eq:normequiv}, and the Parseval-type identity~\eqref{eq:parseval}, we conclude
\begin{equation} \label{eq:sumnormequiv}
 \begin{split}
  \norm{(\norm{\nabla w_{\ell,\nu}}_\Omega)_{\nu \in \JJJ}}_{\ell^2(\JJJ)} 
  &\eqreff{eq:normequivX}{\lesssim} \big\lVert(\norm{w_{\ell,\nu}}_\X)_{\nu \in \JJJ} \big \rVert_{\ell^2(\JJJ)} 
  \eqreff{eq:parseval}{=} \enorm{w_\ell}_0
  \eqreff{eq:normequiv}{\lesssim} \enorm{w_\ell},
 \end{split}
\end{equation}%
where the hidden constants depend only on the diffusion coefficient $a$.
Combining the results~\eqref{eq:fintautildestab} and~\eqref{eq:sumnormequiv} yields estimate~\eqref{eq:stabnrg}. 
This concludes the proof.
\end{proof}
We now prove the statements of Proposition~\ref{prop:stabred} consecutively, starting with~\eqref{eq:taustab}. 

\begin{proof}[\textbf{Proof of stability~(\ref{eq:taustab}).}]
The proof is split into two steps.

\textbf{Step~1 (Estimate for \boldsymbol{$\widetilde{\tau}_\ell$}).}
First, we want to show that
\begin{equation} \label{eq:tautildeesti}
 \widetilde{\tau}_\ell(w_\ell) \leq \widetilde{C}_\tau \const{C}{nrg} \, 2^{K_\ell/d}\,\enorm{w_\ell} \quad \text{for all } w_\ell \in \V_\ell,
\end{equation}%
where $\widetilde{C}_\tau, \const{C}{nrg} > 0$ are the constants from Lemma~\ref{lem:tech}.
Let $w_\ell \in \V_\ell$. 
We apply Lemma~\ref{lem:tech} to obtain
\begin{equation*}
  \begin{split}
   \widetilde{\tau}_\ell(w_\ell)^2 
   &= \sum_{\nu \in \PPP_\ell} \widetilde{\tau}_\ell(\nu;w_\ell)^2
   \eqreff{eq:tech}{\leq} \widetilde{C}_\tau^2 2^{2K_\ell/d} \sum_{\nu \in \PPP_\ell} \Bigl[ \sum_{m \in \N_0} \norm{a_m}_{W^{1,\infty}(\Omega)} L_{m,\nu}(\norm{\nabla w_\ell}_\Omega) \Bigr]^2  \\
   &\leq \widetilde{C}_\tau^2 2^{2K_\ell/d} \sum_{\nu \in \JJJ} \Bigl[ \sum_{m \in \N_0} \norm{a_m}_{W^{1,\infty}(\Omega)} L_{m,\nu}(\norm{\nabla w_\ell}_\Omega) \Bigr]^2
   \eqreff{eq:stabnrg}{\leq} \widetilde{C}_\tau^2 \const{C}{nrg}^2 2^{2K_\ell/d} \, \enorm{w_\ell}^2.
 \end{split}
\end{equation*}%
This concludes the proof of~\eqref{eq:tautildeesti}.

\textbf{Step~2 (Stability of \boldsymbol{$\tau_\ell$}).}
Direct calculations show the triangle-type estimate
\begin{equation} \label{eq:elemesti}
 (a + b)^2 + (c + d)^2 \leq [(a^2 + c^2)^{1/2} + (b^2 + d^2)^{1/2}]^2 
 \quad \text{for all } a,b,c,d \in \R.
\end{equation}
Let $w_\ell, v_\ell \in \V_\ell$ and $\nu \in \PPP_\ell$.
The triangle inequality, $\sigma_\nu(w_\ell) - \sigma_\nu(v_\ell) = \sigma_\nu(w_\ell - v_\ell)$, and~\eqref{eq:elemesti} yield
\begin{align*}
 &{\tau_\ell(\nu;w_\ell)}^2 
 \quad \, 
 \overset{\makebox[0pt]{\mbox{\normalfont\tiny\sffamily \eqref{eq:spatialdef}+\eqref{eq:estnot}}}}{=} \quad \, \sum_{T \in \TT_{\ell,\nu}} \sum_{\substack{T' \in \TT_\ell[\nu] \\ T' \subseteq T}}  \bigl[h_T^2 \norm{(\delta_{\vec{0} \nu} f + \div \sigma_\nu(w_\ell))}_{T'}^2 + h_T\norm{\jump{\sigma_\nu(w_\ell)} \cdot \vec{n}_{\partial T'}}_{\partial T' \cap \Omega}^2\bigr] \\
 &\quad \leq \sum_{T \in \TT_{\ell,\nu}} \sum_{\substack{T' \in \TT_\ell[\nu] \\ T' \subseteq T}}\bigl[ h_T^2(\norm{(\delta_{\vec{0} \nu} f + \div \sigma_\nu(v_\ell))}_{T'} + \norm{\div \sigma_\nu(w_\ell - v_\ell)}_{T'})^2 \\
 &\qquad \qquad \qquad  + h_T(\norm{\jump{\sigma_\nu(v_\ell)} \cdot \vec{n}_{\partial T'}}_{\partial T' \cap \Omega} + \norm{\jump{\sigma_\nu(w_\ell - v_\ell)} \cdot \vec{n}_{\partial T'}}_{\partial T' \cap \Omega})^2 \bigr]  \\
 &\quad \leq \Bigl( \Bigl[ \sum_{T \in \TT_{\ell,\nu}} \sum_{\substack{T' \in \TT_\ell[\nu] \\ T' \subseteq T}} h_T^2 \norm{(\delta_{\vec{0} \nu} f + \div \sigma_\nu(v_\ell))}_{T'}^2 \Bigr]^{1/2} \\
 &\qquad \qquad \qquad + \Bigl[ \sum_{T \in \TT_{\ell,\nu}} \sum_{\substack{T' \in \TT_\ell[\nu] \\ T' \subseteq T}} h_T^2 \norm{\div \sigma_\nu(w_\ell - v_\ell)}_{T'}^2\Bigr]^{1/2} \Bigr)^2 \\
 &\qquad \qquad \qquad  + \Bigl(\Bigl[\sum_{T \in \TT_{\ell,\nu}} \sum_{\substack{T' \in \TT_\ell[\nu] \\ T' \subseteq T}} h_T \norm{\jump{\sigma_\nu(v_\ell)} \cdot \vec{n}_{\partial T'}}_{\partial T' \cap \Omega}^2\Bigr]^{1/2} \\
 &\qquad \qquad \qquad + \Bigl[\sum_{T \in \TT_{\ell,\nu}} \sum_{\substack{T' \in \TT_\ell[\nu] \\ T' \subseteq T}} h_T \norm{\jump{\sigma_\nu(w_\ell - v_\ell)} \cdot \vec{n}_{\partial T'}}_{\partial T' \cap \Omega}^2\Bigr]^{1/2}\Bigr)^2 \\
 &\quad \eqreff{eq:elemesti}{\leq} \Bigl( \Bigl[ \sum_{T \in \TT_{\ell,\nu}} \sum_{\substack{T' \in \TT_\ell[\nu] \\ T' \subseteq T}} \bigl(h_T^2 \norm{(\delta_{\vec{0} \nu} f + \div \sigma_\nu(v_\ell))}_{T'}^2 + h_T \norm{\jump{\sigma_\nu(v_\ell)} \cdot \vec{n}_{\partial T'}}_{\partial T' \cap \Omega}^2\bigr)\Bigr]^{1/2} \\
 &\quad \qquad + \Bigl[\sum_{T \in \TT_{\ell,\nu}} \sum_{\substack{T' \in \TT_\ell[\nu] \\ T' \subseteq T}} \bigl(h_T^2\norm{\div \sigma_\nu(w_\ell - v_\ell)}_{T'}^2 + h_T\norm{\jump{\sigma_\nu(w_\ell - v_\ell)} \cdot \vec{n}_{\partial T'}}_{\partial T' \cap \Omega}^2\bigr)\Bigr]^{1/2}\Bigr)^2 \\
 &\quad = \Bigl(\tau_\ell(\nu;v_\ell) + \widetilde{\tau}_\ell(\nu;w_\ell - v_\ell)\Bigr)^2.
\end{align*}
This leads to
\begin{equation*}
 \begin{split}
  \tau_\ell(w_\ell) 
  &= \Bigl(\sum_{\nu \in \PPP_\ell} \tau_\ell(\nu;w_\ell)^2 \Bigr)^{1/2} \leq \Bigl( \sum_{\nu \in \PPP_\ell} [\tau_\ell(\nu;v_\ell) + \widetilde{\tau}_\ell(\nu;w_\ell - v_\ell)]^2 \Bigr)^{1/2} \\
  &\leq \Bigl(\sum_{\nu \in \PPP_\ell} \tau_\ell(\nu;v_\ell)^2 \Bigr)^{1/2} + \Bigl( \sum_{\nu \in \PPP_\ell} \widetilde{\tau}_\ell(\nu;w_\ell - v_\ell)^2\Bigr)^{1/2} \\
  &= \tau_\ell(v_\ell) + \widetilde{\tau}_\ell(w_\ell - v_\ell)
  \eqreff{eq:tautildeesti}{\leq} \tau_\ell(v_\ell) + C_\tau \, 2^{K_\ell/d} \,\enorm{w_\ell - v_\ell},
 \end{split}
\end{equation*}
where $C_\tau = \widetilde{C}_\tau \const{C}{nrg}$.
This concludes the proof of~\eqref{eq:taustab}.
\end{proof}

\begin{proof}[\textbf{Proof of stability~(\ref{eq:zetastab}).}]
The proof is split into two steps.

\textbf{Step~1.} We show that
\begin{equation} \label{eq:zetaest}
 \zeta_\ell(w_\ell) 
 \leq \const{C}{nrg} \, \enorm{w_\ell} 
 \quad \text{for all } w_\ell \in \V_\ell,
\end{equation}
where $\const{C}{nrg} > 0$ is the constant from Lemma~\ref{lem:tech} in~\eqref{eq:stabnrg}.
We apply Lemma~\ref{lem:tech} to obtain
\begin{equation*}
 \begin{split}
  \zeta_\ell(w_\ell)
  &\eqreff*{eq:paramdef}{=} \Bigl(\sum_{\nu \in \QQQ_\ell} \Bigl[ \sum_{m \in \N} \norm{a_m}_{W^{1,\infty}(\Omega)} L_{m,\nu}(\norm{\nabla w_\ell}_\Omega) \Bigr]^2 \Bigr)^{1/2} \\
  &\leq \Bigl(\sum_{\nu \in \JJJ} \Bigl[ \sum_{m \in \N} \norm{a_m}_{W^{1,\infty}(\Omega)} L_{m,\nu}(\norm{\nabla w_\ell}_\Omega) \Bigr]^2 \Bigr)^{1/2}
  \eqreff{eq:stabnrg}{\leq} \const{C}{nrg} \, \enorm{w_\ell},
 \end{split}
\end{equation*}
which concludes the proof of~\eqref{eq:zetaest}.

\textbf{Step~2 (Stability of \boldsymbol{$\zeta_\ell$}).}
Let $w_\ell, v_\ell \in \V_\ell$. 
The triangle inequality and~\eqref{eq:Kmnot} yield
\begin{equation*}
 \begin{split}
  L_{m,\nu}(\norm{\nabla w_\ell}_\Omega) 
  \leq L_{m,\nu}(\norm{\nabla v_\ell}_{\Omega}) + L_{m,\nu}(\norm{\nabla (w_\ell - v_\ell)}_{\Omega}) 
  \quad \text{for all } m \in \N_0 \text{ and } \nu \in \QQQ_\ell,
 \end{split}
\end{equation*}
which immediately yields
\begin{equation} \label{eq:zetatriangle}
 \zeta_\ell(\nu;w_\ell) 
 \leq \zeta_\ell(\nu;v_\ell) + \zeta_\ell(\nu;w_\ell - v_\ell)
 \quad \text{for all } \nu \in \QQQ_\ell.
\end{equation}
Combining the triangle inequality with the estimates~\eqref{eq:zetaest}--\eqref{eq:zetatriangle} establishes
\begin{equation*}
 \begin{split}
  \zeta_\ell(w_\ell) 
  &\eqreff*{eq:zetatriangle}{\leq} \Bigl(\sum_{\nu \in \QQQ_\ell} [\zeta_\ell(\nu;v_\ell) + \zeta_\ell(\nu;w_\ell - v_\ell)]^2 \Bigr)^{1/2} \\
  &\leq \Bigl( \sum_{\nu \in \QQQ_\ell} \zeta_\ell(\nu;v_\ell)^2 \Bigr)^{1/2} + \Bigl( \sum_{\nu \in \QQQ_\ell} \zeta_\ell(\nu;w_\ell - v_\ell)^2 \Bigr)^{1/2} \\
  &= \zeta_\ell(v_\ell) + \zeta_\ell(w_\ell - v_\ell) 
  \eqreff{eq:zetaest}{\leq} \zeta_\ell(v_\ell) + \const{C}{nrg} \enorm{w_\ell - v_\ell}.
 \end{split}
\end{equation*}
This concludes the proof of~\eqref{eq:zetastab} with $C_\zeta = \const{C}{nrg}$.
\end{proof}
\begin{proof}[\textbf{Proof of reduction~(\ref{eq:taured}).}]
Let $w_\ell \in \V_\ell$. There holds 
\begin{equation} \label{eq:tausplit}
 \tau_{\ell + 1}(w_\ell)^2 
 = \tau_{\ell + 1}(\PPP_\ell;w_\ell)^2 + \tau_{\ell + 1}(\PPP_{\ell + 1} \setminus \PPP_\ell;w_\ell)^2.
\end{equation}
The remainder of the proof is split into three steps.

\textbf{Step~1 (Estimate for \boldsymbol{$\tau_{\ell + 1}(\PPP_\ell;w_\ell)$}).} 
We want to show that 
\begin{equation} \label{eq:FEMred}
 \tau_{\ell + 1}(\PPP_\ell;w_\ell)^2 
 \leq \tau_\ell(w_\ell)^2 - (1 - 2^{-1/d}) \, \tau_\ell((\TT_{\ell,\nu} \setminus \TT_{\ell + 1,\nu})_{\nu \in \PPP_\ell};w_\ell)^2.
\end{equation}
Let $\nu \in \PPP_\ell$. Since $w_\ell \in \V_\ell$, the assumptions on $(a_m)_{m \in \N_0}$ imply that
\begin{equation*}
 \sigma_\nu(w_\ell) 
 \in H(\div,T') 
 \coloneqq \set{\vec{\sigma} \in L^2(T') \given \div \vec{\sigma} \in L^2(T')}
 \quad \text{for all } T' \in \TT_\ell[\nu].
\end{equation*}
Hence, for any $\widetilde{T} \in \TT_{\ell,\nu}$, $\TT_{\ell + 1}[\nu] \preceq \TT_\ell[\nu]$ yields
\begin{equation} \label{eq:divsum}
 \sum_{\substack{T \in \TT_{\ell+1,\nu} \\ T \subseteq \widetilde{T}}} \sum_{\substack{T' \in \TT_{\ell + 1}[\nu] \\ T' \subseteq T}} \norm{\delta_{\vec{0} \nu} f + \div \sigma_\nu(w_\ell)}_{T'}^2 
 = \sum_{\substack{T' \in \TT_\ell[\nu] \\ T' \subseteq \widetilde{T}}} \norm{\delta_{\vec{0} \nu} f + \div \sigma_\nu(w_\ell)}_{T'}^2.
\end{equation}
By similar arguments, $\jump{\sigma_\nu(w_\ell)}$ vanishes on the interior faces of $\widetilde{T}$, leading to
\begin{equation} \label{eq:jumpsum}
 \sum_{\substack{T \in \TT_{\ell + 1,\nu} \\ T \subseteq \widetilde{T}}} \sum_{\substack{T' \in \TT_{\ell + 1}[\nu] \\ T' \subseteq T}} \norm{\jump{\sigma_\nu(w_\ell)} \cdot \vec{n}_{\partial T'}}_{\partial T' \cap \Omega}^2 
 = \sum_{\substack{T' \in \TT_\ell[\nu] \\ T' \subseteq \widetilde{T}}} \norm{ \jump{\sigma_\nu(w_\ell)} \cdot \vec{n}_{\partial T'}}_{\partial T' \cap \Omega}^2.
\end{equation}
Using~\eqref{eq:divsum}--\eqref{eq:jumpsum} for $\widetilde{T} \in \TT_{\ell,\nu} \cap \TT_{\ell + 1,\nu}$, we obtain
\begin{equation} \label{eq:tausum}
 \tau_{\ell + 1}(\nu;w_\ell)^2 = \tau_\ell(\nu,\TT_{\ell,\nu} \cap \TT_{\ell + 1,\nu};w_\ell)^2 + \tau_{\ell + 1}(\nu,\TT_{\ell + 1,\nu} \setminus \TT_{\ell,\nu};w_\ell)^2.
\end{equation}
Refinement by bisection guarantees
\begin{equation} \label{eq:sizered}
 h_T 
 = \abs{T}^{1/d} 
 \leq (\abs{\widetilde{T}}/2)^{1/d} 
 = 2^{-1/d} \, h_{\widetilde{T}} 
\end{equation}
for all $\widetilde{T} \in \TT_{\ell,\nu} \setminus \TT_{\ell + 1,\nu}$ and all $T \in \TT_{\ell + 1,\nu} \setminus \TT_{\ell,\nu}$ with $T \subsetneqq \widetilde{T}$. Arguing as in~\eqref{eq:divsum}--\eqref{eq:jumpsum} for $\widetilde{T} \in \TT_{\ell,\nu} \setminus \TT_{\ell + 1,\nu}$, inequality~\eqref{eq:sizered} yields
\begin{equation} \label{eq:femred}
 \begin{split}
  \tau_{\ell + 1}(\nu,\TT_{\ell + 1,\nu} \setminus \TT_{\ell,\nu};w_\ell)^2  
  &= \sum_{\widetilde{T} \in \TT_{\ell,\nu} \setminus \TT_{\ell + 1,\nu}} \sum_{\substack{T \in \TT_{\ell + 1,\nu}\setminus \TT_{\ell,\nu} \\ T \subset \widetilde{T}}} \tau_{\ell + 1}(\nu,T;w_\ell)^2 \\
  &\leq 2^{-1/d} \, \tau_\ell(\nu,\TT_{\ell,\nu} \setminus \TT_{\ell + 1,\nu};w_\ell)^2. \\
 \end{split}
\end{equation}
This, together with~\eqref{eq:tausum}, leads to
\begin{equation*}
 \begin{split}
  \tau_{\ell + 1}(\PPP_\ell;w_\ell)^2 
  &\eqreff{eq:tausum}{=} \tau_\ell((\TT_{\ell+1,\nu} \cap \TT_{\ell,\nu})_{\nu \in \PPP_\ell};w_\ell)^2 + \tau_{\ell + 1}((\TT_{\ell + 1,\nu} \setminus \TT_{\ell,\nu})_{\nu \in \PPP_\ell};w_\ell)^2 \\
  &\eqreff{eq:femred}{\leq} \tau_\ell(w_\ell)^2 - (1 - 2^{-1/d}) \, \tau_\ell((\TT_{\ell,\nu} \setminus \TT_{\ell + 1,\nu})_{\nu \in \PPP_\ell};w_\ell)^2.
 \end{split}
\end{equation*}
This concludes the proof of~\eqref{eq:FEMred}.

\textbf{Step~2 (Estimate for \boldsymbol{$\tau_{\ell + 1}(\PPP_{\ell + 1} \setminus \PPP_\ell;w_\ell)$}).} We want to show that
\begin{equation} \label{eq:tau2zeta}
\tau_{\ell + 1}(\PPP_{\ell + 1} \setminus \PPP_\ell;w_\ell)^2
 \leq \widetilde{C}_\tau^2 \, 2^{2K_{\ell+ 1}/d} \, \zeta_\ell(\QQQ_\ell \cap \PPP_{\ell + 1};w_\ell)^2,
\end{equation}
where $\widetilde{C}_\tau > 0$ is the constant from Lemma~\ref{lem:tech}.
Since $\vec{0} \notin \PPP_{\ell + 1} \setminus \PPP_\ell$, there holds $\tau_{\ell + 1}(\nu;w_\ell) = \widetilde{\tau}_{\ell + 1}(\nu;w_\ell)$ (see~\eqref{eq:tautildedef}) as well as $L_{0,\nu}(w_\ell) = w_{\ell,\nu} = 0$ (see~\eqref{eq:L0def}) for all $\nu \in \PPP_{\ell + 1} \setminus \PPP_\ell$. Hence, we may apply Lemma~\ref{lem:tech} to obtain
\begin{equation} \label{eq:tau2zetaproof}
 \begin{split}
  \tau_{\ell + 1}(\PPP_{\ell + 1} \setminus \PPP_\ell;w_\ell)^2 
  &= \sum_{\nu \in \PPP_{\ell + 1} \setminus \PPP_\ell} \tau_{\ell + 1}(\nu;w_\ell)^2
  = \sum_{\nu \in \PPP_{\ell + 1} \setminus \PPP_\ell} \widetilde{\tau}_{\ell + 1}(\nu;w_\ell)^2 \\
  &\eqreff*{eq:tech}{\leq} \, \widetilde{C}_\tau^2 \, 2^{2K_{\ell + 1}/d} \sum_{\nu \in \PPP_{\ell + 1} \setminus \PPP_\ell} \Bigl[ \sum_{m \in \N_0} \norm{a_m}_{W^{1,\infty}(\Omega)} L_{m,\nu}(\norm{\nabla w_\ell}_\Omega) \Bigr]^2 \\
  &\eqreff*{eq:paramdef}{=} \, \widetilde{C}_\tau^2 \, 2^{2K_{\ell + 1}/d}  \sum_{\nu \in \PPP_{\ell + 1} \setminus \PPP_\ell} \zeta_\ell(\nu;w_\ell)^2 \\
  &= \widetilde{C}_\tau^2  2^{2K_{\ell + 1}/d}  \zeta_\ell(\PPP_{\ell + 1} \setminus \PPP_\ell;w_\ell)^2.
 \end{split}
\end{equation}
With $\PPP_{\ell + 1} \setminus \PPP_\ell \subseteq \QQQ_\ell \cap \PPP_{\ell + 1}$ and~\eqref{eq:tau2zetaproof}, we conclude~\eqref{eq:tau2zeta}.

\textbf{Step~3.} Combining~\eqref{eq:tausplit}, \eqref{eq:FEMred}, and~\eqref{eq:tau2zeta} yields reduction~\eqref{eq:taured} with $\const{C}{red} = \widetilde{C}_\tau^2$ and $\const{q}{red} = 1 - 2^{-1/d}$. 
This concludes the proof.
\end{proof}

\begin{proof}[\textbf{Proof of reduction~(\ref{eq:zetared})}]
Let $w_\ell \in \V_\ell$. By~\eqref{eq:detaildef}, there holds $\QQQ_\ell \setminus \QQQ_{\ell + 1} = \QQQ_\ell \cap \PPP_{\ell + 1}$. This yields
\begin{equation} \label{eq:seteq}
 \QQQ_\ell \cap \QQQ_{\ell + 1} = \QQQ_\ell \setminus (\QQQ_\ell \setminus \QQQ_{\ell + 1}) = \QQQ_\ell \setminus (\QQQ_\ell \cap \PPP_{\ell + 1}),
\end{equation}
which, together with $\zeta_{\ell + 1}(\nu;w_\ell) = \zeta_\ell(\nu;w_\ell)$ for all $\nu \in \QQQ_\ell$, leads to
\begin{equation} \label{eq:zetaeq}
 \begin{split}
  \zeta_{\ell + 1}(w_\ell)^2 
  &= \zeta_{\ell + 1}(\QQQ_{\ell + 1} \cap \QQQ_\ell;w_\ell)^2 + \zeta_{\ell + 1}(\QQQ_{\ell + 1} \setminus \QQQ_\ell;w_\ell)^2 \\
  &\eqreff*{eq:seteq}{=} \, \zeta_{\ell + 1}(\QQQ_\ell \setminus (\QQQ_\ell \cap \PPP_{\ell + 1});w_\ell)^2 + \zeta_{\ell + 1}(\QQQ_{\ell + 1} \setminus \QQQ_\ell;w_\ell)^2  \\
  &= \zeta_{\ell}(w_\ell)^2 - \zeta_\ell(\QQQ_\ell \cap \PPP_{\ell + 1};w_\ell)^2 + \zeta_{\ell + 1}(\QQQ_{\ell + 1} \setminus \QQQ_\ell;w_\ell)^2.
 \end{split}		
\end{equation}
It remains to estimate the last term on the right-hand side of~\eqref{eq:zetaeq}. We note that $w_{\ell,\nu} = 0$ for all $\nu \notin \PPP_\ell$. Thus, if $M_\ell = M_{\ell + 1}$, there holds $\zeta_{\ell + 1}(\nu;w_\ell) = 0$ for all $\nu \in \QQQ_{\ell + 1} \setminus \QQQ_\ell$. If $M_{\ell + 1} = M_\ell + 1$, then 
\begin{equation*}
 \zeta_{\ell + 1}(\QQQ_{\ell + 1} \setminus \QQQ_\ell;w_\ell)^2 
 = \zeta_{\ell + 1}(\eps_{M_{\ell + 1} + 1};w_\ell) ^2 
 \eqreff{eq:paramdef}{=} (\beta_{M_{\ell + 1} + 1}^0 \norm{a_{M_{\ell + 1} + 1}}_{W^{1,\infty}(\Omega)} \norm{\nabla w_{\ell,\vec{0}}}_\Omega)^2.
\end{equation*}
Using the definition of $\alpha_{\ell}(w_\ell)$ in~\eqref{eq:alphadef}, we obtain $\zeta_{\ell + 1}(\QQQ_{\ell + 1} \setminus \QQQ_\ell;w_\ell)^2 = \alpha_{\ell}(w_\ell)^2$, which, together with~\eqref{eq:zetaeq}, yields~\eqref{eq:zetared}.
This concludes the proof.
\end{proof}

\section{Multilevel refinement algorithm} \label{sec:algo}

Due to the restrictive nature of the mesh compatibility condition~\eqref{eq:defLcomp}, one has to modify the algorithms concerning multimesh refinement and parametric enrichment included in \cite{Bespalov2021,Bespalov2021a}. 
The interdependence of all involved meshes and the choice of triangulations associated with newly added multi-indices are the key difficulties in the implementation.

\subsection{Mesh refinement}
We now restrict ourselves to the case $d=2$ and recall the notion of reference edges (see Section~\ref{subsec:meshes}). 
The present section introduces an NVB-based multilevel refinement algorithm, which preserves the compatibility condition~\eqref{eq:defLcomp} for all meshes in the family. 
The algorithm is based on the edge-based NVB-algorithm from \cite{Karkulik2012}.

All of the subsequent algorithms require the compatibility degree to be a multiple of $2$, i.e., contained in $2\N$.
Hence, for a sequence of compatibility degrees $(K_\ell)_{\ell \in \N_0}$, we assume (without loss of generality) that $K_\ell \in 2\N$ for all $\ell \in \N_0$.
As a first algorithm, we present a multilevel refinement strategy, which generates a $K_{\ell + 1}$-compatible family of meshes $(\TT_{\ell + 1,\nu})_{\nu \in \PPP_{\ell}}$ from a given $K_\ell$-compatible family of meshes $(\TT_{\ell,\nu})_{\nu \in \PPP_\ell}$.
\begin{algorithm}[multilevel refinement] \label{alg:refine}
 \textbf{Input:} Finite index set $\PPP_\ell \subset \JJJ$, compatibility degrees $K_\ell, K_{\ell + 1} \in 2\N$, $K_\ell$-compatible family of meshes $(\TT_{\ell,\nu})_{\nu \in \PPP_\ell}$, and sets of marked elements $(\MM_{\ell,\nu})_{\nu \in \PPP_\ell}$ with $\MM_{\ell,\nu} \subseteq \TT_{\ell,\nu}$ for all $\nu \in \PPP_\ell$.
 \begin{enumerate}[label = \rm(\roman*)]
  \item Initialize $\MM_\nu^0 \coloneqq \set{E_T \given T \in \MM_{\ell,\nu}}$ for all $\nu \in \PPP_\ell$.
  \item \rm{\textbf{For}} all $k = 0,1,2,\ldots,$ repeat the following steps \rm{(a)--(c)}:
  \begin{enumerate}[label = \rm(\alph*)]
   \item \rm{\textbf{For}} all $\nu \in \PPP_\ell$, set
   \begin{equation*}
    \begin{split}
     \UU_{\nu}^k 
     \coloneqq \set{E_T \given T \in \TT_{\ell,\nu} \colon \EE_T \cap \MM_\nu^k \neq \emptyset}.
	\end{split}
   \end{equation*}
   \rm{\textbf{For}} all $\mu \in \PPP_\ell[\nu]$, set
   \begin{equation*}
	 \VV_{\nu,\mu}^k 
	 \coloneqq \{E_T \, | \, T \in \TT_{\ell,\nu} \colon \exists \, T' \in \TT_{\ell,\mu} \colon E_{T'} \in \MM_{\mu}^k, T' \subsetneqq T, \gen{T'} = \gen{T} + K_{\ell + 1}\}.
   \end{equation*}
   \item \rm{\textbf{For}} all $\nu \in \PPP_\ell$ and $\mu \in \PPP_\ell[\nu]$ set $\MM_\nu^{k + 1} \coloneqq \MM_\nu^k \cup \UU_{\nu}^k \cup \VV_{\nu,\mu}^k$.
   \item \rm{\textbf{If}} $\MM_\nu^{k + 1} = \MM_\nu^k$ for all $\nu \in \PPP_\ell$, \rm{\textbf{stop}} the loop.
  \end{enumerate}
  \item For every $\nu \in \PPP_\ell$, generate a refinement $\TT_{\ell + 1,\nu}$ of $\TT_{\ell,\nu}$ by bisecting every edge in $\MM_\nu^{k + 1}$.
 \end{enumerate}
 \textbf{Output:} $K_{\ell + 1}$-compatible family of triangulations $(\TT_{\ell + 1,\nu})_{\nu \in \PPP_\ell}$, where for all $\nu \in \PPP_\ell$, $\TT_{\ell + 1,\nu}$ is a refinement of $\TT_{\ell,\nu}$ such that at least all elements in $\MM_{\ell,\nu}$ are bisected. 
\end{algorithm}

The following result states that Algorithm~\ref{alg:refine} terminates after a finite number of iterations and generates the desired output. Its proof can be found in Appendix~\ref{sec:append}.
\begin{proposition} \label{prop:refine}
Given the specified input, Algorithm~\ref{alg:refine} terminates and generates a $K_{\ell + 1}$-compatible family of meshes $(\TT_{\ell + 1,\nu})_{\nu \in \PPP_{\ell}}$, where at least all elements in $\MM_{\ell,\nu} \subseteq \TT_{\ell,\nu}$ have been refined for all $\nu \in \PPP_\ell$.
\end{proposition}

\subsection{Parametric enrichment}
A key difficulty in the implementation is the choice of suitable meshes associated with newly added multi-indices, since the new meshes need to be $K_\ell$-compatible with the existing ones. 
In the following, we propose a strategy to associate suitable meshes to newly added multi-indices, which preserves the $K_\ell$-compatibility condition~\eqref{eq:defLcomp}. 
Throughout the remainder of this paper, we suppose that $\PPP_\ell$ is \emph{monotone}, i.e.,
\begin{equation*}
 \set{\nu - \eps_m \given \nu \in \PPP_\ell, m \in \supp{\nu}} \subseteq \PPP_\ell.
\end{equation*}
We note that monotonicity also implies that for all $\nu \in \PPP_\ell$ and $(b_m)_{m \in \supp{\nu}} \subset \N_0$ with $b_m \leq \nu_m$ for all $m \in \supp{\nu}$, there holds
\begin{equation*}
 \nu - \sum_{m \in \supp{\nu}} b_m \eps_m \in \PPP_\ell.
\end{equation*}
Given a set of marked multi-indices $\MMM_\ell \subseteq \QQQ_\ell$, we define
\begin{equation} \label{eq:markdef}
 \widetilde{\MMM}_\ell 
 \coloneqq \set{\nu \in \JJJ \setminus \PPP_\ell\given \exists \, \mu \in \MMM_\ell \colon \supp{\nu} \subseteq \supp{\mu} \text{ and } \nu_m \leq \mu_m \ \forall \, m \in \supp{\nu}},
\end{equation}
and note that the monotonicity of $\PPP_\ell$ implies that
\begin{equation} \label{eq:monotoneQ}
 \begin{split}
  \mathfrak{Q}_\ell 
  = \big\{ \nu \in \JJJ \setminus \PPP_\ell \; \vert \; &\exists \, m \in \set{1,2,\ldots,M_\ell} \text{ such that } \nu - \eps_m \in \PPP_\ell \big\} \cup \set{\eps_{M_\ell + 1}}.
 \end{split}
\end{equation}
It is easy to see that $\widetilde{\MMM}_\ell \subseteq \QQQ_\ell$. Indeed, let $\nu \in \widetilde{\MMM_\ell}$ and $\mu \in \MMM_\ell$ with $\supp{\nu} \subseteq \supp{\mu}$ and $\nu_m \leq \mu_m$ for all $m \in \supp{\nu}$. This yields the existence of a set $\II \subseteq \supp{\mu}$ and natural numbers $(b_m)_{m \in \II}$ such that
\begin{equation*}
 \nu + \sum_{m \in \II} b_m \eps_m = \mu.
\end{equation*} 
Since $\mu \in \MMM_\ell \subseteq \QQQ_\ell$, there is $m_0 \in \supp{\mu}$ such that %
\begin{equation*}
 \nu + \sum_{m \in \II} b_m \eps_m - \eps_{m_0} = \mu - \eps_{m_0} \in \PPP_\ell.
\end{equation*}
If $m_0 \in \II$, $b_{m_0} \geq 1$ yields that
\begin{equation*}
 \nu = \Big( \nu + \sum_{m \in \II} (b_m - \delta_{m m_0}) \eps_m \Big) - \sum_{m \in \II} (b_m - \delta_{m m_0}) \eps_m \in \PPP_\ell,
\end{equation*}
which contradicts the fact that $\nu \in \widetilde{\MMM}_\ell$. Hence, $m_0 \notin \II$ and
\begin{equation*}
 \nu - \eps_{m_0} = \Bigl( \nu + \sum_{m \in \II} b_m \eps_m - \eps_{m_0} \Bigr) - \sum_{m \in \II} b_m \eps_m = \mu - \eps_{m_0} - \sum_{m \in \II} b_m \eps_m \in \PPP_\ell,
\end{equation*}
yielding $\nu \in \QQQ_\ell$ by~\eqref{eq:monotoneQ}.

\begin{remark} \label{rem:monotone}
 Monotonicity of $\PPP_\ell$ is required to ensure that Algorithm~\ref{algo:compilize} below produces the desired output (cf. Proposition~\ref{prop:compilize} and its proof in Appendix~\ref{sec:append}) and has been considered in earlier works (see, e.g.,~\cite{Eigel2014a,Chkifa2014}). Additionally, the monotonicity of $\PPP_\ell$ is desirable from an implementation point of view, since it allows for a simple construction of the set $\QQQ_\ell$ (see~\eqref{eq:monotoneQ}) and imposes a tree structure on the index set $\PPP_\ell$ while retaining the approximation properties of the corresponding truncated polynomial chaos expansion as observed in~\cite[Section~2]{Chkifa2013}.
\end{remark}
Given a triangulation $\TT$ and a set of marked elements $\MM \subseteq \TT$, we denote by $\mathtt{NVB}(\TT,\MM)$ the usual $2d$-NVB bisection routine (see, e.g.,~\cite{Stevenson2007,Karkulik2012}).

Algorithm~\ref{algo:compilize} will later be called with the set $\RRR \coloneqq \PPP_\ell[\nu]$ consisting of all neighboring multi-indices of a newly added multi-index $\nu \in \QQQ_\ell \cap \PPP_{\ell + 1}$, to which a mesh has already been assigned.
\begin{algorithm}[Compatibilization]\label{algo:compilize}
 \textbf{Input$\colon$} Initial mesh $\TT_0$, compatibility degree $K_\ell \in 2\N$, finite index set $\mathfrak{R} \subseteq \PPP_\ell$, and pairwise $2K_\ell$-compatible family of meshes $(\TT_{\ell,\nu})_{\nu \in \RRR}$.
 \begin{enumerate}[label = \rm(\roman*)]
 \item Define the overlay mesh
 \begin{equation*}
  \TT_\infty \coloneqq \bigoplus_{\nu \in \RRR} \TT_{\ell,\nu}.
 \end{equation*}

 \item \rm{\textbf{For}} all $k = 0,1,\dots,$ repeat the following steps \rm{(a)--(c)} starting with the initial mesh $\TT_0$:
\begin{enumerate}[label = \rm(\alph*)]
\item Set
\begin{equation} \label{eq:compmarked}
 \MM_k 
 \coloneqq \set{T \in \TT_k \given \exists \, T' \in \TT_\infty \colon T' \subsetneqq T \text{ and } \gen{T'} > \gen{T} + K_\ell}.	
\end{equation}
\item \rm{\textbf{If}} $\MM_k = \emptyset$, \rm{\textbf{stop}} the loop and set $\underline{k} \coloneqq k$.
\item Apply the usual NVB-refinement algorithm to $\TT_k$ and $\MM_k$ to obtain the refinement $\TT_{k + 1} \coloneqq \mathtt{NVB}(\TT_k,\MM_k)$ such that all elements in $\MM_k$ are bisected.
\end{enumerate}
\end{enumerate}
\textbf{Output:} Mesh $\TT_{\underline{k}}$ such that $\TT_{\underline{k}}$ and $\TT_{\ell,\nu}$ are $K_\ell$-compatible for all $\nu \in \RRR$.
\end{algorithm}
\begin{figure}[!ht]
  \resizebox{\textwidth}{!}{
   \subfloat{
    \stackunder{{\Huge $\TT_1$}}{\includegraphics{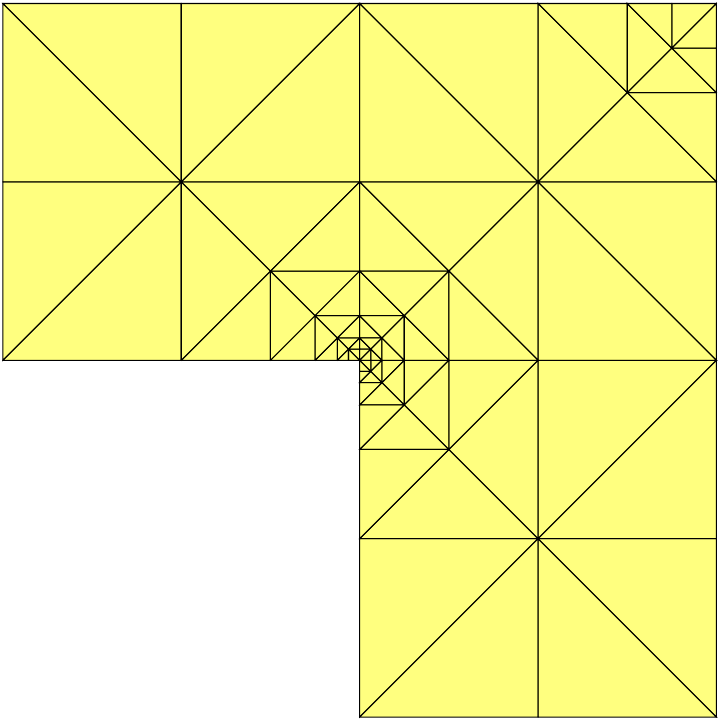}}
   }
   \subfloat{
    \stackunder{{\Huge $\TT_2$}}{\includegraphics{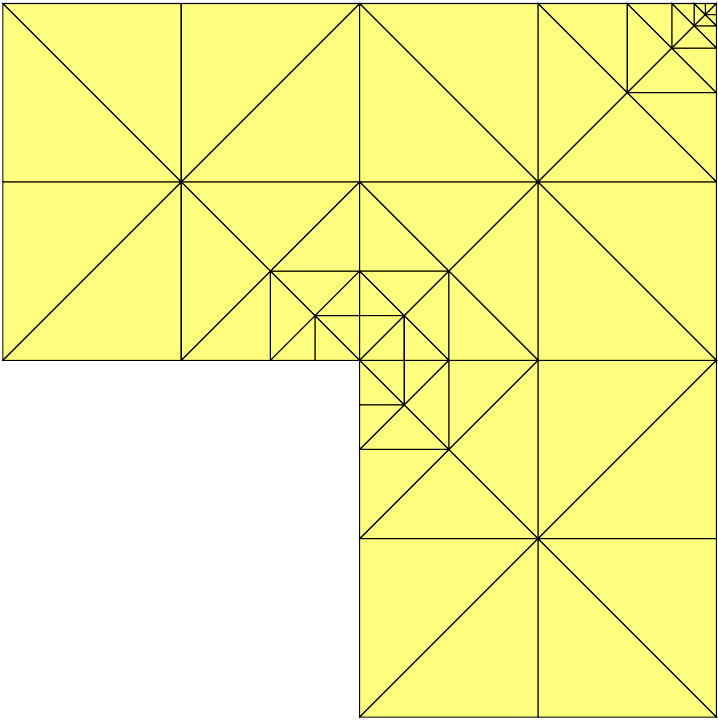}}
   }

   \subfloat{
    \stackunder{{\Huge $\mathtt{Compatibilize}(\TT_1,\TT_2,K = 2)$}}{\includegraphics{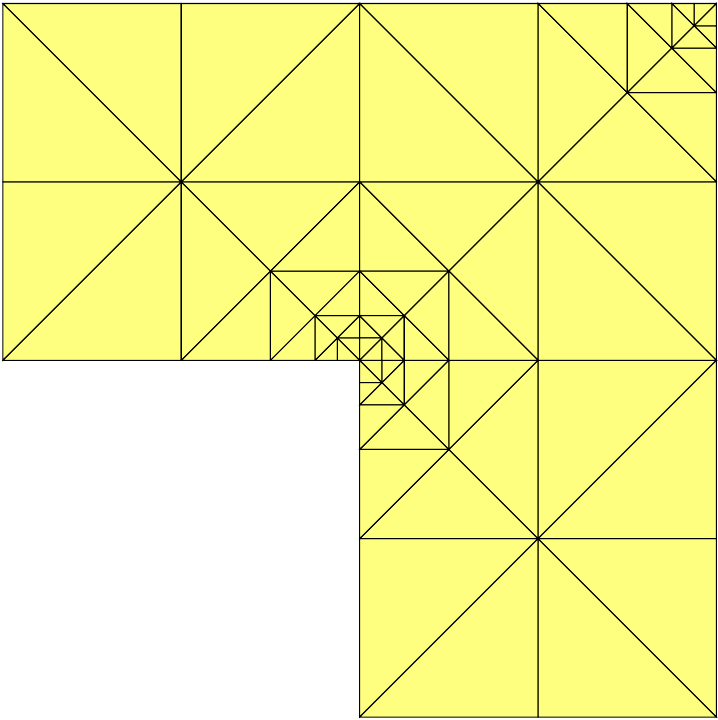}}
   }
  }
   \caption{Illustration of the output of Algorithm~\ref{algo:compilize}~(right) with $2K$-compatible input meshes $\TT_1$~(left), $\TT_2$~(center) and $K = 2$.}

   \label{fig:meshes}
 \end{figure}
An illustration of the output of Algorithm~\ref{algo:compilize} is given in Figure~\ref{fig:meshes}.
The following result states that Algorithm~\ref{algo:compilize} terminates after a finite number of iterations and generates the correct output. 
Its proof can be found in Appendix~\ref{sec:append}.
\begin{proposition} \label{prop:compilize}
Given the specified input, Algorithm~\ref{algo:compilize} terminates and $\TT_{\underline{k}}$ and $\TT_{\ell,\nu}$ are $K_\ell$-compatible for all $\nu \in \RRR$.
\end{proposition}
Now we are in a position to propose a multilevel spatial refinement and parametric enrichment algorithm, which can be used in the adaptive Algorithm~\ref{alg:adapt}.
\begin{algorithm}[multilevel refinement and enrichment] \label{alg:multilevel}
 \textbf{Input:} Finite index set $\PPP_\ell \subset \JJJ$, compatibility degrees $K_\ell, K_{\ell + 1} \in 2\N$, $K_\ell$-compatible family of meshes $(\TT_{\ell,\nu})_{\nu \in \PPP_\ell}$, a family of sets of marked elements $(\MM_{\ell,\nu})_{\nu \in \PPP_\ell}$ with $\MM_{\ell,\nu} \subseteq \TT_{\ell,\nu}$ for all $\nu \in \PPP_\ell$, and a set of marked multi-indices $\MMM_\ell \subseteq \QQQ_\ell$.
 \begin{enumerate}[label = \rm(\roman*)]
 \item Call Algorithm~\ref{alg:refine} with $(\TT_{\ell,\nu})_{\nu \in \PPP_\ell}$ and $(\MM_{\ell,\nu})_{\nu \in \PPP_\ell}$ to obtain a $K_{\ell + 1}$-compatible family of meshes $(\TT_{\ell + 1,\nu})_{\nu \in \PPP_\ell}$.
 \item With $\widetilde{\MMM}_\ell$ defined in~\eqref{eq:markdef}, set $\PPP_{\ell + 1} \coloneqq \PPP_\ell \cup \widetilde{\MMM}_\ell$ and construct the set $\QQQ_{\ell + 1}$.
 \item Sort $\widetilde{\MMM}_\ell = \set{\nu^1, \ldots , \nu^m}$ with respect to $\abs{\cdot}$, i.e., $\abs{\nu^{j-1}} \leq \abs{\nu^j}$ for all $j = 2, \ldots , m$.
 \item \rm{\textbf{For}} $j = 1, \ldots , m$ repeat the following steps \rm{(a)--(b)}:
 \begin{enumerate}[label = \rm(\alph*)]
 \item Set $\RRR_j \coloneqq \PPP_{\ell + 1}[\nu^j] \setminus \set{\nu^{j + 1}, \ldots ,\nu^m}$.
 \item Call Algorithm~\ref{algo:compilize} with $K_{\ell + 1}$, $\RRR \coloneqq \RRR_j$ and $(\TT_{\ell+1,\nu})_{\nu \in \RRR_j}$ to obtain $\TT_{\ell + 1,\nu^j} \coloneqq \TT_{\underline{k}}$.

 \end{enumerate} 
\end{enumerate}
\textbf{Output:} Sets $\PPP_{\ell + 1}, \QQQ_{\ell + 1}$, and a $K_{\ell + 1}$-compatible family of meshes $(\TT_{\ell + 1,\nu})_{\nu \in \PPP_{\ell + 1}}$.
\end{algorithm}
The following theorem states that Algorithm~\ref{alg:multilevel} generates the desired output, i.e., the output satisfies the properties of the routine $\mathtt{refine}$ (see~\eqref{eq:refine}) in the adaptive Algorithm~\ref{alg:adapt}.
\begin{theorem}
 Given the specified input, Algorithm~\ref{alg:multilevel} generates a family of meshes $(\TT_{\ell + 1,\nu})_{\nu \in \PPP_{\ell+1}}$ that is $K_{\ell+1}$-compatible and at least all elements in $\MM_{\ell,\nu} \subseteq \TT_{\ell,\nu}$ have been refined in $\TT_{\ell + 1,\nu}$ for all $\nu \in \PPP_\ell$.
\end{theorem}
\begin{proof}
The proof is split into two steps.

\textbf{Step~1.}
Since $(\TT_{\ell,\nu})_{\nu \in \PPP_\ell}$ is $K_\ell$-compatible, Algorithm~\ref{alg:multilevel}(i) and Proposition~\ref{prop:refine} yield that $(\TT_{\ell + 1,\nu})_{\nu \in \PPP_\ell}$ is $K_{\ell + 1}$-compatible. 

\textbf{Step~2.} We show that $(\TT_{\ell + 1,\nu})_{\nu \in \PPP_{\ell + 1}}$ is $K_{\ell + 1}$-compatible. For $j \in \set{1, \ldots , m}$, we set $\SSS_j \coloneqq \{\nu^1, \ldots , \nu^j\}$. Since $\PPP_{\ell + 1} = \PPP_\ell \cup \SSS_m$, it suffices to prove that $(\TT_{\ell + 1,\nu})_{\nu \in \PPP_\ell \cup \SSS_j}$ is $K_{\ell + 1}$-compatible for all $j = 1, \ldots , m$. We proceed by induction over $j$.

\textbf{Base case (\boldsymbol{$j = 1$}).} Since $\abs{\nu^1} = \min \set{\abs{\nu} \given \nu \in \widetilde{\MMM}_\ell}$, we obtain that $\set{\nu \in \PPP_{\ell + 1} \given \abs{\nu} < \abs{\nu^1}} \subseteq \PPP_\ell$. 
Hence, $(\TT_{\ell + 1,\nu})_{\nu\in \PPP_{\ell + 1}, \abs{\nu} < \abs{\nu^1}}$ is $K_{\ell + 1}$-compatible and it suffices to prove that $\TT_{\ell + 1,\nu^1}$ and $\TT_{\ell + 1,\nu}$ are $K_{\ell + 1}$-compatible for all $\nu \in \PPP_\ell[\nu^1] = \RRR_1$. 
Due to Proposition~\ref{prop:compilize}, it remains to show that the family $(\TT_{\ell + 1,\nu})_{\nu \in \RRR_1}$ is pairwise $2K_{\ell + 1}$-compatible. 
First we note that in the case of $\abs{\RRR_1} = 1$, the family $(\TT_{\ell + 1,\nu})_{\nu \in \RRR_1}$ is trivially pairwise $2K_{\ell + 1}$-compatible. 
We now consider the case $\abs{\RRR_1} > 1$. Let $\mu^1,\mu^2 \in \RRR_1 \subseteq \PPP_\ell$ with $\mu^1 \neq \mu^2$. 
By~\eqref{eq:monotoneQ}, there exist $m_1,m_2 \in \supp{\nu^1}$ such that $\mu^1 = \nu^1 - \eps_{m_1}$ and $\mu^2 = \nu^1 - \eps_{m_2}$.
Hence, $\mu \coloneqq \nu^1 - \eps_{m_1} - \eps_{m_2}$ is a common neighbor of $\mu^1$ and $\mu^2$ and, due to monotonicity of $\PPP_\ell$, also contained in $\PPP_\ell$. 
Since $\TT_{\ell + 1,\mu^1}$ and $\TT_{\ell + 1, \mu}$ as well as $\TT_{\ell + 1,\mu^2}$ and $\TT_{\ell + 1,\mu}$ are $K_{\ell + 1}$-compatible, $\TT_{\ell + 1,\mu^1}$ and $\TT_{\ell + 1,\mu^2}$ are $2K_{\ell + 1}$-compatible. 
Since $\mu^1,\mu^2 \in \RRR_1$ were arbitrary, we obtain that the meshes $(\TT_{\ell + 1,\mu})_{\mu \in \RRR_1}$ are pairwise $2K_{\ell + 1}$-compatible. This concludes the proof of $K_{\ell + 1}$-compatibility of $(\TT_{\ell + 1,\nu})_{\nu \in \PPP_\ell \cup \SSS_1}$.

\textbf{Induction hypothesis.} Let $j \in \set{1, \ldots , m-1}$ and assume that $(\TT_{\ell + 1,\nu})_{\nu \in \PPP_\ell \cup \SSS_j}$ is $K_{\ell + 1}$-compatible. 

\textbf{Induction step.} We want to prove that $(\TT_{\ell + 1,\nu})_{\nu \in \PPP_\ell \cup \SSS_{j + 1}}$ is $K_{\ell + 1}$-compatible. Noting that, due to the induction hypothesis, $(\TT_{\ell + 1,\nu})_{\nu \in \PPP_{\ell + 1}, \abs{\nu} < \abs{\nu^{j+1}}}$ is $K_{\ell + 1}$-compatible, analogous arguments for $\nu^{j+1}$ instead of $\nu^1$ as in the base case yield that also $(\TT_{\ell + 1,\nu})_{\nu \in \PPP_\ell \cup \SSS_{j+1}}$ is $K_{\ell + 1}$-compatible.
Hence, we proved that $(\TT_{\ell + 1,\nu})_{\nu \in \PPP_\ell \cup \SSS_j}$ is $K_{\ell + 1}$-compatible for all $j = 1, \ldots, m$. This concludes the proof.
\end{proof}

\section{Numerical experiments} \label{sec:numerics}
This section presents numerical experiments in $2$D that illustrate the performance of the proposed adaptive algorithm. We conducted all experiments using the software package T-IFISS from~\cite{Bespalov2019} for the lowest-order SGFEM with $p = 1$. Throughout this section, we compare the performance of Algorithm~\ref{alg:adapt} for uniform ($\theta = 1$) and adaptive ($0 < \theta < 1$) multilevel refinement as well as its performance for different choices of compatibility degrees $(K_\ell)_{\ell \in \N_0
}$. All experiments are carried out with either 
\begin{subequations}
\begin{align} 
 K_\ell &\coloneqq K_0 + 2 \cdot \lfloor \log_2(\ell + 1)/4 \rfloor \quad \text{ for all } \ell \in \N_0 \text{ or} \label{eq:Kvary} \\ 
 K_\ell &\coloneqq K_0 \quad \text{ for all } \ell \in \N_0, \label{eq:Kconst}
\end{align}
\end{subequations}
for different values of $K_0 \in 2\N_0$. It is readily verified that in both these cases, the sequence $(K_\ell)_{\ell \in \N_0}$ satisfies~\eqref{eq:sumdiv} and that the case~\eqref{eq:Kconst} with $K_0 = 0$ ($K_\ell \equiv 0$) corresponds to a single-level SGFEM; see~\cite{Eigel2014a,Bespalov2019a}.
\subsection{Test problem 1 (square domain, finite expansion)} \label{subsec:squareex}
\begin{figure}[!ht]
  \subfloat{
  \resizebox{\textwidth}{!}{
    \includegraphics[scale = 0.02]{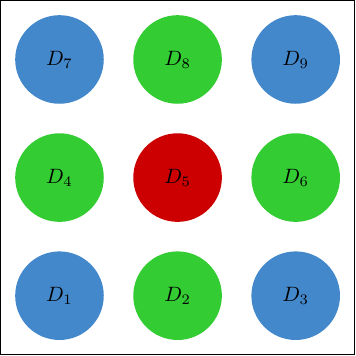}
  }}
   \caption{Computational domain for the test problem in §\ref{subsec:squareex}. Subdomains $D_m$, $m = 1, \ldots , 9$, are depicted in different colors, where each color corresponds to a different value of $k_m$.}
   \label{fig:cookie}
 \end{figure}

We consider the so-called \emph{cookie problem} (see, e.g.,~\cite{Ballani2015,Bespalov2021,Bespalov2021a}), i.e., model problem~\eqref{eq:modelproblem} with computational domain $\Omega = (0,1)^2$, right-hand side $f \equiv 1$, and the diffusion coefficient given by~\eqref{eq:affine} with $a_0 = 1$, $a_m = k_m \chi_{D_m}$ for $m = 1, \ldots , 9$, and $a_m = 0$ for $m \geq 10$, where $D_m$, $m = 1, \ldots , 9$, are pairwise disjoint circles with radius $1/8$ and centers $((2i - 1)/6,(2j - 1)/6)$ for $i,j = 1,2,3$ as depicted in Figure~\ref{fig:cookie}, and
\begin{equation*}
 k_m = \begin{cases}
        0.5 & \text{if } m = 1,3,7,9, \\
        0.7 & \text{if } m = 2,4,6,8, \\
        0.9 & \text{if } m = 5.
       \end{cases}
\end{equation*}

\begin{figure}[!ht]
  \resizebox{\textwidth}{!}{
   \subfloat{
    \stackunder{{\Huge $\nu=\eps_4$}}{\includegraphics{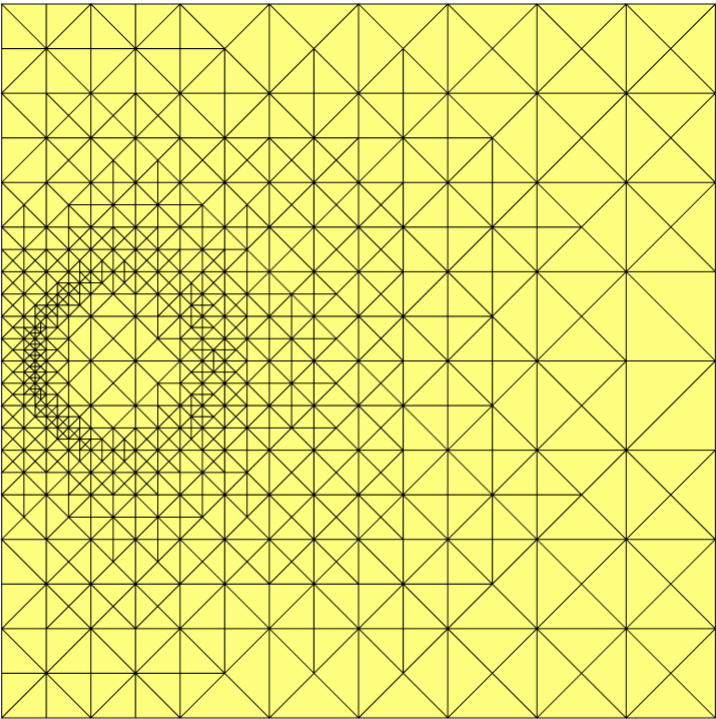}}
   }
   \subfloat{
    \stackunder{{\Huge $\nu=\eps_5$}}{\includegraphics{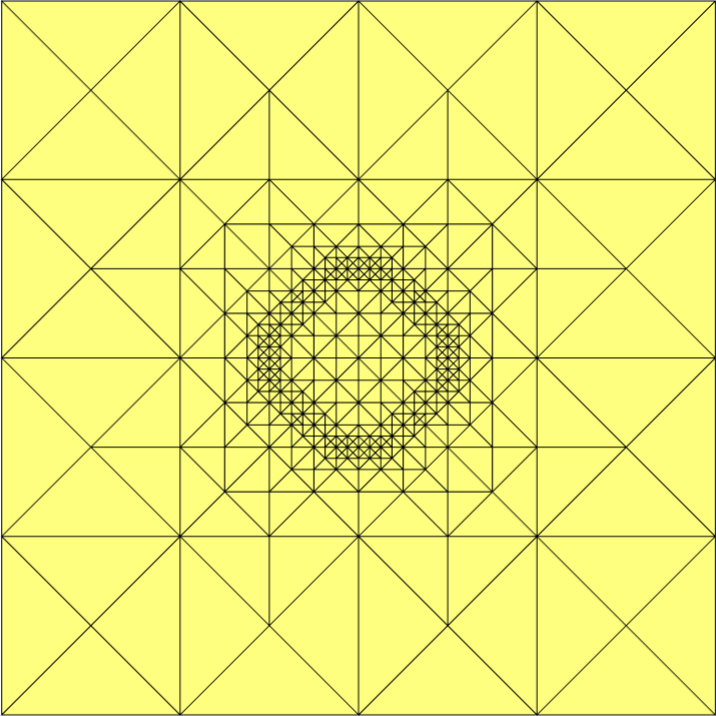}}
   }

   \subfloat{
    \stackunder{{\Huge $\nu=\eps_9$}}{\includegraphics{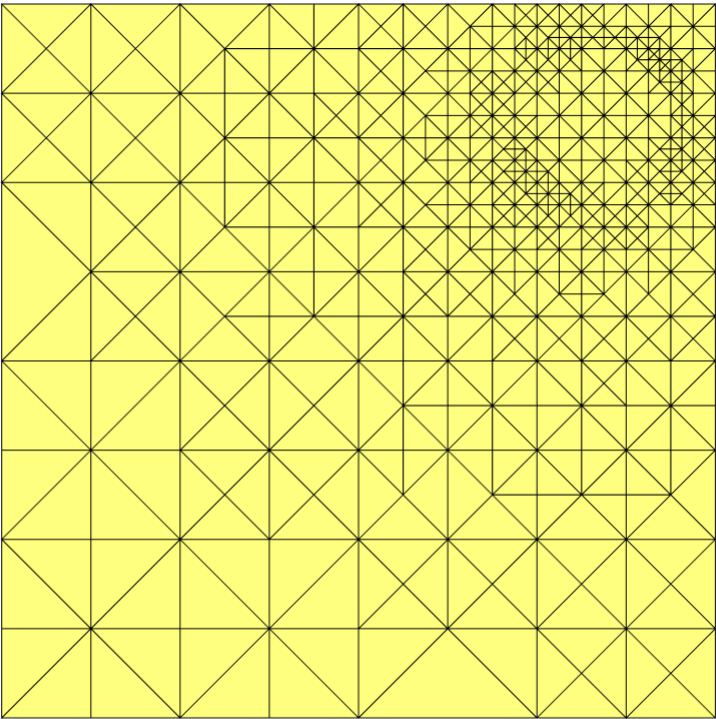}}
   }
  }
   \caption{Adaptively generated meshes $\TT_{\ell,\nu}$ after $\ell = 35$ iterations for $\nu = \eps_4$~(left), $\nu = \eps_5$~(center), and $\nu = \eps_9$~(right) in the test problem in §\ref{subsec:squareex}.}

   \label{fig:cookiemeshes}
 \end{figure}

In our first experiment, we set $(K_\ell)_{\ell \in \N_0}$ as in~\eqref{eq:Kvary} with $K_0 = 8$ and run Algorithm~\ref{alg:adapt} with an initial multilevel structure
\begin{equation*}
 [\PPP_0,(\TT_{0,\nu})_{\nu \in \PPP_0}] = [\set{\vec{0}} \cup \set{\eps_j \given j = 1, \ldots, 9},(\TT_0)_{j = 0}^9]
\end{equation*}
for different marking parameters $\theta \in \set{0.2,0.4,0.6,0.8,1}$; here, the initial mesh $\TT_0$ consists of $32$ triangles. 
 \begin{figure}[!ht]
 \resizebox{\textwidth}{!}{
  \subfloat{
   \includegraphics{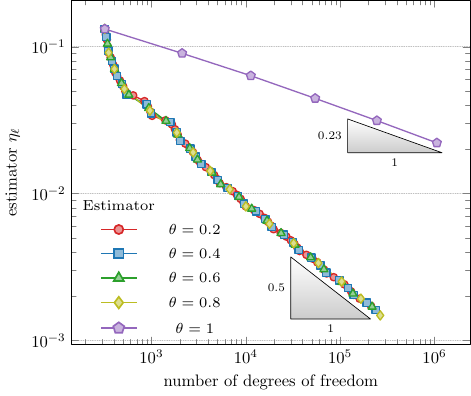}
  }
  \subfloat{
   \includegraphics{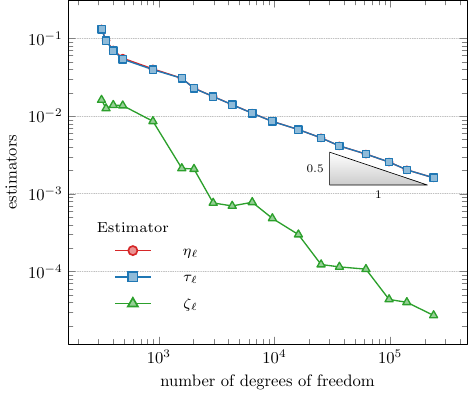}
  }
 }
  \caption{Evolution of the total error estimator $\eta_\ell$ for the SGFEM with $K_\ell$ as in~\eqref{eq:Kvary} with $K_0 = 8$ and different marking parameters $0 < \theta \leq 1$~(left), and the total error estimator $\eta_\ell$ and its contributions $\tau_\ell$ and $\zeta_\ell$~(right) in the test problem in §\ref{subsec:squareex}.}
  \label{fig:cookieerror}
\end{figure}
 The evolution of the total error estimator $\eta_\ell = \eta_\ell(u_\ell^\star)$ (see~\eqref{eq:estnot}) is depicted in Figure~\ref{fig:cookieerror}~(left). We observe that the total error estimators $\eta_\ell$ decrease with optimal rate $\eta_\ell = \OO \Bigl( \bigl( \sum_{\nu \in \PPP_\ell} \#\TT_{\ell,\nu} \bigr)^{-1/2} \Bigr)$ for $\theta \in \set{0.2,0.4,0.6,0.8}$, whereas uniform refinement $(\theta = 1)$ leads to suboptimal convergence. Figure~\ref{fig:cookieerror}~(right) shows the evolution of the total error estimator $\eta_\ell$ alongside its contributions $\tau_\ell = \tau_\ell(u_\ell^\star)$ and $\zeta_\ell = \zeta_\ell(u_\ell^\star)$ from~\eqref{eq:spatialdef}--\eqref{eq:estnot} for the case $\theta = 0.4$. Note that $\norm{a_m}_{W^{1,\infty}(\Omega)} = \norm{a_m}_{L^\infty(\Omega)}$ in the definition of $\zeta_\ell$ (see~\eqref{eq:paramdef}), since $\nabla a_m(x) = 0$ almost everywhere for all $m = 1, \ldots , 9$.
 Figure~\ref{fig:cookiemeshes} shows the adaptively refined meshes $\TT_{\ell,\eps_4}$, $\TT_{\ell,\eps_5}$, and $\TT_{\ell,\eps_9}$ in the case of $\theta = 0.4$ for $\ell = 35$, and we observe that the meshes are able to adequately resolve the individual subdomains $D_m$.

\begin{figure}[!ht]
 \resizebox{\textwidth}{!}{
   \includegraphics[scale=0.04]{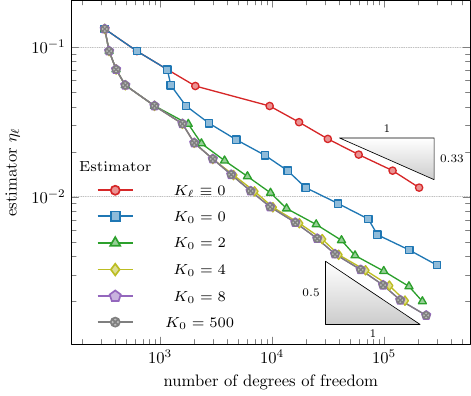}
 }
  \caption{Evolution of the total error estimator $\eta_\ell$ for the SGFEM with $\theta = 0.4$ and different choices of $K_\ell$ as in~\eqref{eq:Kvary} with $K_0 \in \set{0,2,4,8,500}$ and~\eqref{eq:Kconst} with $K_0 = 0$ (single-level SGFEM) in the test problem in §\ref{subsec:squareex}.}
  \label{fig:cookieK}
\end{figure}

In the second experiment, we set $\theta = 0.4$ and run Algorithm~\ref{alg:adapt} with the same initial multilevel structure as in the first experiment for different choices of $(K_\ell)_{\ell \in \N_0}$ as in~\eqref{eq:Kvary} with $K_0 \in \set{0,2,4,8,500}$ and in~\eqref{eq:Kconst} with $K_0 = 0$ (single-level SGFEM).
Figure~\ref{fig:cookieK} shows the evolution of the total error estimator $\eta_\ell$. We observe very similar convergence behavior for $K_0 \in \set{4,8,500}$, whereas choosing $K_0 \in \set{0,2}$ leads to the same reduced convergence rates as observed for the single-level SGFEM in~\cite{Eigel2014a}.
This indicates that the choice of compatibility degrees $K_\ell$ has a significant influence on the performance of the adaptive algorithm and that, in general, larger values of $K_\ell$ are desirable. However, the negligible difference between the convergence behavior for $K_0 = 8$ and $K_0 = 500$ suggests that, in practice, it is not necessary to choose $K_\ell$ very large to achieve optimal convergence rates.
\subsection{Test problem 2 (L-shaped domain, infinite expansion)} \label{subsec:lshapeex}

We consider problem~\eqref{eq:modelproblem} with computational domain $\Omega = (-1,1)^2 \setminus [-1,0]^2$, right-hand side $f \equiv 1$, and the diffusion coefficient given by~\eqref{eq:affine} with $(a_m)_{m \in \N_0}$ being the planar Fourier modes, i.e.,
\begin{equation} \label{eq:fourier}
 a_0 \equiv 1,
 \quad a_m(x_1,x_2) = \alpha m^{-\sigma} \cos( 2\pi \beta_1(m) x_1)\cos(2 \pi \beta_2(m) x_2)
\end{equation}
for $m \in \N$ and $(x_1,x_2) \in \Omega$.
Here, the functions $\beta_1$ and $\beta_2$ are defined as $\beta_1(m) \coloneqq m - k(m)(k(m) + 1)/2$ and $\beta_2(m) \coloneqq k(m) - \beta_1(m)$ with $k(m) \coloneqq \lfloor \sqrt{1/4 + 2m} - 1/2 \rfloor$ for $m \in \N$, $\sigma > 1$, and $\alpha > 0$ chosen such that~\eqref{eq:asumbound} holds with $\gamma = 0.9$. We will consider~\eqref{eq:fourier} with fast ($\sigma = 4$) and slow ($\sigma = 2$) decay of the Fourier coefficients.
We note that~\eqref{eq:ampliesti} is satisfied for all $\sigma > 2$. In order to include the case $\sigma = 2$, we artificially truncate the series expansion~\eqref{eq:affine} by setting $a_m = 0$, $m > m_0$, for some large $m_0 \in \N$ (in this case, $m_0 = 1000$).

In the first experiment for this problem, we set $(K_\ell)_{\ell \in \N_0}$ as in~\eqref{eq:Kvary} with $K_0 = 8$ and run Algorithm~\ref{alg:adapt} with an initial multilevel structure $[\set{\vec{0}},(\TT_0)]$ and with different marking parameters $\theta \in \set{0.2,0.4,0.6,0.8,1}$; here, the initial mesh $\TT_0$ consists of $24$ triangles. 
 \begin{figure}[!ht]
 \resizebox{\textwidth}{!}{
  \subfloat{
   \includegraphics{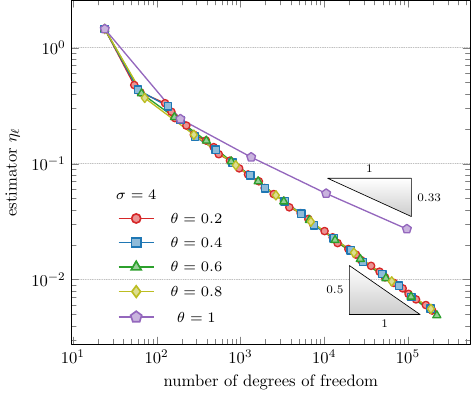}
  }
  \subfloat{
   \includegraphics{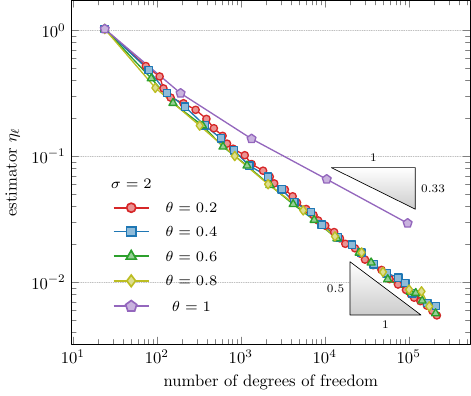}
  }
 }
  \caption{Evolution of the total error estimator $\eta_\ell$ for the SGFEM with $K_\ell$ as in~\eqref{eq:Kvary} with $K_0 = 8$ and different marking parameters $0 < \theta \leq 1$ for fast ($\sigma = 4$, left) and slow ($\sigma = 2$, right) decay of the Fourier coefficients in the test problem in §\ref{subsec:lshapeex}.}
  \label{fig:lrates}
\end{figure}
 \begin{figure}[!ht]
 \resizebox{\textwidth}{!}{
  \subfloat{
   \includegraphics{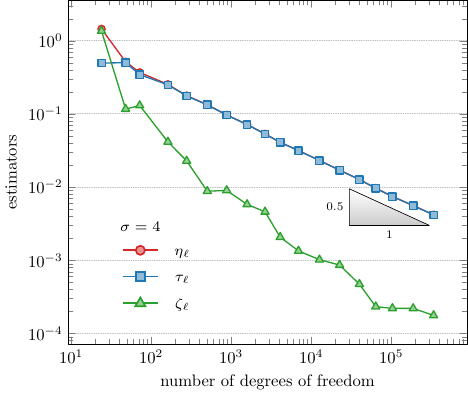}
  }
  \subfloat{
   \includegraphics{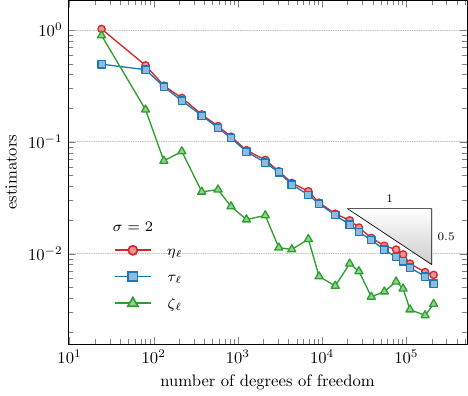}
  }
 }
  \caption{Evolution of the total error estimator $\eta_\ell$ and its contributions $\tau_\ell$ and $\zeta_\ell$ for the SGFEM with $K_\ell$ as in~\eqref{eq:Kvary} with $K_0 = 8$ and $\theta = 0.4$ for fast ($\sigma = 4$, left) and slow ($\sigma = 2$, right) decay of the Fourier coefficients in the test problem in §\ref{subsec:lshapeex}.}
  \label{fig:lcont}
\end{figure}
The evolution of the total error estimator $\eta_\ell$ (see~\eqref{eq:estnot}) for fast ($\sigma = 4$, left) and slow ($\sigma = 2$, right) decay of the Fourier coefficients is depicted in Figure~\ref{fig:lrates}. 
As in §\ref{subsec:squareex}, we observe that the total error estimators $\eta_\ell$ decrease with the optimal rate for $\theta \in \set{0.2,0.4,0.6,0.8}$, whereas uniform refinement $(\theta = 1)$ leads to a suboptimal convergence rate.
Figure~\ref{fig:lcont} shows the evolution of the total error estimator $\eta_\ell$ alongside its contributions $\tau_\ell = \tau_\ell(u_\ell^\star)$ and $\zeta_\ell = \zeta_\ell(u_\ell^\star)$ from~\eqref{eq:spatialdef}--\eqref{eq:estnot} for the case $\theta = 0.4$ as well as for fast ($\sigma = 4$, left) and slow ($\sigma = 2$, right) decay of the Fourier coefficients.
We note that in earlier works (see, e.g.,~\cite{Eigel2014,Eigel2014a}), the parametric error indicator $\zeta_\ell$ is computed using the $L^{\infty}(\Omega)$-norm of the coefficients $a_m$, whereas here it is evaluated using their $W^{1,\infty}(\Omega)$-norm (see~\eqref{eq:paramdef}); this leads to a larger contribution of $\zeta_\ell$ to the total error estimator $\eta_\ell$. Nevertheless, as in the existing literature, we observe that the spatial error indicators $\tau_\ell$ still dominate the parametric error indicators $\zeta_\ell$. 
\begin{figure}[!ht]
  \resizebox{\textwidth}{!}{
   \subfloat{
    \stackunder{{\Huge $\nu=(1,0,0,\ldots)$}}{\includegraphics{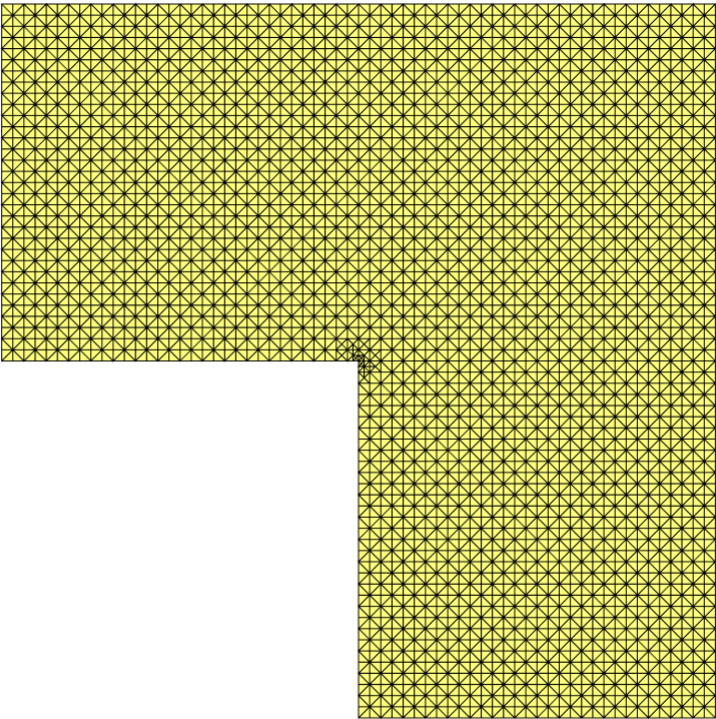}}
   }
   \subfloat{
    \stackunder{{\Huge $\nu=(2,0,0,\ldots)$}}{\includegraphics{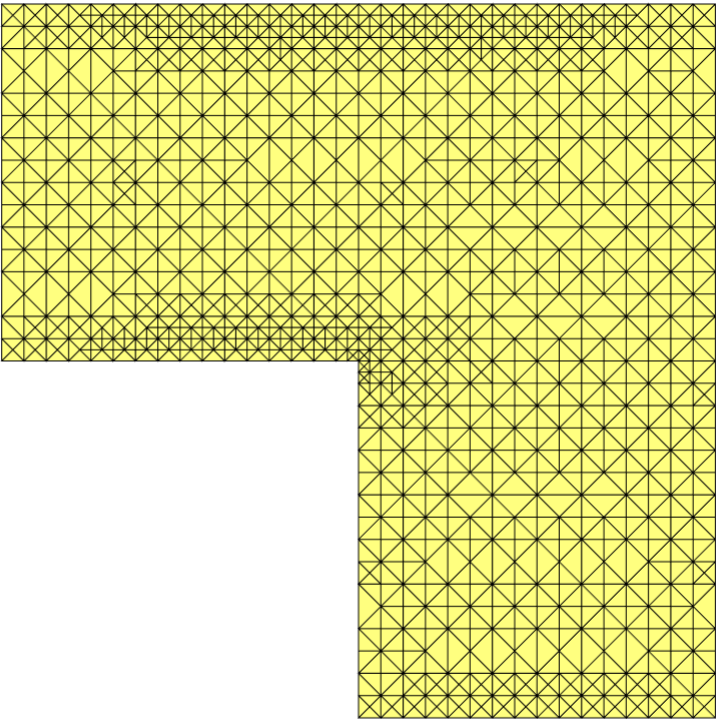}}
   }

   \subfloat{
    \stackunder{{\Huge $\nu=(3,0,0,\ldots)$}}{\includegraphics{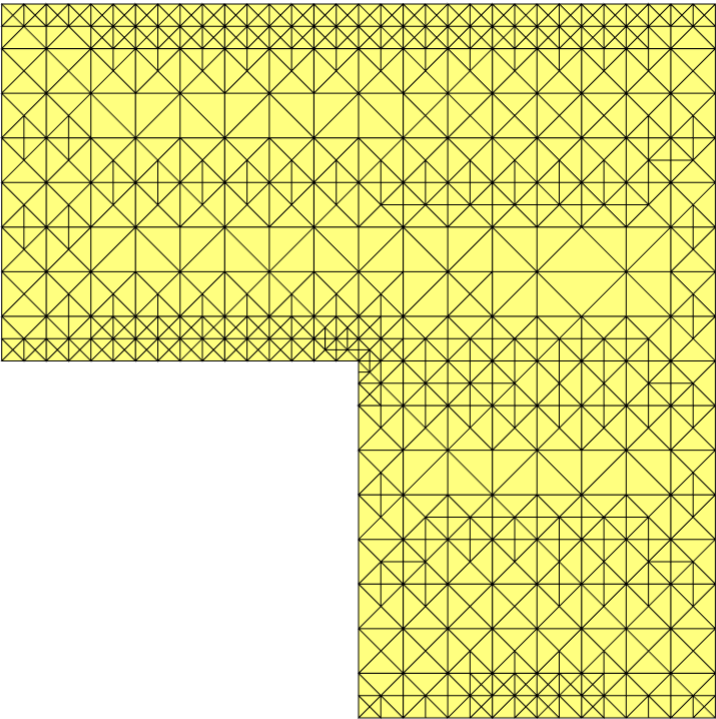}}
   }
  }
   \caption{Adaptively generated meshes $\TT_{\ell,\nu}$ after $\ell = 35$ iterations for $\nu = (1,0,0,\ldots)$~(left), $\nu = (2,0,0,\ldots)$~(center), and $\nu = (3,0,0,\ldots)$~(right) in the test problem in §\ref{subsec:lshapeex}.}

   \label{fig:lmeshes}
\end{figure}
\begin{table}[!ht]
  \resizebox{\textwidth}{!}{
   \includegraphics[scale = 0.1]{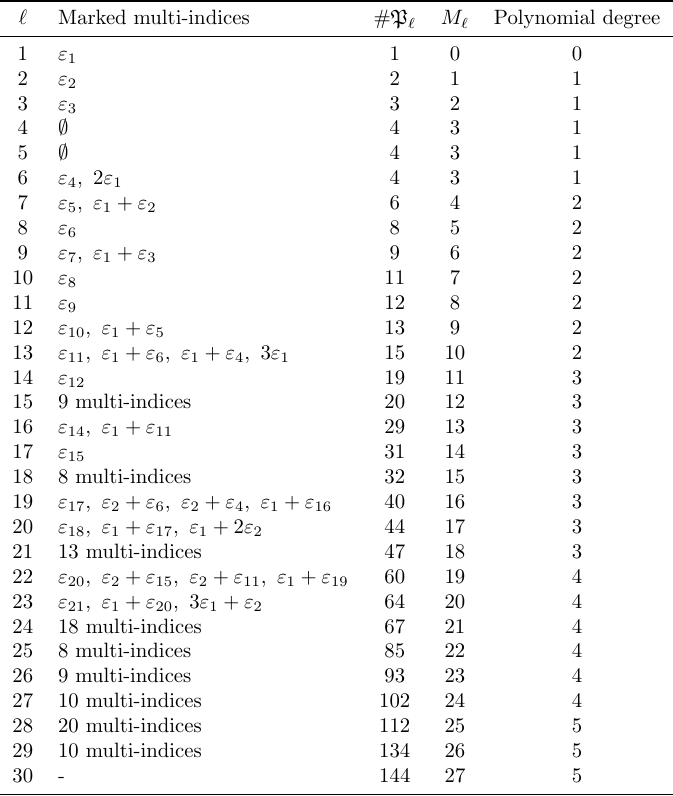}}
   \caption{Evolution of the index set $\PPP_\ell$, the number $M_\ell$ of active components, and the highest total polynomial degree for  the SGFEM with $(K_\ell)_{\ell \in \N_0}$ as in~\eqref{eq:Kvary} with $K_0 = 8$ and $\theta = 0.4$ for slow ($\sigma = 2$) decay of the Fourier coefficients in the test problem in §\ref{subsec:lshapeex}.}
   \label{fig:lindex}
\end{table}

Figure~\ref{fig:lmeshes} shows the adaptively refined meshes $\TT_{\ell,(1,0,\ldots)}$, $\TT_{\ell,(2,0,\ldots)}$, and $\TT_{\ell,(3,0,\ldots)}$ in the case of $\theta = 0.
4$ and $(K_\ell)_{\ell \in \N_0}$ as in~\eqref{eq:Kvary} with $K_0 = 2$ for $\ell = 35$, and we observe that the meshes are able to resolve the singularity at the reentrant corner while verifying $K_\ell$-compatibility with $K_{35} = 2 + 2 \cdot \lfloor \log_2(36)/4 \rfloor = 4$. 
Table~\ref{fig:lindex} shows the evolution of the index set $\PPP_\ell$ along with the number $M_\ell$ of active components and the highest total polynomial degree $\max_{\nu \in \PPP_\ell} \abs{\nu}$ in the case of $\sigma = 2$ in~\eqref{eq:fourier} and for the D\"orfler marking parameter $\theta = 0.4$. 
We observe that, due to the larger contribution of the parametric error indicator together with combined D\"orfler marking~\eqref{eq:doerfler}, the adaptive algorithm adds more multi-indices to $\PPP_\ell$ than in earlier works (see, e.g.,~\cite{Bespalov2021,Bespalov2021a}). 
Table~\ref{fig:lindex} also indicates that for all but one iteration, the index set $\PPP_{\ell + 1}$ ($\ell \geq 1$) generated by adding only the marked multi-indices from $\QQQ_\ell$ is monotone; it is only iteration $29$ where monotonicity of $\PPP_{29}$ has to be enforced by adding two extra multi-indices to the corresponding marking set of $20$ multi-indices ($\#\PPP_{29} = \#\PPP_{28} + 22$).

\begin{figure}[!ht]
 \resizebox{\textwidth}{!}{
  \subfloat{
   \includegraphics{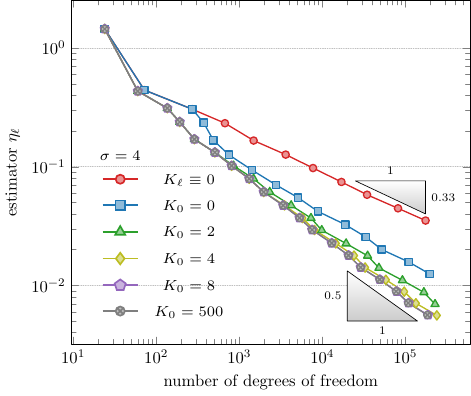}
  }
  \subfloat{
   \includegraphics{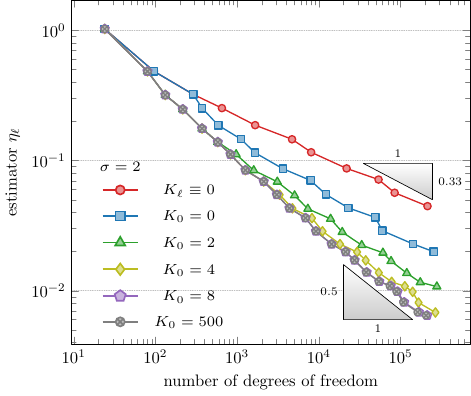}
  }
 }
  \caption{Evolution of the total error estimator $\eta_\ell$ for the SGFEM with $\theta = 0.4$ and different choices of $K_\ell$ as in~\eqref{eq:Kvary} with $K_0 \in \set{0,2,4,8,500}$ and~\eqref{eq:Kconst} with $K_0 = 0$ (single-level SGFEM) for fast ($\sigma = 4$, left) and slow ($\sigma = 2$, right) decay of the Fourier coefficients in the test problem in §\ref{subsec:lshapeex}. }
  \label{fig:lK}
\end{figure}

In our last experiment, we set $\theta = 0.4$ and run Algorithm~\ref{alg:adapt} with the same initial multilevel structure as in the first experiment for the test problem in §\ref{subsec:lshapeex} for different choices of $(K_\ell)_{\ell \in \N_0}$ as in~\eqref{eq:Kvary} with $K_0 \in \set{0,2,4,8,500}$ as well as in~\eqref{eq:Kconst} with $K_0 = 0$. Figure~\ref{fig:lK} shows the evolution of the total error estimator $\eta_\ell$ for fast ($\sigma = 4$, left) and slow ($\sigma = 2$, right) decay of the Fourier coefficients. As in §\ref{subsec:squareex}, we observe that there is little to no difference in the convergence behavior in the case of $(K_\ell)_{\ell \in \N_0}$ given by~\eqref{eq:Kvary} with $K_0 \in \set{4,8,500}$. 
However, in the cases of $(K_\ell)_{\ell \in \N_0}$ given by~\eqref{eq:Kvary} with $K_0 = 2$ and $(K_\ell)_{\ell \in \N_0}$ given by~\eqref{eq:Kvary} and~\eqref{eq:Kconst} with $K_0 = 0$ we observe reduced convergence rates.

\printbibliography

\makeatletter
\@setaddresses
\makeatother
 \appendix 

\section{Proofs of Proposition~\ref{prop:refine} and Proposition~\ref{prop:compilize}} \label{sec:append}

\subsection{Preliminary results}

In the following, we assume that the initial triangulation $\TT_0$ is admissible, conforming and $\kappa$-shape regular. The following lemmata are key to the proofs of the remaining results. The first lemma is an immediate consequence of~\cite[Corollary~4.6]{Stevenson2008}.
\begin{lemma} \label{lem:nvb}
Let $\TT$ be a triangulation and $T,T' \in \TT$ with $E \coloneqq T \cap T' \in \EE$.
Then, exactly one of the following cases holds for the reference edges $E_T$ and $E_{T'}$, respectively$\colon$
\begin{enumerate}[label = \rm(\roman*)]
\item $E_T = E = E_{T'} \quad$ and $\quad \gen{T} = \gen{T'}$.
\item $E_T \neq E \neq E_{T'} \quad$ and $\quad \gen{T} = \gen{T'}$.
\item $E_T = E \neq E_{T'} \quad$ and $\quad \gen{T} = \gen{T'} + 1$.
\item $E_T \neq E = E_{T'} \quad$ and $\quad \gen{T} + 1 = \gen{T'}$.
\end{enumerate}
\end{lemma} 
The proof of Proposition~\ref{prop:refine} requires the following auxiliary result.
\begin{lemma} \label{lem:nvb2}
 Let $\TT_1,\TT_2$ be triangulations and $T_1 \in \TT_1$, $T_2 \in \TT_2$ with $T_1 \subseteq T_2$ and $\gen{T_1} - \gen{T_2} \in 2\N - 1$. Then, there holds $\abs{E_{T_1} \cap E_{T_2}} = 0$, where $\abs{\, \cdot \,}$ denotes the $1$-dimensional Hausdorff measure.
\end{lemma}
\begin{proof}
Let $n \in \N$ such that $\gen{T_1} - \gen{T_2} = 2n - 1$. We proceed to prove the statement by induction over $n$.

\textbf{Base case (\boldsymbol{$n = 1$}).} Since $\gen{T_1} = \gen{T_2} + 1$, $T_1$ is generated by bisecting $T_2$ once. Hence, $E_{T_2}$ is bisected and $E_{T_1} \in \EE_{T_2} \setminus \set{E_{T_2}}$. This proves $\abs{E_{T_1} \cap E_{T_2}} = 0$.

\textbf{Induction hypothesis.} Let $n \in \N$ and suppose the claim holds for all $1 \leq j \leq n$.

\textbf{Induction step.} Let $T_1 \in \TT_1$, $T_2 \in \TT_2$ with $T_1 \subseteq T_2$ and $\gen{T_1} - \gen{T_2} = 2(n + 1) - 1 = 2n + 1$. Since $\gen{T_1} - \gen{T_2} \geq 3$, there exist triangles $T,T'$ with $T_1 \subsetneqq T \subsetneqq T' \subsetneqq T_2$ and
\begin{equation*}
 \gen{T_1}
 = \gen{T} + 1
 = \gen{T'} + 2
\end{equation*}
and there exist edges $E,E'$  with  $E_{T_1} \subseteq E \subseteq E'$.
Due to the induction hypothesis, there holds $\abs{E_{T_1} \cap E_T} = 0$, which yields $E \neq E_T$. Since $T$ is generated by bisecting $T'$ once and $E' \supseteq E \neq E_T$, we obtain $E_{T'} = E'$. Finally, due to the induction hypothesis, that is applicable because $\gen{T'} - \gen{T_2} = 2n - 1$, we have $0 = \abs{E_{T'} \cap E_{T_2}} = \abs{E' \cap E_{T_2}} \geq \abs{E_{T_1} \cap E_{T_2}}$. This concludes the proof.
\end{proof}

\begin{proof}[\textbf{Proof of Proposition~\ref{prop:refine}}]
Since $\MM_{\nu}^k \subseteq \EE_{\ell,\nu}$ for all $\nu \in \PPP_\ell$ and $k \in \N_0$ (here, $\EE_{\ell,\nu}$ denotes the set of all edges in $\TT_{\ell,\nu}$), the algorithm terminates after at most $\displaystyle{\sum_{\nu \in \PPP_\ell} \abs{\EE_{\ell,\nu}}} < +\infty$ iterations. 

Since $\MM_\nu^0$ contains all reference edges of the marked elements $\MM_{\ell,\nu}$ for all $\nu \in \PPP_\ell$, all elements in $\MM_{\ell,\nu}$ are bisected in $\TT_{\ell + 1,\nu}$.

To show that the generated family $(\TT_{\ell + 1,\nu})_{\nu \in \PPP_\ell}$ is $K_{\ell + 1}$-compatible, we note that, due to the assumption on $(K_\ell)_{\ell \in \N_0}$, there holds $K_{\ell + 1} - K_\ell \in 2\N_0$. 
Since one iteration of refinement increases the level of a triangle by at most $2$, usual NVB-refinement of every triangulation in the $K_\ell$-compatible family $(\TT_{\ell,\nu})_{\nu \in \PPP_\ell}$ yields a $(K_{\ell} + 2)$-compatible family $(\TT_{\ell + 1,\nu})_{\nu \in \PPP_\ell}$. 
Therefore, if $K_{\ell + 1} - K_\ell \in 2\N$, $(\TT_{\ell + 1,\nu})_{\nu \in \PPP_\ell}$ is $K_{\ell + 1} ( \geq K_\ell + 2)$-compatible. 
Thus, we only need to consider the case $K_{\ell + 1} = K_\ell$.

To this end, let $\nu \in \PPP_\ell$, $\mu \in \PPP_\ell[\nu]$, $T_{\mu} \in \TT_{\ell + 1,\mu}$, and $T_{\nu} \in \TT_{\ell + 1,\nu}$ with $|T_{\mu} \cap T_{\nu}| > 0$. 
Without loss of generality, we can assume that $T_{\nu} \subseteq T_{\mu}$ and we need to show that $\gen{T_{\nu}} - \gen{T_{\mu}} \leq K_{\ell + 1}$. 
The remainder of the proof is split into four cases.

\textbf{Case~1 (\boldsymbol{$T_\nu \in \TT_{\ell,\nu}$} and \boldsymbol{$T_\mu \in \TT_{\ell,\mu}$}).}
Since both triangles are in the meshes which are $K_\ell$-compatible by assumption, there holds $\gen{T_{\nu}} - \gen{T_{\mu}} \leq K_\ell = K_{\ell + 1}$.

\textbf{Case~2 (\boldsymbol{$T_\nu \in \TT_{\ell,\nu}$} and \boldsymbol{$T_\mu \notin \TT_{\ell,\mu}$}).}
Since $T_{\mu} \in \TT_{\ell + 1,\mu} \setminus \TT_{\ell,\mu}$, there exists $T_{\mu}' \in \TT_{\ell,\mu}$ with $T_{\mu} \subsetneqq T_{\mu}'$ and 
\begin{equation*}
 \gen{T_{\mu}} - 2
 \leq \gen{T_{\mu}'} 
 \leq \gen{T_{\mu}} - 1.
\end{equation*}
This, together with $T_{\nu} \subseteq T_{\mu} \subseteq T_{\mu}'$ and $K_\ell$-compatibility of $(\TT_{\ell,\nu})_{\nu \in \PPP_\ell}$ yields
\begin{equation*}
 \gen{T_{\nu}} 
 \leq \gen{T_{\mu}'} + K_\ell 
 \leq \gen{T_{\mu}} + K_\ell - 1
 < \gen{T_{\mu}} + K_{\ell + 1}.
\end{equation*}

\textbf{Case~3 (\boldsymbol{$T_\nu \notin \TT_{\ell,\nu}$} and \boldsymbol{$T_\mu \in \TT_{\ell,\mu}$}).}
We will argue by contradiction and assume that $\gen{T_{\nu}} > \gen{T_{\mu}} + K_{\ell + 1} = \gen{T_\mu} + K_\ell$. 
Since $T_{\nu} \in \TT_{\ell + 1,\nu} \setminus \TT_{\ell,\nu}$, there exists $T_{\nu}' \in \TT_{\ell,\nu}$ with $T_{\nu} \subsetneqq T_{\nu}'$ and NVB ensures that
\begin{equation} \label{eq:childest}
 \gen{T_{\nu}} - 2
 \leq \gen{T_{\nu}'} 
 \leq \gen{T_{\nu}} -1.
\end{equation}
Since $|T_{\nu}' \cap T_{\mu}| \geq |T_{\nu} \cap T_{\mu}| > 0$ and $K_\ell \geq 2$, there holds
\begin{equation} \label{eq:Lchildest}
 \gen{T_{\nu}'} 
 \geq \gen{T_{\nu}} - 2 
 > \gen{T_{\mu}} + K_\ell - 2 
 \geq \gen{T_{\mu}}.
\end{equation}
Hence, $T_{\nu}' \subseteq T_{\mu}$, and, due to $K_\ell$-compatibility of $(\TT_{\ell,\nu})_{\nu \in \PPP_\ell}$, we obtain $\gen{T_{\nu}'} \leq \gen{T_{\mu}} + K_\ell$, which, together with \eqref{eq:Lchildest}, implies
\begin{equation*}
 \gen{T_{\mu}} + K_\ell - 1 
 \eqreff{eq:Lchildest}{\leq} \gen{T_{\nu}'} 
 \leq \gen{T_{\mu}} + K_\ell.
\end{equation*}
Thus, $\gen{T_\nu'} - \gen{T_\mu} \in \set{K_\ell - 1,K_\ell}$. If $\gen{T_{\nu}'} - \gen{T_{\mu}} = K_\ell = K_{\ell + 1}$, then, since $T_{\nu}'$ is bisected, $T_{\mu}$ must also be bisected in Algorithm~\ref{alg:refine}(a), which contradicts the fact that $T_{\mu}$ is contained in $\TT_{\ell + 1,\mu} \cap \TT_{\ell,\mu}$. 

Now we consider the case
\begin{equation} \label{eq:leveleq}
 \gen{T_{\nu}'} - \gen{T_{\mu}} 
 = K_\ell - 1.
\end{equation}	 
Then, the assumption $\gen{T_{\nu}} > \gen{T_{\mu}} + K_\ell$ implies
\begin{equation*}
 \gen{T_{\nu}'}
 \eqreff{eq:leveleq}{=} \gen{T_{\mu}} + K_\ell - 1
 \leq \gen{T_{\nu}} - 2.
\end{equation*}
Together with \eqref{eq:childest}, this yields $\gen{T_{\nu}'} = \gen{T_{\nu}} - 2$. 
Hence, $T_\nu$ is generated by bisecting $T_{\nu}'$ twice. Since two refinements take place, at least two edges must be marked, i.e., there exists an edge $E_{\nu} \in \EE_{T_{\nu}'} \setminus \set{E_{T_{\nu}'}}$ such that $E_{T_{\nu}'},E_{\nu} \in \MM_{\nu}^k$ for some $k \in \N_0$. 
Since $E_{\nu} \in \MM_\nu^k$ and $\MM_\nu^k$ contains only reference edges of triangles in $\TT_{\ell,\nu}$, there exists $N_{\nu}' \in \TT_{\ell,\nu}$ with $E_{T_\nu'} \neq E_{N_\nu'} = E_{\nu} = T_{\nu}' \cap N_{\nu}'$. 
Hence, Lemma \ref{lem:nvb} and~\eqref{eq:leveleq} establish
\begin{equation} \label{eq:levelNnu}
 \gen{N_{\nu}'} 
 = \gen{T_{\nu}'} + 1 
 \eqreff{eq:leveleq}{=} \gen{T_{\mu}} + K_\ell
 = \gen{T_{\mu}} + K_{\ell + 1}.
\end{equation}
We first consider the case $N_{\nu}' \subseteq T_{\mu}$. 
Since $N_{\nu}'$ is bisected, also $T_\mu$ is bisected by Algorithm~\ref{alg:refine}(a), contradicting the fact that $T_\mu$ is contained in $\TT_{\ell + 1,\mu} \cap \TT_{\ell,\mu}$.

It remains to consider the situation when $N_{\nu}' \nsubseteq T_{\mu}$. 
Since $T_{\nu}' \subseteq T_{\mu}$ and $E_{\nu} = T_{\nu}' \cap N_{\nu}'$, there exists $E_{\mu} \in \EE_{T_{\mu}}$ and $N_{\mu} \in \TT_{\ell,\mu}$ with $E_{\mu} \supseteq E_{\nu}$, $E_{\mu} = T_{\mu} \cap N_{\mu}$ and $N_{\nu}' \subseteq N_{\mu}$. 
Since no initial statement about the reference edge can be made, we consider separately the following three cases (according to Lemma~\ref{lem:nvb})$\colon$
\begin{equation*}
 \gen{T_{\mu}}  - \gen{N_{\mu}}
 \in \set{-1,0,1}.
\end{equation*}

\textbf{Case 3.1 (\boldsymbol{$\gen{T_{\mu}} - \gen{N_{\mu}} = 1$}).}
By~\eqref{eq:leveleq} and~\eqref{eq:levelNnu} we obtain
\begin{equation*}
 \gen{N_{\mu}} 
 = \gen{T_{\mu}} - 1 
 \eqreff{eq:leveleq}{=} \gen{T_{\nu}'} - K_\ell 
 \eqreff{eq:levelNnu}{=} \gen{N_{\nu}'} - K_\ell - 1,
\end{equation*}
which, together with $N_{\nu}' \subseteq N_{\mu}$, contradicts $K_\ell$-compatibility of the family $(\TT_{\ell,\nu})_{\nu \in \PPP_\ell}$.

\textbf{Case 3.2 (\boldsymbol{$\gen{T_{\mu}} - \gen{N_{\mu}} = -1$}).}
By~\eqref{eq:leveleq} and~\eqref{eq:levelNnu} we obtain
\begin{equation} \label{eq:case32}
 \gen{N_{\mu}} 
 = \gen{T_{\mu}} + 1
 \eqreff{eq:leveleq}{=} \gen{T_{\nu}'} - K_\ell + 2 
 \eqreff{eq:levelNnu}{=} \gen{N_{\nu}'} - K_\ell +1.
\end{equation}
Since $K_\ell - 1 \in 2\N - 1$, we can apply Lemma~\ref{lem:nvb2} to $N_\nu'$ and $N_\mu$ and obtain that $\abs{E_{N_\nu'} \cap E_{N_\mu}} = 0$. 
Since $E_{N_\nu'} = E_\nu \subseteq E_\mu$, we obtain $E_\mu \neq E_{N_\mu}$. 
However, since $\gen{N_\mu} = \gen{T_\mu} + 1$, Lemma~\ref{lem:nvb} yields that $E_{N_\mu} = N_\mu \cap T_\mu = E_\mu$, which contradicts $E_\mu \neq E_{N_\mu}$. 

\textbf{Case 3.3 (\boldsymbol{$\gen{T_{\mu}} - \gen{N_{\mu}} = 0$}).}
By~\eqref{eq:leveleq} and~\eqref{eq:levelNnu} we obtain
\begin{equation} \label{eq:case33}
 \gen{N_{\mu}} 
 = \gen{T_{\mu}} 
 \eqreff{eq:leveleq}{=} \gen{T_{\nu}'} - K_\ell + 1 
 \eqreff{eq:levelNnu}{=} \gen{N_{\nu}'} - K_\ell,
\end{equation}
which, by Algorithm~\ref{alg:refine}(a), implies that $N_{\mu}$ is bisected. 
If $E_{\mu} = E_{N_{\mu}}$, since $E_\mu = T_\mu \cap N_\mu$, $T_{\mu}$ is also bisected, which contradicts $T_{\mu} \in \TT_{\ell + 1,\mu} \cap \TT_{\ell,\mu}$. 

Therefore, $E_{\mu} \neq E_{N_{\mu}}$. 
Since $N_\mu$ is bisected, there exists $N_{\mu}' \in \TT_{\ell + 1,\mu}$ with $N_{\mu}' \subsetneqq N_{\mu}$ and $E_{\mu} = E_{N_{\mu}'}$. 
Since $N_\mu'$ is generated by bisecting $N_\mu$ exactly once, this yields 
\begin{equation*}
 \gen{N_{\mu}'} = \gen{N_{\mu}} + 1 \eqreff{eq:case33}{=} \gen{N_\nu'} - K_\ell + 1. 
\end{equation*}
Hence, in view of~\eqref{eq:case32}, we are able to argue analogously to Case~3.2 for $N_\mu$ to obtain a contradiction.

\textbf{Case~4 (\boldsymbol{$T_\nu \notin \TT_{\ell,\nu} \text{ and } T_\mu \notin \TT_{\ell,\mu}$}).}
Let $T_{\nu}' \in \TT_{\ell,\nu}$ and $T_{\mu}' \in \TT_{\ell,\mu}$ with $T_{\nu} \subsetneqq T_{\nu}'$ and $T_{\mu} \subsetneqq T_{\mu}'$. 
There holds 
\begin{equation*}
 \gen{T_{\nu}} - \gen{T_{\nu}'} \in \set{1,2} \quad \text{and} \quad \gen{T_{\mu}} - \gen{T_{\mu}'} 
 \in \set{1,2}
\end{equation*}
We consider two subcases.

\textbf{Case~4.1 (\boldsymbol{$\gen{T_{\nu}} - \gen{T_{\nu}'} \leq \gen{T_{\mu}} - \gen{T_{\mu}'}$}).}
There holds $|T_{\nu}' \cap T_{\mu}'| \geq |T_{\nu} \cap T_{\mu}| > 0$ and, since $T_{\nu} \subseteq T_{\mu}$, we obtain $\gen{T_\nu} \geq \gen{T_\mu}$ and
\begin{equation*}
 \gen{T_{\mu}'} 
 \leq \gen{T_{\mu}} - \gen{T_{\nu}} + \gen{T_{\nu}'} 
 \leq \gen{T_{\nu}'},  
\end{equation*}
which implies $T_{\nu}' \subseteq T_{\mu}'$. 
Thus, $K_\ell$-compatibility of the family $(\TT_{\ell,\nu})_{\nu \in \PPP_\ell}$ yields $\gen{T_{\nu}'} \leq \gen{T_{\mu}'} + K_\ell$.
We obtain
\begin{equation*} 
 \gen{T_{\nu}} 
 \stackrel{\text{Case 4.1}}{\leq}  \gen{T_{\mu}} + \gen{T_{\nu}'} - \gen{T_{\mu}'} 
 \leq \gen{T_{\mu}} + K_\ell
 = \gen{T_\mu} + K_{\ell + 1}.
\end{equation*}

\textbf{Case~4.2 (\boldsymbol{$\gen{T_{\nu}} = \gen{T_{\nu}'} + 2 \text{ and } \gen{T_{\mu}} = \gen{T_{\mu}'} + 1$}).}
To show $\gen{T_\nu} \leq \gen{T_\mu} + K_\ell$, we argue by contradiction and assume that $\gen{T_{\nu}} > \gen{T_{\mu}} + K_\ell$. 
Since
\begin{equation*}
 \gen{T_{\nu}'} 
 \stackrel{\text{Case 4.2}}{=} \gen{T_{\nu}} - 2 
 \stackrel{\text{Assump.}}{\geq} \gen{T_{\mu}} + K_\ell - 1 
 \stackrel{\text{Case 4.2}}{=} \gen{T_{\mu}'} + K_\ell,
\end{equation*}
$K_\ell$-compatibility of the family $(\TT_{\ell,\nu})_{\nu \in \PPP_\ell}$ yields 
\begin{equation} \label{eq:case42}
 \gen{T_{\nu}'} = \gen{T_{\mu}'} + K_\ell.
\end{equation}
Since $K_\ell \geq 2$ and $|T_{\nu}' \cap T_{\mu}'| \geq |T_{\nu} \cap T_{\mu}| > 0$, we obtain $T_{\nu}' \subsetneqq T_{\mu}'$. 
Since $T_{\nu}'$ is bisected at least twice, there exists $E_\nu \in \EE_{T_\nu'} \setminus \set{E_{T_\nu'}}$ with $E_\nu \in \MM_{\nu}^k$ for some $k \in \N_0$. 
Since $\MM_\nu^k$ contains only reference edges,  there exists $N_{\nu}' \in \TT_{\ell,\nu}$ with $E_{T_\nu'} \neq E_{N_{\nu}'} = E_{\nu} = T_{\nu}' \cap N_{\nu}'$. 
Lemma~\ref{lem:nvb} yields
\begin{equation} \label{eq:case42a}
\gen{N_{\nu}'} 
= \gen{T_{\nu}'} + 1 
\eqreff{eq:case42}{=} \gen{T_{\mu}'} + K_\ell + 1.
\end{equation}
This implies $N_\nu' \nsubseteq T_\mu'$, since $N_\nu' \subseteq T_\mu'$ would contradict $K_\ell$-compatibility of the family $(\TT_{\ell,\nu})_{\nu \in \PPP_\ell}$.
Since $T_{\nu}' \subsetneqq T_{\mu}'$ and $E_{\nu} = T_{\nu}' \cap N_{\nu}'$, there exists $E_{\mu} \in \EE_{T_{\mu}'}$ and $N_{\mu}' \in \TT_{\ell,\mu}$ with $E_{\mu} \supseteq E_{\nu}$, $E_{\mu} = T_{\mu}' \cap N_{\mu}'$ and $N_{\nu}' \subseteq N_{\mu}'$. 
Hence, Lemma~\ref{lem:nvb} yields
\begin{equation*}
 \gen{N_{\mu}'} - \gen{T_{\mu}'} 
 \in \set{-1,0,1}.
\end{equation*} 
If $\gen{N_{\mu}'} - \gen{T_{\mu}'} \leq 0$, we obtain
\begin{equation*}
 \gen{N_{\mu}'} 
 \leq \gen{T_\mu'} 
 \eqreff{eq:case42}{=} \gen{T_{\nu}'} - K_\ell 
 \eqreff{eq:case42a}{=} \gen{N_{\nu}'} - K_\ell - 1,
\end{equation*} 
which contradicts $K_\ell$-compatibility of the family $(\TT_{\ell,\nu})_{\nu \in \PPP_\ell}$. 

The remaining case is $\gen{N_{\mu}'} - \gen{T_{\mu}'} = 1$. 
In this case, we have
\begin{equation*}
 \gen{N_{\mu}'} 
 = \gen{T_\mu'} + 1 
 \eqreff{eq:case42a}{=} \gen{N_{\nu}'} - K_\ell.
\end{equation*}
Since $N_\nu' \subseteq N_\mu'$ is bisected and $\gen{N_{\nu}'} = \gen{N_{\mu}'} + K_\ell$, Algorithm~\ref{alg:refine}(a) yields that also $N_{\mu}'$ is bisected. 
However, the fact that $\gen{T_\mu'} + 1 = \gen{N_\mu'}$ (see~\eqref{eq:case42a}) together with Lemma~\ref{lem:nvb} yields that $E_{N_\mu'} = E_\mu \neq E_{T_\mu'}$. 
This implies that $T_\mu'$ is bisected at least twice. 
Since $T_\mu \subsetneqq T_\mu'$, $T_\nu \subsetneqq T_\mu$, and $N_\nu' \nsubseteq T_\mu'$, we obtain that $T_\mu$ is generated by bisecting $T_\mu'$ twice, i.e., $\gen{T_\mu} = \gen{T_{\mu'}} + 2$, contradicting $\gen{T_\mu} = \gen{T_\mu'} + 1$. 

Thus, all cases lead to a contradiction, which implies that the initial assumption $\gen{T_{\nu}} > \gen{T_{\mu}} + K_\ell$ is false. 
This concludes the proof.
\end{proof}

\subsection{Proof of Proposition~\ref{prop:compilize}}

\begin{proof}[\textbf{Proof of Proposition~\ref{prop:compilize}}]
Since marked elements are bisected at least once in each iteration, the algorithm terminates after at most $\max(\displaystyle{\max_{T \in \TT_\infty}} \gen{T} - K_\ell,0) < +\infty$ iterations of the $k$-loop.
The remainder of the proof is split into four steps.

\textbf{Step~1 (Proof of \boldsymbol{$\TT\mkern-5mu_{\underline{k}} \succeq \TT\mkern-5mu_\infty$}).}
We proceed to prove $\TT_k \succeq \TT_\infty$ for all $k = 0, \ldots , \underline{k}$ by induction over $k$.
For $k = 0$, this is trivially true. 
Let $k \in \set{0, \ldots , \underline{k}-1}$ and suppose that $\TT_k \succeq \TT_\infty$. 
For every triangle $T \in \MM_k$, there exists $T' \in \TT_\infty$ such that $T' \subsetneqq T$ and $\gen{T'} > \gen{T} + K_\ell$. 
Hence, by the induction hypothesis, $\TT_\infty$ is a refinement of $\TT_k$ such that all elements in $\MM_k$ have been refined. 
Since $\TT_{k + 1} = \mathtt{NVB}(\TT_k,\MM_k)$ is the coarsest refinement of $\TT_k$ such that all elements in $\MM_k$ are refined (see, e.g.,~\cite[Theorem~5.1]{Stevenson2008}), we obtain that $\TT_{k+1} \succeq \TT_\infty$. 
In particular, we showed that $\TT_{\underline{k}} \succeq \TT_\infty$.

\textbf{Step~2 (\boldsymbol{$\TT\mkern-5mu_{\underline{k}} \text{ and } \TT\mkern-5mu_\infty \text{ are } K_\ell$}-compatible).}
By Step~1, $\TT_{\underline{k}} \succeq \TT_\infty$. 
It remains to show that for all $T \in \TT_{\underline{k}}$ and $T' \in \TT_\infty$ with $|T \cap T'| > 0$ there holds $\gen{T'} \leq \gen{T} + K_\ell$. 
Since $\MM_{\underline{k}} = \emptyset$ for the set in~\eqref{eq:compmarked}, the claim follows.

\textbf{Step~3 (Definition of \boldsymbol{$\widehat{\TT}$}).}
We recall the $\mathtt{bisec3}$-refinement procedure, where each edge is bisected exactly once, i.e., every triangle is partitioned into four triangles of the same size. 
Hence, given a triangulation $\TT$, there holds $\abs{\mathtt{bisec3}(\TT)} = 4\abs{\TT}$ and $\gen{T'} = \gen{T} + 2$ for all $T' \in \mathtt{bisec3}(\TT)$ and $T \in \TT$ with $T' \subsetneqq T$. 
For $n \in \N$, we set 
\begin{equation*}
 (\mathtt{bisec3})^n(\TT) 
 \coloneqq \underbrace{\mathtt{bisec3}(\mathtt{bisec3}(\ldots \mathtt{bisec3}(\TT) \ldots))}_{n \text{ times}}.
\end{equation*}

We define the finest common coarsening of all triangulations $\TT_{\ell,\nu}$ with $\nu \in \RRR$ to be
\begin{equation*}
 \TT' 
 \coloneqq  \Big\{ \bigcup_{\nu \in \RRR} T_\nu \ | \ T_\nu \in \TT_{\ell,\nu} \text{ with } \nu \in \RRR \text{ and } \Big|\bigcap_{\nu \in \RRR} T_\nu\Big| > 0\Big \}.
\end{equation*}
By analogous arguments as for the overlay in~\cite[Lemma~3.7]{Cascon2008}, we obtain that $\TT'$ is a conforming triangulation, which is coarser than $\TT_{\ell,\nu}$ for all $\nu \in \RRR$. 
We set
\begin{equation*}
 \widehat{\TT} 
 \coloneqq (\mathtt{bisec3})^{K_\ell/2}(\TT').
\end{equation*}

\textbf{Step~4 (Proof of \boldsymbol{$\TT\mkern-5mu_{\underline{k}} \succeq \widehat{\TT}$}).}
We proceed to prove $\TT_k \succeq \widehat{\TT}$ for all $k = 0, \ldots , \underline{k}$ by induction over $k$. 

\textbf{Base case (\boldsymbol{$k = 0$}).} For $k = 0$, this is trivially true.

\textbf{Induction hypothesis (on \boldsymbol{$k$}).}
Let $k \in \set{0, \ldots , \underline{k}-1}$ and suppose that $\TT_k \succeq \widehat{\TT}$. 

\textbf{Induction step (on \boldsymbol{$k$}).}
To show that $\TT_{k + 1} = \mathtt{NVB}(\TT_k, \MM_k) \succeq \widehat{\TT}$, due to~\cite[Theorem~5.1]{Stevenson2008}, it suffices to show that $\TT_k^M \coloneqq \mathtt{NVB}(\TT_k,\set{M}) \succeq \widehat{\TT}$ for all $M \in \MM_k$. 
To this end, let $M \in \MM_k$, $T \in \TT_k^M$, and $\widehat{T} \in \widehat{\TT}$ with $\abs{T \cap \widehat{T}} > 0$. 
If $T \in \TT_k \cap \TT_k^M$, the induction hypothesis yields that $\widehat{T} \subseteq T$. 
It remains to consider the case $T \in \TT_k^M \setminus \TT_k$. 
Since $\TT_k^M = \mathtt{NVB}(\TT_k,\set{M})$, there exist $m \in \N$ and $T_1, \ldots ,T_m \in \TT_k$ with $T_1 = M$, $T_{j + 1} \cap T_j = E_{T_j}$ for all $j = 1, \ldots , m - 1$, $E_{T_j} \neq E_{T_{j + 1}}$ for all $j = 1, \ldots , m-2$, and $\TT_k \setminus \TT_k^M = \set{T_1, \ldots , T_m}$. 
Since $T \in \TT_{k}^M \setminus \TT_{k}$, there exists $j \in \set{1, \ldots , m}$ such that $T \subsetneqq T_j$. 

We proceed by induction over $j$ to show that, for all $T \in \TT_k^M$ and $\widehat{T} \in \widehat{\TT}$ with $T \subsetneqq T_j$ and $\abs{T \cap \widehat{T}}>0$, there holds $\widehat{T} \subseteq T$. 

\textbf{Base case (\boldsymbol{$j = 1$}).} For $j = 1$, we obtain $T \subsetneqq T_1 = M$. 
Due to~\cite[Theorem~5.1]{Stevenson2008}, we have $\gen{M} + 1 = \gen{T}$.  
In order to prove $\widehat{T} \subseteq T$, we argue by contradiction and assume that $T \subsetneqq \widehat{T}$. 
Since $M \in \MM_k$, there exists $N \in \TT_\infty$ such that $N \subsetneqq M$ with $\gen{N} > \gen{M} + K_\ell$. 
Since $M \cap \widehat{T} \supseteq T$ and because of the induction hypothesis (on $k$), i.e., $\TT_k \succeq \widehat{\TT}$, we obtain $T \subsetneqq \widehat{T} \subseteq M$. 
Since $T$ is generated by bisecting $M$ exactly once, this implies $M = \widehat{T}$. 
Because $M = \widehat{T} \in (\mathtt{bisec3})^{K_\ell/2}(\TT')$, there exists $T' \in \TT'$ with $M \subsetneqq T'$, $\gen{T'} = \gen{M} - K_\ell$ and $T \subsetneqq T'$. 
Since $(\TT_{\ell,\nu})_{\nu \in \RRR}$ are pairwise $2K_\ell$-compatible, the meshes $\TT'$ and $\TT_\infty$ are $2K_\ell$-compatible. 
Since $N \subsetneqq M = \widehat{T} \subsetneqq T'$, there holds
\begin{equation*}
 \gen{N}
 \leq \gen{T'} + 2K_\ell
 = \gen{\widehat{T}} + K_\ell
 = \gen{M} + K_\ell,
\end{equation*}
which contradicts the fact that $\gen{N} > \gen{M} + K_\ell$. 
Hence, $\widehat{T} \subseteq T$, which concludes the proof for $j = 1$.

\textbf{Induction hypothesis (on \boldsymbol{$j$}).} Let $j \in \set{1, \ldots , m-1}$ and suppose that for all $T \in \TT_k^M$ and $\widehat{T} \in \widehat{\TT}$ with $T \subsetneqq T_j$ and $\abs{T \cap \widehat{T}}>0$ there holds $\widehat{T} \subseteq T$. 

\textbf{Induction step (on \boldsymbol{$j$}).}
Let $T \in \TT_k^M$ and $\widehat{T} \in \widehat{\TT}$ with $T \subsetneqq T_{j + 1}$ and $\abs{T \cap \widehat{T}}>0$. 
We argue by contradiction and assume $T \subsetneqq \widehat{T}$.
The induction hypothesis on $k$, i.e., $\TT_k \succeq \widehat{\TT}$, and $\widehat{T} \cap T_{j + 1} \supseteq T$ establish that $\widehat{T} \subseteq T_{j + 1}$. 
We first consider the case $\widehat{T} = T_{j + 1}$. 
This yields $\widehat{T} = T_{j + 1} \in \widehat{\TT} \cap \TT_k$. Since $E_{T_j} = T_{j + 1} \cap T_j$, there exists $\widehat{N} \in \widehat{\TT}$ with $\widehat{N} \cap \widehat{T} = E_{T_j}$, which implies $T_j \subseteq \widehat{N}$. 
However, by the induction hypothesis on $j$, $\widehat{N}$ is contained in a child $\widetilde{T} \subsetneqq T_j$. 
The contradiction $\widehat{N} \subseteq \widetilde{T} \subsetneqq T_j \subseteq \widehat{N}$ proves that $\widehat{T} \subseteq T$.

We now consider the situation when $\widehat{T} \neq T_{j + 1}$. 
Since $T \subsetneqq \widehat{T} \subsetneqq T_{j + 1}$, we obtain that $\gen{T} = \gen{\widehat{T}} + 1 = \gen{T_{j + 1}} + 2$. 

By~\cite[Theorem~5.1]{Stevenson2008}, $T_{j + 1}$ is bisected exactly twice and there exists $T' \in \TT_k^M$ with $T' \subsetneqq T_j$ and $T' \cap T$ is an edge in $\TT_k^M$. 
This, $T \subsetneqq \widehat{T} \subsetneqq T_{j + 1}$, and $T_{j + 1} \cap T_j = E_{T_j} \neq E_{T_{j + 1}}$ yields $\widehat{T} \cap T_j = E_{T_j} = E_{\widehat{T}}$. 
Hence, there exists $\widehat{N} \in \widehat{\TT}$ with $\widehat{T} \cap \widehat{N} = E_{T_j}$, which implies $T_j \subseteq \widehat{N}$. 
Similar to the previous case, the induction hypothesis yields that $\widehat{N}$ is contained in a child $\widetilde{T} \subsetneqq T_j$. 
This contradiction proves that $\widehat{T} \subseteq T$.

\textbf{Step~5.}
Finally, we prove that $\TT_{\underline{k}}$ and $\TT_{\ell,\nu}$ are $K_\ell$-compatible for all $\nu \in \RRR$.
Let $\nu \in \RRR$, $T \in \TT_{\underline{k}}$, and $T_\nu \in \TT_{\ell,\nu}$ with $|T \cap T_\nu| > 0$. 
By Step~1, $\TT_{\underline{k}} \succeq \TT_\infty$, which implies the existence of $T_\infty \in \TT_\infty$ with $T_\infty \subseteq T$, $\gen{T} \leq \gen{T_\infty}$, and $|T_\nu \cap T_\infty| > 0$. 
Since $\TT_\infty = \bigoplus_{\mu \in \RRR} \TT_{\ell,\mu} \preceq \TT_{\ell,\nu}$, there holds $T_\infty \subseteq T_\nu$ and $\gen{T_\nu} \leq \gen{T_\infty}$. 
By Step~2, $\TT_{\underline{k}}$ and $\TT_\infty$ are $K_\ell$-compatible, which yields
\begin{equation} \label{eq:meshcoreq1}
 \gen{T_\nu} 
 \leq \gen{T_\infty} 
 \leq \gen{T} + K_\ell.	
\end{equation}
By Step~4, $\TT_{\underline{k}}$ is coarser than $\widehat{\TT}$. 
Hence, there exists $\widehat{T} \in \widehat{\TT}$ with $\widehat{T} \subseteq T$, $\gen{T} \leq \gen{\widehat{T}}$, and $\abs{\widehat{T} \cap T_\nu} > 0$. 
Since $\widehat{\TT} = (\mathtt{bisec3})^{K_\ell/2}(\TT')$, there exists $T' \in \TT'$ with $\widehat{T} \subsetneqq T'$ and $\gen{T'} = \gen{\widehat{T}} - K_\ell$. 
There holds $\abs{T' \cap T_\nu} \geq \abs{\widehat{T} \cap T_\nu} > 0$. 
Since $\TT'$ is coarser than $\TT_{\ell,\nu}$, we obtain $T_\nu \subseteq T'$ and $\gen{T'} \leq \gen{T_\nu}$. 
With $\gen{T} \leq \gen{\widehat{T}} = \gen{T'} + K_\ell$ and $\gen{T'} \leq \gen{T_\nu}$,
\begin{equation*}
 \gen{T}
 \leq \gen{\widehat{T}}
 = \gen{T'} + K_\ell
 \leq \gen{T_\nu} + K_\ell.
\end{equation*}
Overall, this and~\eqref{eq:meshcoreq1} conclude the proof of $\abs{\gen{T_\nu} - \gen{T}} \leq K_\ell$.
\end{proof}

\end{document}